\documentclass{amsart}

\usepackage{amsmath}
\usepackage{amssymb}
\usepackage{enumerate}
\usepackage{graphicx}
\usepackage{amscd}

\newtheorem{theorem}{Theorem}[section]
\newtheorem{cor}[theorem]{Corollary}
\newtheorem{lem}[theorem]{Lemma}
\newtheorem{prop}[theorem]{Proposition}

\newtheorem{rem}[theorem]{Remark}
\numberwithin{equation}{section}
\newcommand{\Complex}{\mathbb{C}}
\newcommand{\m}{\mathfrak{E}}
\newcommand{\ep}{\epsilon}
\newcommand{\lam}{\lambda}
\newcommand{\gam}{\gamma}
\newcommand{\al}{\alpha}
\newcommand{\la}{\langle}
\newcommand{\ra}{\rangle}
\newcommand{\f}{\mathbb{F}}
\newcommand{\C}{\mathcal{C}}
\newcommand{\Z}{\mathbb{Z}}

\renewcommand{\Gamma}{\varGamma}
\renewcommand{\epsilon}{\varepsilon}
\renewcommand{\bar}{\overline}
\renewcommand{\hat}{\widehat}
\renewcommand{\leq}{\leqslant}
\renewcommand{\geq}{\geqslant}

\begin{document}


\author{Hung P. Tong-Viet}
\email{Tong-Viet@ukzn.ac.za}
\address{School of Mathematical Sciences,
University of KwaZulu-Natal, Pietermaritzburg 3209, South Africa}

\title[Rank $3$ permutation characters]{Rank $3$ permutation characters and maximal subgroups}

\begin{abstract}
Let $G$ be a transitive permutation group acting on a finite set
$\m$ and let $P$ be the stabilizer in $G$ of a point in $\m.$ We say
that $G$ is primitive rank $3$ on $\m$ if $P$ is maximal in $G$ and
$P$ has exactly three orbits on $\m.$ For any subgroup $H$ of $G,$
we denote by $1_H^G$ the permutation character (or permutation
module) over $\Complex$ of $G$ on the cosets $G/H.$ Let $H$ and $K$
be subgroups of $G.$ We say $1_H^G\leq 1_K^G$ if $1_K^G-1_H^G$ is a
character of $G.$ Also a finite group $G$ is called nearly simple
primitive rank $3$ on $\m$ if there exists a quasi-simple group $L$
such that $L/Z(L)\unlhd G/Z(L)\leq Aut(L/Z(L))$ and $G$ acts as a
primitive rank $3$ permutation group on the cosets of some subgroup
of $L.$ In this paper we classify all maximal subgroups $M$ of a
nearly simple primitive rank $3$ group $G$ of type
$L=\Omega_{2m+1}(3), m\geq 3,$ acting on an $L$-orbit $\m$ of
non-singular points of the natural module for $L$ such that
$1_P^G\leq 1_M^G,$ where $P$ is the stabilizer of a non-singular
point in $\m.$ This result has an application to the study of
minimal genera of algebraic curves which admit group actions.
\end{abstract}
\subjclass{Primary 20B15; secondary 20C20}

\keywords{permutation characters, nearly simple groups, maximal subgroups}
\date{\today}
\maketitle

\section{Introduction} Let $G$ be a transitive permutation group
acting on a finite set $\m$ and let $P$ be the stabilizer of a point
in $\m.$ We say that $G$ is \emph{primitive} on $\m$ if and only if
$P$ is maximal in $G.$ We define the \emph{rank} of $G$ on $\m$ to
be the number of $P$-orbits on $\m.$ For any subgroup $H$ of $G,$ we
denote by $1_H^G$  the permutation character of $G$ on the cosets
$G/H$ over $\Complex.$   We also use the same notation $1_H^G$ for
the permutation module. Let $H, K$ be subgroups of $G.$ Consider the
permutation characters $1_H^G$ and $1_K^G,$ we say  $1_H^G\leq
1_K^G$ if $1_K^G-1_H^G$ is a character of $G.$ In terms of
permutation modules, $1_H^G\leq 1_K^G$ if  and only if $1_H^G$ is
isomorphic to a submodule of $1_K^G.$ A finite group $L$ is said to
be \emph{quasi-simple} if $L$ is perfect and $L/Z(L)$ is simple. A
finite group $G$ is called \emph{nearly simple  of type $L$} if
$L\unlhd G$ and $L/Z(L)\leq G/Z(L)\leq Aut(L/Z(L))$ for some
quasi-simple group $L.$  Moreover a finite group $G$ is called
\emph{almost simple of type $L$} or {\em almost simple with socle
$L$} if $L\unlhd G\leq Aut(L)$ for some finite simple group $L.$ It
follows from definitions that if $G$ is nearly simple of type $L$
then $G/Z(L)$ is almost simple of type $L/Z(L).$ Assume that $m\geq
3$ is an integer. Let $L$ be one of the following quasi-simple
groups $\Omega_{2m+1}(3),\Omega_{2m}^\ep(3),\Omega_{2m}^\ep(2)$ or
$SU_m(2).$ Let $G$ be a nearly simple group of type $L$ such that
$G$ acts on the $L$-orbit $\m(V)$ of non-singular points in the
natural module $V$ for $L$ and let $P$ be the stabilizer of a
non-singular point in $\m(V).$ Then $G$ is a primitive rank $3$
group on $\m(V).$  In this situation, we say that $G$ is a
\emph{nearly simple primitive rank $3$ group of type $L.$} In this
paper, we shall classify all maximal subgroups $M$ of a nearly
simple primitive rank $3$ group $G$ of type $\Omega_{2m+1}(3)$ such
that $1_P^G\leq 1_M^G.$  The remaining types of nearly simple
primitive rank $3$ groups of type $L$ have been considered in
\cite{TV}.

\begin{table}
\centering
\caption{ $M\in \mathcal{C}$}\label{in2}
\begin{tabular}{l|c|c|c|r }
\hline $L$& type of $M$ &conditions&orbits&Ref\\\hline
$\Omega_{2m+1}(3)$&$O_1(3)\perp O_{2m}^\ep(3)$&&&\ref{odd1}\\
&$P_\al$& $1\leq \al \leq m$& $\leq 2$&\ref{odd2} \\
& $O_1(3)\wr S_{2m+1}$ &$(2m+1,\xi,r)=(5,\pm,t)$ &$2$&\ref{odd3}\\
&&$(7,+,t)$&$2$&\\
& $O_{\al}(3^3)$& $2m+1=3\al$ &$3$&\ref{odd5}\\\hline
\end{tabular}
\end{table}

\begin{table}
\centering
\caption{ $M\in \mathcal{S}$}\label{in3}
\begin{tabular}{l|c|c|c|r}
\hline
$L$&socle of $M$ & modules  &orbits &Ref\\\hline
$\Omega_{7}(3)$&$A_9$&$\lam=(8,1)$&$\leq 2$&\ref{odd10}\\
$\Omega_7(3)$&$PSp_6(2)$&&$2$&\ref{odd11}\\
$\Omega_{13}(3)$&$PSp_6(3)$&$\lam_2$&$\leq 2$&\ref{odd14}\\
$\Omega_{7}(3)$&$G_2(3)$&$\lam_1,\lam_2$&$1$&\ref{odd16}\\
$\Omega_{25}(3)$&$F_4(3)$&$\lam_4$&$\leq 2$&\\\hline
\end{tabular}
\end{table}

\begin{table}
\centering
\caption{Possible exceptions}\label{in4}
\begin{tabular}{l|c|c|r }
\hline
$L$&socle of $M$ & modules  & Ref\\\hline
$\Omega_{41}(3)$&$S_8(3)$&$\lam_1$&\ref{odd14}\\
$\Omega_{77}(3)$&$E_6(3)$&adjoint module&\ref{odd16}\\
$\Omega_{133}(3)$&$E_7(3)$&adjoint module&\\\hline
\end{tabular}
\end{table}

\begin{theorem}\label{in1}
Let $L=\Omega_{2m+1}(3),$ with $m\geq 3,$ and let $G$ be a nearly
simple primitive rank $3$ group of type $L.$ Let $P$ be the
stabilizer of a non-singular point in $V.$ Let $M$ be any maximal
subgroup of $G.$  Then $1_P^G\leq 1_M^G$ unless the pairs $(L,M)$
appear in Tables \ref{in2}-\ref{in4}.
\end{theorem}
\begin{rem}
For all pairs $(L,M)$ in Tables \ref{in2} and \ref{in3}, $1_P^G$ is
not contained in $1_M^G.$ The pairs in Table \ref{in4} are the cases
that we have not determined whether or not there is a containment.
\end{rem}
The main motivation for this work comes from algebraic curves which
admit group actions. From Riemann's Existence Theorem, we know that
for every finite group there are infinitely many Riemann surfaces
with automorphism group $G.$ We would like to identify $G$-curves of
smallest possible genus. Let $X$ be a compact Riemann surface with
$G\cong Aut(X).$ Let $\ell$ be a prime such that $(\ell,|G|)=1.$ Let
$J_X$ be the Jacobian of $X$ and let $J_X[\ell]$ be the
$\ell$-torsion points of $J_X.$ It is well known that $J_X[\ell]$ is
a characteristic $\ell$ module and the genus of the quotient curve
$X/P$ can be computed using the relation $2g(X/P)=dim(J_X[\ell]^P),$
where $P\leq G$ and $J_X[\ell]^P$ is the fixed point space of
$J_X[\ell]$ under $P$ (see \cite[Lemma $8.1$]{FGM}). Now combining
the latter fact with the Frobenius reciprocity, we deduce that if
$P,M$ are proper subgroups of $G$ such that $1_P^G\leq 1_M^G,$ then
$g(X/P)\leq g(X/M).$ For a proof of this result see \cite[Corollary
$8.2$]{FGM}. Theorem \ref{in1} with the discussion above now yield
the following.

\begin{cor}\label{in5} Assume the assumption and
notation of Theorem \ref{in1}. Let $X$ be a compact Riemann surface with $G\leq Aut(X),$ and
inertia groups $\la g_1\ra, \dots,\la g_r\ra$ over $X/G.$ Then $g(X/P)\leq g(X/M)$ for any maximal
subgroup $M$ of $G$ which does not appear in Tables \ref{in2}-\ref{in4}.
\end{cor}

If $G$  is doubly transitive on a non-empty set $\m$ with point stabilizer $P,$ then either
$1_P^G\leq 1_M^G$ or $G=PM$ for any maximal subgroup $M$ of $G.$ The maximal factorization of
almost simple groups was classified completely by M. Liebeck, C. Praeger and J. Saxl in \cite{lps},
so for almost simple groups, we can tell exactly for which $M$ we have the containment $1_P^G\leq
1_M^G.$  In case of rank $3,$ M. Aschbacher, R. Guralnick and K. Magaard \cite{asc} have a
criterion in terms of the Higman rank $3$ parameters. In that paper, they consider the case when
$G$ is a nearly simple classical group acting on the set of singular points on its natural module.
Also a partial result of this case has been dealt with by D. Frohardt, the second and the third
authors above in \cite{FGM}. In the case when $G$ is a nearly simple primitive rank $3$ group of
sporadic type, the containment of the permutation characters of $G$ is completely determined since
all the permutation characters of maximal subgroups of $G$ are stored in \cite{GAP}, except $HS:2,
Fi_{22}:2$ and $Fi_{24}':2.$

We now describe our strategy. Let $L$ be a finite simple classical
group of degree $d\geq 2,$ defined over a finite field ${\f}_q, q$ a
prime power, and let $V$ be the natural module for $L.$ Assume that
$G$ is an almost simple group with simple socle $L.$ We have a
powerful theorem on the subgroup structure of $G$ by M. Aschbacher.
The theorem says that if $M$ is a subgroup of $G$ then $M$
belongs to a collection $\mathcal{C}(G)$ of geometric subgroups of
$G$ or $M\in \mathcal{S}(G),$ that is, $M$ is an almost simple group
and the full covering group of the socle of $M$ acts absolutely
irreducible on the natural module $V$ for $G$ and cannot be realized
over any proper subfield. Thus if $M$ is a maximal subgroup of $G$
then either $M\in \mathcal{C}(G)$ or $M\in \mathcal{S}(G).$ The
subgroup structure and the maximality among members of
$\mathcal{C}(G)$ have been determined by P. Kleidman and M. Liebeck
in \cite{kl} when the degree is at least $13.$ For this case, using
the geometrical properties of the groups, we can solve the problem
completely. When $M$ is not a geometric subgroup, that is, $M\in
\mathcal{S}(G),$ the problem is much more complicated as we still do
not know which members of $\mathcal{S}(G)$ are maximal. Now assume
that $M\in \mathcal{S}(G).$ Denote by $S$ the socle of $M.$ So $S$
is a non-abelian finite simple group. According to the
Classification of Finite Simple Groups, $S$ is an alternating group
of degree at least $5,$ a finite group of Lie type or one of the
$26$ sporadic groups. By way of contradiction, we assume that
$1_P^G\not\leq 1_M^G.$ From this assumption, we will get an upper
bound for the dimension of $V$ in terms of the size of the
automorphism group of $S.$ From the  definition of members in
$\mathcal{S}(G),$ the full covering group $\hat{S}$ of $S$ acts
absolutely irreducible on $V.$ Now using the information on the
lower bound for the dimension of the absolutely irreducible
representations of finite simple groups, we will get a finite list
of cases that we can handle either by constructing the
representations or by computer program GAP \cite{GAP}.

As in the almost simple doubly transitive case, we can get a list of
maximal subgroups $M$ such that $1_P^G\not\leq 1_M^G.$ In Table
\ref{in2}, we list all the cases when $M\in \mathcal{C}(G),$ Table
\ref{in3} contains all the cases when $M\in \mathcal{S}(C),$ and in the
last table, we list the cases that we have not determined whether or
not   there is a containment. Notice that we only have a finite
number of exceptions in Table \ref{in4}. Also there is a finite
number of cases in Table \ref{in3}.

For the notation in the tables of Theorem \ref{in1}, the columns `orbits' give the number of orbits
of $M$ on $\m(V)$ and this is also the number of double cosets of $G$ on $P$ and $M.$ The first
columns are the type of the nearly simple group $G.$ The last columns `Ref' give the references for
the result. For example `$5.3$' means that the case is dealt with in Proposition $5.3.$

\section{Preliminaries}\label{pre}
We adopt the constructions and notation of \cite{kl}. Fix a finite field ${\f}_q,$ $q$ a prime
power. Let $V$ be an ${\f}_q$ vector space of dimension $n.$ The {\em  general linear group}
$GL(V)$ of $V$ over ${\f}_q$ is a group of all non-singular ${\f}_q$-linear transformations of $V.$
The {\em special linear group} of $V$ over ${\f}_q,$ $SL(V),$ is the group of all elements of
$GL(V)$ with determinant $1.$ The \emph{projective linear groups} $PGL(V)$ and $PSL(V)$ are
obtained  by factoring out the scalar matrices in the corresponding linear groups. For any subgroup
$X$ of $GL(V),$ we write $PX$ or $\bar{X}$ for the corresponding projective group $X/X\cap
{\f}_q^*.$

A map $g$ from $V$ to itself is called an ${\f}_q$-{\em semilinear
transformation} of $V$ if there is a field automorphism
$\sigma(g)\in Aut({\f}_q)$ such that for all $v,w\in V$ and
$\lambda\in {\f}_q,$
\begin{equation}
\label{pre1}
(v+w)g=vg+wg\mbox{ and } (\lambda v)g=\lambda^{\sigma(g)}(vg)
\end{equation}

The \emph{general semilinear group } of $V$ over ${\f}_q,$ $\Gamma
L(V)$ consists  of all non-singular ${\f}_q$-semilinear
transformations of $V.$ As ${\f}_q^*\trianglelefteq\Gamma L(V),$ we
can factor out the scalars to get the {\em projective general
semilinear group } $P\Gamma L(V).$ Let $\kappa$ be a left linear or
a quadratic form on $V.$ Observe $\kappa$ is a map from $V^k$ to
${\f}_q$ where $k=1,2.$ Define
$$I(V,{\f}_q,\kappa)=\{g\in GL(V)\;|\; \kappa({\bf v}g)=\kappa({\bf v})\mbox{ for all ${\bf v}\in V^k$} \};$$
$$S(V,{\f}_q,\kappa)=I(V,{\f}_q,\kappa)\cap SL(V);$$
$$\Xi(V,{\f}_q,\kappa)=\{g\in \Gamma L(V)\;|\; \kappa({\bf v}g)=\tau(g)\kappa({\bf v})^{\sigma(g)}\mbox{ for all ${\bf v}\in V^k$}\}$$
$\mbox{where $\tau(g)\in {\f}_q^*$}, \sigma(g)\in Aut({\f}_q),$ and
$$\Lambda(V,{\f}_q,\kappa)=\{g\in \Xi(V,{\f}_q,\kappa)\;|\;\sigma(g)=1 \}.$$
Define $$A=\left \{\begin{array}{lll}
                              \Xi\langle \iota\rangle& \mbox{in case
                              $\mathbf{L}$ with $n\geq 3;$}\\
                              \Xi& \mbox{otherwise;}
                              \end{array}\right.$$ and
$$\Omega=\left \{\begin{array}{lll}
                              \mbox{certain subgroup of index 2 in S} & \mbox{in case
                            $\mathbf{O};$ }\\
                             S & \mbox{otherwise;}
                             \end{array}\right.$$
where $\iota$ is an inverse-transpose automorphism of $GL(V,{\f}_q).$ We get a sequence of groups:
$\Omega\leq S\leq I\leq \Lambda\leq \Xi\leq A.$ Note that in \cite{kl}, $\Xi(V,{\f}_q,\kappa),$ and
$\Lambda(V,{\f}_q,\kappa)$ are denoted by $\Gamma(V,{\f}_q,\kappa),$ and $\Delta(V,{\f}_q,\kappa),$
respectively. For more details, see \cite[Chapter $2$]{kl}.

Let $(V,{\f}_q,Q)$ be a classical orthogonal geometry  with $q$ odd,
and $dimV=n.$ Let $\ep=sgn(Q)$ be the sign of the quadratic
form $Q.$ Note that $\ep=\circ$ when $n$ is odd, otherwise,
$\ep$ is either $+$ or $-.$ A square or a non-square in
${\f}_q^*$ will be denoted by $\square$ or $\boxtimes,$
respectively. If $W$ is a non-degenerate subspace of $V$ then we
write $sgn(W)=sgn(Q_W),$ where $Q_W$ is the restriction of $Q$ to
$W.$ When $n$ is odd, for any non-zero vector $x$ in $V,$ we denote
by $S(n,x)$ the number of all vectors $v\in V$ with $Q(v)=Q(x).$
When $n$ is even, for any $\gam\in {\f}_q,$ we denote by
$S^\ep(n,\gam)$ the number of all vectors $v\in V$ with $Q(v)=\gam.$
For $x\in V\setminus \{0\},$ a one-space with representative $x$
will be called a {\em point} in $V$ and denoted by $\la x\ra.$ We
now define a type function $\rho=\rho_V$ on $V\setminus \{0\}$ as
follows: if $x$ is a singular vector in $V,$ that is, $x\in
V\setminus \{0\}$ and $Q(x)=0,$ then $\rho(x)=0.$ If $dimV$ is even,
then $\rho(x)=Q(x).$ If $dimV$ is odd, then
$\rho(x)=\mbox{sgn}(x^\perp).$ Assume that $dimV$ is odd. Let $x$ be
a non-singular vector in $V.$ We say $x$ is a {\em plus vector} if
$\rho(x)=+;$ and $x$ is a {\em minus vector} if $\rho(x)=-.$ We also
say that $\la x\ra$ is of plus or minus type according to whether
its representative $x$ is a plus or a minus vector. Let $x\in V$ be a
non-singular vector with $\rho(x)=\xi.$ Define $\m_{\xi}^{\ep}(V)$
to be the set of all non-singular points of type $\xi$ in $V,$ where
$\ep=sgn(Q).$ When $\ep=\circ,$ we will write $\m_\xi(V)$ instead of
$\m_\xi^\circ(V). $

\begin{lem}\label{pre9} Let $(V,{\f}_q,Q)$ be a classical orthogonal
geometry with $dimV$ odd. Two non-singular vectors $x, y$ $($two
non-singular points $\la x\ra, \la y\ra)$ have the same type if and
only if $Q(x)\equiv Q(y)$ $(\mbox{mod $({\f}_q^*)^2$}).$
Hence for any non-singular vector $z,$ we have $\m_{\rho(z)}(V)=
\{\la v\ra\subseteq V\;|\; Q(v)\equiv Q(z)\;  (\mbox{mod $({\f}_q^*)^2$}) \}.$
 In particular if $q=3,$ then
$\m_{\rho(z)}(V)=\{\la v\ra \subseteq V\;|\;Q(v)=Q(z)\}.$
\end{lem}

\begin{proof}
Assume  that $Q(x)\equiv Q(y)\;(\mbox{mod $({\f}_q^*)^2$}).$ By \cite[Proposition $2.5.4(ii)$]{kl},
$\langle x\rangle$ and $\langle y\rangle$ are isometric. By Witt's lemma, this isometry extends to
an isometry $g$ of $V$ such that $\la x\ra g=\la y\ra.$ As $\la x\ra, \la y\ra$ are non-degenerate,
$x^\perp g=y^\perp.$ It follows that $x^\perp$ and $y^\perp$ are isometric, and hence
$sgn(x^\perp)=sgn(y^\perp),$ so that $\rho(x)=\rho(y).$ Now assume that $x, y$ have the same type.
By Witt's lemma and \cite[Proposition $2.5.4(i)$]{kl}, there exists an isometry between $x^{\perp}$
and $y^{\perp}.$ This isometry can extend to an isometry $g$ of $V$ such that
$(x^{\perp})g=y^{\perp}.$ Since $(x^{\perp})^{\perp}=\langle x\rangle,$ and
$(y^{\perp})^{\perp}=\langle y\rangle,$ $(\langle x\rangle)g=\langle y\rangle .$ Thus $xg=\mu y$
for some $\mu\in {\f}_q^*.$ Therefore $Q(x)=Q(xg)=Q(\mu y)=\mu^2Q(y).$ The other statements are
obvious.
 \end{proof}

The following lemma will be used to compute the Higman rank $3$
parameters for orthogonal groups.

\begin{lem}\label{pre10} Let $(V,{\f}_q,Q)$ be a classical orthogonal geometry  with $q$ odd and $\ep=sgn(Q).$

$(1)$ if $dimV=2k$ and $\gam\in {\f}_q,$ then $S^\ep(2k,\gam)=\left
\{\begin{array}{lcr}
                              q^{2k-1}+\ep(q^k-q^{k-1})& \mbox{if
                              $\gam=0,$}\\
                              q^{2k-1}-\ep q^{k-1} & \mbox{if $\gam\neq 0;$}
                              \end{array}\right.$

$(2)$ if $dimV=2k+1$ and $x\in V-\{0\}$  then
$S(2k+1,x)=q^{2k}+\rho(x) q^k.$
\end{lem}

\begin{proof}
Statement $(1)$ follows from \cite[Proposition $9.10$]{Grov}. For $(2),$ let $\gam=Q(x)$ and
$\xi=\rho(x).$ Assume first that $x$ is a non-singular vector. Then $V=\la x\ra \perp x^\perp,
\mbox{sgn}(x^\perp)=\xi$ and $\mbox{dim}(x^\perp)=2k.$ For any vector $v\in V$ with $Q(v)=\gam,$
write $v=\varphi x+v_0,$ where $\varphi\in {\f}_q$ and $v_0\in x^\perp.$ We have
$Q(v_0)=Q(v)-\varphi^2Q(x)=\gam(1-\varphi^2).$ If $\varphi=\pm 1,$ then $Q(v_0)=0,$ hence by $(1),$
there are $2S^\xi(2k,0)=2(q^{2k-1}+\xi(q^k-q^{k-1}))$ such $v.$ If $\varphi\not=\pm 1,$ then
$Q(v_0)=\gam(1-\varphi^2)\not=0,$ hence by $(1),$ again, there are
$(q-2)S^\xi(2k,\gam(1-\varphi^2))=(q-2)(q^{2k-1}-\xi  q^{k-1}))$ such $v.$ Thus
$S(2k+1,x)=2S^\xi(2k,0)+(q-2)S^\xi(2k,\gam(1-\varphi^2))=q^{2k}+\xi q^k.$ Assume that $x$ is a
singular vector. Observe that $V$ always contains a non-singular vector $y.$ Let $\eta=\rho(y)$ and
$\mu=Q(y).$ We have $V=\la y\ra \perp y^\perp,$ $sgn(y^\perp)=\eta, Q(y)=\mu \in {\f}_q^*$ and
dim$(y^\perp)=2k.$ For any $v\in V$ with $Q(v)=0,$ write $v=\varphi y+v_0,$ where $\varphi\in
{\f}_q$ and $v_0\in y^\perp.$ Then $Q(v_0)=Q(v)-\varphi^2Q(y)=-\mu\varphi^2.$ If $\varphi=0,$ then
$Q(v_0)=0,$ hence by $(1),$ there are $S^\eta(2k,0)=q^{2k-1}+\eta  (q^k-q^{k-1})$ such $v.$ If
$\varphi\not=0,$ then $Q(v_0)=-\mu\varphi^2\not=0,$ hence by $(1),$ again, there are $(q-1)S^\eta
(2k,-\mu\varphi^2)=(q-1)(q^{2k-1}-\eta q^{k-1}))$ such $v.$ Thus
$S(2k+1,x)=S^\eta(2k,0)+(q-1)S^\eta(2k,-\mu\varphi^2)=q^{2k}.$ The proof is complete.
\end{proof}

By the classification of the finite simple groups, every non-abelian
finite simple group is either an alternating group $A_n, n\geq 5,$ a
finite simple group of Lie type, or one of the $26$ sporadic simple
groups. The next lemma gives an upper bound for the full
automorphism groups of the non-abelian finite simple groups of Lie
type. The proof will be omitted.

\begin{lem}\label{pre17}
If $L$ is a finite simple group of Lie type, then $|Aut(L)|\leq
f(L),$ where $f(L)$ is given in Table \ref{pre18}.
\end{lem}
\begin{table}
\centering
 \caption{Upper bounds for the size of $Aut(L)$.}
 \begin{tabular}{l|c|c|r}\hline
$L$ &$f(L)$&$L$&$f(L)$\\\hline
$L_n(q), n\geq 3$&$q^{n^2}$&$E_8(q)$&$q^{249}$\\
$PSp_{2n}(q)$&$q^{2n^2+n+1}$&$F_4(q)$&$q^{53} $ \\
$U_n(q), n\geq 3$ &$q^{n^2}$&$^2E_6(q)$&$q^{79}$\\
$P\Omega_{2n}^+(q),n\not=4$&$q^{2n^2-n+1}$ &$G_2(q)$&$q^{15}$\\
$P\Omega_8^+(q)$&$3q^{29}$&$^3D_4(q)$&$3q^{29}$\\
$P\Omega_{2n}^-(q)$&$q^{2n^2-2n+3}$&$^2F_4(q)$&$q^{27}$\\
$\Omega_{2n+1}(q)$& $q^{2n^2+n+1}$&$Sz(q)$&$q^6$\\
$E_6(q)$&$q^{79}$&$^2G_2(q)$&$q^8$\\
$E_7(q)$&$q^{134}$&$P\Omega_{2n}^-(q)$&$2q^{2n^2-2n+2}$\\
\hline
\end{tabular}
\label{pre18}\end{table} Let $G$ be a finite group and $ P, M$ be
subgroups of $G.$ Denote by $M \backslash G/ P,$ a set of
representatives for the double cosets of $G$ on $P$ and $M.$ Let
$\m=G/P,$ the  right cosets of $G$ on $P.$ Then $M$ acts on $\m$ by
right multiplication.

\begin{lem}\label{pre19} Let $M, P$ be subgroups of a finite group $G.$ Then

$(i)$ $M$ has $|M\backslash G/ P|$ orbits on $G/P,$ where $M$ acts
on $G/P$ by right multiplication.

$(ii)$ $(1_M^G,1_P^G)=|M\backslash G/ P|.$
\end{lem}

\begin{proof} $(i)$ Suppose that $M \backslash G/ P=\{x_1,\dots, x_k\}.$ If $i,j\in\{1,\dots,
k\}$ with $i\neq j,$ then $Px_iM\cap Px_jM=\emptyset.$  Clearly
$Px_iM$ are distinct orbits of $M$ on $G/P.$ As $G=\cup_{i=1}^k
Px_iM,$ $\{Px_iM\}_{i=1}^k$
 is a complete set of orbits of $M$ on $ G/P.$

$(ii)$ By \cite[Corollary $5.5$]{isaacs}, the number of orbits of
$M$ on $G/P$ is the inner product
$(1_M,(1_P^G)_M)=\frac{1}{|M|}\sum_{m\in M}1_P^G(m),$ and so
$(1_M,(1_P^G)_M)=|M\backslash G/P|$ by $(i).$ Now by the Frobenius
reciprocity \cite[Lemma $5.2$]{isaacs}, we have
$(1_M^G,1_P^G)=(1_M,(1_P^G)_M),$ so that $(1_M^G,1_P^G)=|M\backslash
G/P|$ as required.
  \end{proof}
Note that $(ii)$ could follow easily from the Mackey's formula.

For the representation theory of symmetric groups and the Mullineux
conjecture, we refer to \cite{jam}, \cite{ford} or \cite{fordkles}.
Let $\lambda=(\lambda_1,\lambda_2,\dots, \lambda_k)\vdash n$ be a
partition of $n,$ and let $p$ be a prime. We denote by $[\lambda]$
the Young diagram of the partition $\lambda,$ which consists of $n$
nodes placed in decreasing rows. The {\em rim} of $[\lambda]$ is its
south-east border. The {\em $p$-rim} of $[\lambda]$ is defined as
follows: Beginning at the top right-hand corner of $[\lambda],$ the
first $p$ nodes of the rim are in the $p$-rim. Then skip to the next
row, and take the next $p$ nodes of the rim. Continue until we reach
the end of the rim. The last $p$-segment may contain fewer than $p$
nodes (see \cite{ford} or \cite{fordkles}). Let $h_1$ be the number
of nodes in the $p$-rim of $\lambda,$ and let $r_1$ be the number of
rows in $\lambda.$ Delete the $p$-rim and repeat the process to get
$h_1, r_1,\dots, h_k, r_k,$ where $h_{k+1}=r_{k+1}=0,$ but $h_k\neq
0\neq r_k.$ The {\em Mullineux symbol} is a $2\times k$ matrix,

$$M(\lambda)=\begin{pmatrix}
     h_1 & h_2  &  \dots& h_k  \\
      r_1& r_2 & \dots&r_k
\end{pmatrix}.$$
Now the $p$-regular partition $m(\lam)$ of $n$ is defined via  $$M(m(\lambda))=\begin{pmatrix}
     h_1 & h_2  &  \dots& h_k  \\
      s_1& s_2 & \dots&s_k
\end{pmatrix},$$ where $$\ep_i = \left\{\begin{array}{rr}
   0,  & \mbox{if $p \mid h_i$}    \\
     1, & \mbox{ if  $p\nmid h_i$}
\end{array}\right.$$ and $s_i=h_i-r_i+\ep_i.$ Note that the partition $m(\lam)$
can be reconstructed from the Mullineux symbol $M(m(\lam)).$ The following notation will be useful
if we just want to know whether a given partition is fixed under the Mulllineux map. The {\em
$p$-modular Frobenius symbol } for $\lambda,$ denoted by $Fr_p(\lambda),$ is a $3\times k$ matrix
$$Fr_p(\lambda)=\begin{pmatrix}
      a_1 & a_2  &  \dots& a_k  \\
      b_1& b_2 & \dots&b_k \\
      \ep_1&\ep_2& \dots&\ep_k
\end{pmatrix}$$ where $a_i=h_i-r_i$ and $b_i=r_i-\ep_i.$
If $\lambda$ has  $p$-modular Frobenius symbol $Fr_p(\lambda)$ as constructed above then the Mullineux map $m$ is defined by

$$Fr_p(m(\lambda))=\begin{pmatrix}
      b_1 & b_2  &  \dots& b_k  \\
      a_1& a_2 & \dots& a_k \\
      \ep_1&\ep_2& \dots&\ep_k
\end{pmatrix}$$
which means that we interchange the first two rows of
$Fr_p(\lambda).$ Therefore we see that $\lambda$ is a fixed point of
the Mullineux map if and only if the first two rows of
$Fr_p(\lambda)$ are the same (see \cite{ford}). Denote by  $sgn_n,$
the sign character of $S_n,$ which takes value $1$ at even
permutation and $-1$ at odd permutation. Now if $\lambda$ is a
$p$-regular partition of $n$ then $D^{\lambda}\bigotimes
sgn_n=D^{m(\lambda)}$ (see \cite{fordkles}), where $D^\lambda$
denote the irreducible $S_n$-module in characteristic $p,$
corresponding to the $p$-regular partition $\lambda.$ Let $\lam$ be
a $3$-regular partition of $n$ with $\lam_1\geq n-2.$ We will show
that $\lam$ is rarely  a fixed point of the Mullineux map.
\begin{lem}\label{smallpart} Assume that $p=3, n\geq 5$ and $\lam$ is a $p$-regular partition of $n.$
Suppose that $\lam_1\geq n-2.$  Then $m(\lam)=\lam$ if and only if $5\leq n\leq  6$ and $\lam=(n-2,1^2).$
\end{lem}

\begin{proof} As $\lam_1\geq n-2,$  $\lam_1=n, n-1$ or $n-2.$ It follows that $\lam=(n), (n-1,1), (n-2,2)$ or
$(n-2,1^2).$ If one of the first three cases holds then the result follows from \cite[Lemma $1.8$]{kleshsheth}.
Assume that  $\lam=(n-2,1^2).$ We first compute the Mullineux symbol $M(\lam)$ of $\lam.$  We have
$$M(\lam)=\begin{pmatrix}
      5&3& \dots&3&a    \\
      3&1& \dots&1&1
\end{pmatrix} ,$$ where  $3,$ and so $1,$ occurs $t$ times with $t=\left [\frac{n-2}{3} ]\right.-1,$ and $0\leq a=n-2-3(t+1)\leq 2.$
Hence $$Fr_3(\lambda)=\begin{pmatrix}
      2 &2  &  \dots&2&a-1  \\
      2& 1 & \dots& 1&0 \\
      1&0& \dots&0&1
\end{pmatrix}.$$ If $t\geq 1,$ or equivalently $n\geq 8,$ then clearly, the first two rows of $Fr_3(\lam)$ cannot
be equal, so that $\lam\not=m(\lam).$ Thus $5\leq n\leq 7$ and
$t=0.$ Then
$$Fr_3(\lambda)=\begin{pmatrix}
      2 &a-1  \\
      2& 0 \\
      1&1
\end{pmatrix}.$$ Observe that the first two rows of $Fr_3(\lam)$ are equal if and only if $a=0$ or $a=1.$
Since $t=0,$ $a=n-2-3=n-5,$ so that $n=5$ or $6.$
\end{proof}

We next prove a gap result between the minimal module and the second
minimal module for alternating groups in characteristic $3.$ For
$k\geq 1,$ denote by $R_n(k)$ the set of irreducible $S_n$-modules
$D$ such that $D\cong D^{\lambda}$ or $D\cong D^{m(\lambda)}$ for
some $p$-regular partition $\lambda\vdash n,$ with $\lambda_1\geq
n-k.$
\begin{lem}\label{lowalt}  Let ${\f}$ be a splitting field for $A_n$ of characteristic $p=3.$ Suppose that $n\geq 12$ and $V$
is an irreducible ${\f}A_n$-module with $dimV> n.$ Then $dimV\geq (n^2-5n+2)/2.$
\end{lem}

\begin{proof} It follows from \cite[Theorem $2.1$]{ford} that $V=D^{\lam}\downarrow_{A_n}$ with $m(\lam)\not=\lam$ or
$V=D_{\pm}^{\lam},$ where $m(\lam)=\lam.$ Let $U=D^\lam$ for the partition $\lam$ obtained above.
By \cite[Proposition $2.2$]{MagaardMalle}, either dim$U>(n-2)(n-3),$ or $U\in R_n(2).$ Observe that
$dimV=dimU$ if $m(\lam)\not=\lam$ and if $m(\lam)=\lam$ then $dimV=\frac{1}{2}dimU.$ Thus if
$dimU>(n-2)(n-3),$ then clearly, dim$V\geq
\frac{1}{2}dimU>\frac{1}{2}(n-2)(n-3)>\frac{1}{2}(n^2-5n+2).$ Therefore we can assume $U\in
R_n(2).$ Hence there exists a $3$-regular partition $\mu$ with $\mu(1)\geq n-2$ such that
$\lam=\mu$ or $\lam=m(\mu).$  As $n>10$ and $dimV>n,$ it follows from \cite[Theorem $6$]{jam} that
$\mu$ is neither $(n)$ nor $(n-1,1),$ and so $\mu=(n-2,2),$ or $(n-2,1^2).$ By Lemma
\ref{smallpart}, $\mu$ is not fixed under the Mullineux map. Also, as the Mullineux map is an
involutory map, $m(m(\mu))=\mu\not=m(\mu).$ We conclude that $D^\mu$ and $D^{m(\mu)}$ are
irreducible upon restriction to $A_n.$ Since $dimD^{m(\mu)}=dimD^\mu\otimes sgn_n=dimD^\mu,$ we
have $dimV=dimD^\lam=dimD^\mu.$ Finally, the result follows from \cite[Theorem $24.1,
24.15$]{jam1}. \end{proof}

\section{Higman rank 3 parameters and the equation}\label{higman3}
Let $G$ be a finite group acting transitively on a non-empty
finite set $\m.$ Let $P$ be a stabilizer of
a point $x\in \m$ in $G.$ Recall that the action is primitive if and only if $P$
is maximal in $G.$ Also the rank of $G$ is the number of orbits
of $P$ on $\m.$
Now suppose that $G$ is of even order and acts primitively rank
$3$ on $\m.$ So $P$ has exactly three orbits on $\m,$ namely,
$\{x\},~\Delta(x)~\text{and}~ \Gamma(x).$ Define $k= |\Delta (x)|,\:\:l= |\Gamma (x)|$ and

 $$|\Delta (x)\cap\Delta(y)|=\left \{ \begin{array}{ll}
               \mu &\mbox{if $y \in   \Gamma(x)$}\\
               \lambda &\mbox{if $y\in \Delta(x)$}
               \end{array}
               \right. $$

 $$|\Gamma (x)\cap\Gamma(y)|=\left \{ \begin{array}{ll}
               \lam_1 &\mbox{if $y \in   \Gamma(x)$}\\
               \mu_1 &\mbox{if $y\in \Delta(x).$}
               \end{array}
               \right.$$
Suppose that $k\leq l.$
\begin{lem}{\em (\cite[Lemma $5,7$]{hig}).}\label{ma7} Let G act primitively rank 3 on $\m.$ Then

(i) $|\m|=k+l+1;$

(ii) $\mu l=k(k-1-\lambda);$

(iii) $D=(\lambda-\mu)^2+4(k-\mu)$ is a square;

(iv) $\lam_1=l-k+\mu-1;$

(v) $\mu_1=l-k+\lam+1.$
\end{lem}
Let $V$ be the permutation module for $G$ on $\m$ over $\Complex,$
hence $\m$ is a basis for $V.$ Further $\Delta$ and $\Gamma$ can be
viewed as linear maps on $V,$ via the corresponding  $x\mapsto
\sum_{y\in\Delta(x)}y$ and $x\mapsto \sum_{y\in \Gamma(x)}y,$ for
$x\in\m,$ and extending linearly. We have that $\Sigma_{y\in \m} y$
is an eigenvector for $\Delta$ and $\Gamma,$ with eigenvalues $k$
and $l,$ respectively. Other eigenvalues for $\Delta$ and $\Gamma$
are as follows:

\begin{lem}{\em (\cite[Lemma $6$]{hig}).}\label{ma8} The eigenvalues different from $k$ of $\Delta$ are:\\
$$s=\dfrac{\lambda-\mu+\sqrt{D}}{2}~\mbox{ and }~
t=\dfrac{\lambda-\mu-\sqrt{D}}{2}.$$
\end{lem}

\begin{lem}{\em (\cite[$1.4.3$]{asc})}\label{ma8bis}
For $r\in \{s, t\},$ the eigenvalues for $\Gamma$ are $-(r+1).$
\end{lem}
Let $V_r,$ $r=s$ or $r=t,$ be the irreducible modules for $G$ on
$V$ which is the eigenspace for $\Delta$ with eigenvalue $r.$
Let $f_r=dim(V_r).$
\begin{lem}{\em (\cite{hig}).}\label{ma9} We have $$f_s=\dfrac{k+t(k+l)}{t-s}\mbox{ and } f_t=\dfrac{k+s(k+l)}{s-t}.$$
\end{lem}
Let $M$ be any subgroup of $G.$ Fix $x\in\m,$ and let $P=G_x.$ By
identifying $\m$ with $G/P,$ $x$ corresponds to  the coset $P$ in
$G/P.$ Consider the action of $G$ on the right cosets  $G/M$ and
form the permutation module $V_M.$ Denote by $y$ the coset $M$ as a
point in $G/M.$ Define $$d=d_x=|xM\cap \Delta(x)| ~\mbox{and~
$c=c_x=|xM \cap\Gamma(x)|$}.$$
As $G$ acts transitively on $\m$ and $G/M,$ $xG=\m$ and $G/M=yG.$
Define $\alpha:V\rightarrow V_M$ by $\alpha(xg)= \Sigma_{p\in P}ypg
~\mbox{for $g\in G$}$ and $\beta:V_M\rightarrow V$ by $\beta(yg)=
\Sigma_{m\in M}xmg ~\mbox{for~ $g\in G$}.$ Then $\alpha$ and $\beta$
are $\Complex G$ maps and $\theta=\frac{1}{|P||P\cap
M|}\beta\circ\alpha\in End_{\Complex G}(V).$ The map $\theta$ can be
written in terms of the linear maps $\Delta$ and $\Gamma$ as
follows:

\begin{lem}{\em (\cite[$2.1$]{asc}).}\label{ma10} We have $$\theta=I+\dfrac{d\Delta}{k}+\dfrac{c\Gamma}{l}.$$
\end{lem}
\begin{proof} For any $g\in G,$ we have

$ \begin{array}{lll}
               \theta(xg) &= &\dfrac{1}{|P||P\cap M|}\beta(\al(xg))\\
               &= &\dfrac{1}{|P||P\cap M|}\sum_{p\in P}\beta(ypg)\\
               &=&\dfrac{1}{|P||P\cap M|}\sum_{p\in P,m\in M}xmpg.
               \end{array}$

Let $\mathcal{O}$ be one of the three orbits of $P$ on $\m$ and
$v_{\mathcal{O}}=\sum_{u\in \mathcal{O}}u\in V.$ As $P$ acts on
$\mathcal{O},$ $xmp\in \mathcal{O}$ if and only if $xm\in
\mathcal{O},$ in which case since $P$ is transitive on
$\mathcal{O},$ $\sum_{p\in
P}xmp=\frac{|P|}{|\mathcal{O}|}v_{\mathcal{O}}.$ Moreover there are
$|M_x|=|P\cap M|$ elements in $M$ fixing $x$ and $|P\cap M|d,|P\cap
M|c$ elements $m\in M$ with $xm\in \mathcal{O}$ for
$\mathcal{O}=\Delta, \Gamma,$ respectively. Therefore
$$\theta(xg)=(v_x+\frac{dv_\Delta}{k}+\frac{cv_\Gamma}{l})(x)g=(I+\frac{d\Delta}{k}+\frac{c\Gamma}{l})(xg).$$
This proves the lemma.
 \end{proof}
For $r=s, t,$ let $\pi_r$ be
the projection of $V$ on $V_r.$

\begin{lem}{\em (\cite[$2.2, 2.3$]{asc}).}\label{ma11} Let $r=s$ or $t,$ then\\
(1) $\theta\circ\pi_r=0$ if and only if $V_r\leq ker(\theta);$\\
(2) If $\theta\circ\pi_r\neq 0$ then $\alpha:V_r\rightarrow V_M$
is an injective $\Complex G$ map;\\
(3) For $r=s,t:$ $\theta\circ\pi_r=0$ if and only if
\begin{equation}\label{eq1}
1+\dfrac{dr}{k}=\dfrac{(r+1)c}{l}.
\end{equation}
\end{lem}
\begin{proof}
Statement $(1)$ is clear as $\pi_r(V)=V_r.$ By Lemma \ref{ma9}, $V_s$ is not $\Complex
G$-isomorphic to $V_t,$ so that $\theta:V_r\rightarrow V_r.$ As $V_r$ is an irreducible $\Complex
G$-module, $\theta$ is an isomorphism if $\theta$ is non-zero. Thus $(2)$ follows from $(1).$ For
$(3),$ let $r=s,t$ and $v\in V.$ Let  $\Sigma=\Delta$ or $\Gamma.$ Then $(\Sigma\circ
\pi_r)(v)=e(\Sigma,r)\pi_r(v),$ where $e(\Sigma,r)$ is the eigenvalue of $\Sigma$ on $V_r.$ Thus
$\Sigma\circ \pi_r=e(\Sigma,r)\pi_r.$ From definition, $e(r,\Delta)=r$ and by Lemma \ref{ma8bis},
$e(r,\Gamma)=-(r+1).$ Therefore by Lemma \ref{ma10}, we obtain:

$$\theta\circ \pi_r=(I+\frac{d\Delta}{k}+\frac{c\Gamma}{l})\circ \pi_r=(1+\dfrac{dr}{k}-\dfrac{(r+1)c}{l})\pi_r.$$
Thus $\theta\circ \pi_r=0$ if and only if  $1+\dfrac{dr}{k}-\dfrac{(r+1)c}{l}=0.$ This finishes the proof.
 \end{proof}
As $G$ is a primitive rank $3$ group on $\m$ with $P$ the
stabilizer of a point $x$ in $\m,$ the permutation character
$1_P^G$ has a decomposition $1_P^G=1+\chi_s+\chi_t,$ where $1$
is the trivial character, and $\chi_s, \chi_t$ are irreducible
characters of $G,$ afforded by the irreducible modules $V_s, V_t,$ with degrees $f_s, f_t,$
respectively.
From definition, $1_P^G\leq 1_M^G$ if $1_M^G-1_P^G$ is a character of
$G.$ This is equivalent to $(\chi_r,1_M^G)>0$ for any
$r\in\{s,t\}.$ By Lemma \ref{ma11}, for $r\in\{s,t \},$  $(\chi_r,1_M^G)=0$
 if and only if equation $(\ref{eq1})$
holds for any $M$-orbits in $\m.$ Note that the parameters $c$ and
$d$ depend on $x,$ or equivalently, on the conjugate of $P.$ When we
pick a different conjugate of $P,$ parameters $c$ and $d$ change.
Thus for $r\in \{s,t\},$ if equation (\ref{eq1}) does not hold for
some point $x\in \m$ or some conjugate of $P,$ then by Lemma
\ref{ma11}, there is an injective $\Complex G$ map from $V_r$ to
$V_M$ and hence $1_P^G\leq 1_M^G.$ Otherwise we need to change to
different conjugate of $P$ or different point in $\m.$ See
Proposition \ref{odd5} for such an example.

The following result will be used frequently to show that there is no containment if $M$ has at most $2$ orbits on $\m.$
\begin{cor}\label{smallorbits}
Let $G$ be a primitive  rank $3$ group acting on a finite set $\m.$
Let $P$ be the stabilizer of a point in $\m,$ and let $M$ be any
subgroup of $G.$ If $M$ has at most two orbits on $\m,$ then
$1_P^G\not\leq 1_M^G.$
\end{cor}

\begin{proof}
By way of contradiction, suppose that $1_P^G\leq 1_M^G.$ Then $\phi=1_M^G-1_P^G$ is a character of
$G.$  Since $\phi$ and $1_P^G$ are characters of $G,$ we have that
$$(1_P^G,1_M^G)=(1_P^G,\phi+1_P^G)=(1_P^G,1_P^G)+(1_P^G,\phi)=3+(1_P^G,\phi)\geq
3,$$ By Lemma \ref{pre19}, $(1_P^G,1_M^G)$ is the number of orbits of $M$ on $G/P.$ Now, by
identifying $\m$ with $G/P,$ $M$ has $(1_P^G,1_M^G)\geq 3$ orbits on $\m,$ a contradiction.
 \end{proof}


\section{Main Hypothesis and Notations }
\label{hypo}

From now on, we assume the following set up. Let
$L=\Omega_{2m+1}(3)$ with $m\geq 2.$ Let $V$ be the natural module
for $L$ over ${\f}={\f}_{3}.$ Denote by $\m_\xi(V)$ an $L$-orbit of
non-singular points of type $\xi$ in $V.$ Let $G$ be a nearly simple
primitive rank $3$ group of type $L$ acting on $\m_\xi(V).$ Observe
that $L\unlhd G\leq I(V),$ where $I(V)$ is the full isometric group
of $V.$ Denote by $P$ the stabilizer of a point in $\m_\xi(V).$ The
letter $M$ will be reserved for the maximal subgroup of $G.$ If
$M\in \mathcal{C}(G)$ and $X$ is a group satisfying $\Omega(V)\leq
X\leq \Xi(V),$ then there exists a subgroup $H\in \mathcal{C}(\Xi)$
such that $H\cap G=M$ and the subgroup $H\cap X\in\mathcal{C}(X)$ is
called the {\em $X$-associate} of $M$ and is denoted by $M_X.$ If
$M\in \mathcal{S}(G)$ then we denote the socle of $M$ by $S.$ Then
$S$ is a non-abelian finite simple group and the full covering group
$\hat{S}$ of $S$ acts absolutely irreducible on $V$ and is not
realizable over a proper subfield of ${\f}.$ Moreover $\hat{S}$
fixes a non-degenerate quadratic form on $V$ so that the
Frobenius-Schur indicator $ind(V)=+.$
\section{Character containment for nearly simple groups of type $L$}
\subsection{Higman rank $3$ parameters for $L$}\label{charcon}
We now assume the hypotheses and notation in Section \ref{hypo} with $L=\Omega_{2m+1}(3),m\geq 2.$
There are two types of non-singular points in $V,$ namely $+$ and $-$ points. For $\xi\in \{\pm\},$
denote by $\m_\xi(V)$ the set of all non-singular points of type $\xi.$ For any $\la x\ra
\in\m_\xi(V),$ define
$$\Delta(x)=\m_\xi(V)\cap x^\perp \mbox{ and
}\Gamma(x)=\m_\xi(V)\cap(V-x^\perp-\{\la x\ra\}).$$ Then $P$ has exactly three orbits $\{\la
x\ra\}, \Delta(x)$ and $\Gamma(x)$ on the set $\m_\xi(V).$  Recall the Higman rank $3$ parameters
defined in Section \ref{higman3}. For $\xi=\pm,$ we write $\xi  1$ to denote $+1$ or $-1$ when
$\xi=+$ or $\xi=-,$ respectively.
\begin{lem}\label{ma17}
Let $\xi\in \{\pm\}$ and $\la x\ra\in \m_\xi(V).$  Then\\
(i) $|\m_\xi(V)|=\frac{1}{2}3^m(3^m+\xi  1);$\\
(ii) $k=\frac{1}{2}3^{m-1}(3^m-\xi  1)$;\\
(iii) $l=(3^m-\xi)(3^{m-1}+\xi  1);$\\
(iv) $\mu=\lambda=\frac{1}{2}3^{m-1}(3^{m-1}-\xi  1);$\\
(v) $\sqrt{D}=2\cdot 3^{m-1};$\\
(vi) $s=3^{m-1};$\\
(vii) $t=-3^{m-1}.$
\end{lem}
\begin{proof} Let $\gam=Q(x).$ We have $\rho(x)=sgn(x^\perp_V)=\xi.$ By Lemma \ref{pre9},
we obtain $\m_{\xi}(V)=\{\la v\ra\subseteq V \;|\; Q(v)={\gam}\}.$ Observe that each point $\la
v\ra$ has $2$ representatives $v$ and $-v,$ by Lemma \ref{pre10},
$|\m_\xi(V)|=\frac{1}{2}S(2m+1,x)=\frac{1}{2}(3^{2m}+\xi 3^m)=\frac{1}{2}3^m(3^m+\xi  1),$ which
gives $(i).$

From definition we get $$k=|\Delta(x)|=|\m_\xi(V)\cap x^\perp|=|\{\la v\ra\subseteq
x^\perp\;|\;Q(v)=Q(x)\}|.$$ By Lemma \ref{pre10} again,
$k=\frac{1}{2}S^\xi(2m,\gam)=\frac{1}{2}S^{\xi}(2m,\gam)=\frac{1}{2}(q^{2m-1}-\xi q^{m-1}).$ The
parameter $l$ can be computed from the relation $1+k+l=|\m_\xi(V)|.$ This proves $(ii)$ and
$(iii).$

To compute $\lambda,$ take $\la y\ra \in \Delta(x).$ Then $\rho(y)=\rho(x)=\xi, (x,y)=0$ and
$Q(y)=Q(x).$ We have $$\lambda=|\Delta(x)\cap \Delta(y)|=\frac{1}{2}|\{v\in x^\perp\cap
y^\perp\;|~Q(v)=Q(x)\}|=\frac{1}{2}S(2m-1,z),$$ where $z\in W:=x^\perp\cap y^\perp=\la
x,y\ra^\perp$ with $Q(z)=Q(x).$ We need to determine $\rho_W(z)=\mbox{sgn}(z^\perp_W).$ Since
$x^\perp=\la y\ra \perp (y^\perp\cap x^\perp)=\la y\ra \perp W$ and $W=z_W^\perp\perp \la z\ra,$ we
deduce that $x^\perp=\la y, z\ra \perp z_W^\perp,$ where $\mbox{dim}(z_W^\perp)=2m-2,
\mbox{dim}W=2m-1$ and $\mbox{dim}\la y, z\ra=2.$ As $D\la
y,z\ra=det(\mbox{diag}(-\gam,-\gam))=\square,$ by \cite[Proposition $2.5.13$]{kl}, $\mbox{sgn}\la
y,z\ra=(-)^1=-.$ It follows from \cite[Proposition $2.5.10$]{kl} that
$\mbox{sgn}(z_W^\perp)=-\mbox{sgn}(x^\perp)=-\xi.$ Thus by Lemma \ref{pre10},
$\lam=\frac{1}{2}(3^{2m-2}-\xi  3^{m-1}),$ which gives $(iv).$ The remaining parameters follow from
Lemmas $\ref{ma7}$ and $\ref{ma8}.$
 \end{proof}

\begin{cor}\label{ma18} Let $M$ be a subgroup of $G$ and $\xi\in\{\pm\}.$ Suppose that equation
(\ref{eq1}) holds for some $r\in\{s,t\}, $ and for
some $M$-orbit $\la x\ra M$ with $x\in\m_\xi(V).$ Then\\
$(i)$ If $(\xi,r)=(+,s)$ or $(-,t)$ then equation (\ref{eq1}) has
the form
\begin{equation}\label{eq2}
c-2d=\xi  3^m-1.
\end{equation}
$(ii)$ If $(\xi,r)=(+,t)$ or $(-,s)$ then equation (\ref{eq1}) has
the form
\begin{equation}\label{eq3} (\xi  3^{m-1}+1)(\xi  3^m-1+
c-2d)=2c
\end{equation}
$(iii)$ If $m\geq 2$ then
\begin{equation}\label{eq4}1+c+d\geq \frac{3^m+1}{2}\geq 3^{m-1}.
\end{equation}
\end{cor}
\begin{proof} From definitions we have  $c\geq 0$ and $d\geq 0.$ Let $A=1+c+d.$
Multiply both sides of equation (\ref{eq1}) by  $l,$ we have $(r+1)c=l+drl/{k}.$ Subtracting
$drl/k$ from both sides,
\begin{equation}
\label{eq5}
(r+1)c-\frac{drl}{k}=l.
\end{equation}

$(1)$ If $\xi=+,$ $r=s,$ then by Lemma \ref{ma17}, we have
$r=3^{m-1},l=(3^m-1)(3^{m-1}+1),$ and
$l/k={2(3^{m-1}+1)}/{3^{m-1}}.$  From (\ref{eq5}),
$$(3^{m-1}+1)c-2d(3^{m-1}+1)=(3^m-1)(3^{m-1}+1).$$ Dividing both sides
by $3^{m-1}+1,$ we get $c-2d=3^m-1=\xi  3^m-1.$ Hence  $c=2d+3^m-1.$
Thus $A=3^m+3d\geq 3^m\geq (3^m+1)/2.$

$(2)$ If $\xi=-,$ $r=t,$ then by Lemma \ref{ma17},
$r=-3^{m-1},l=(3^m+1)(3^{m-1}-1),$ and
${l}/{k}=2(3^{m-1}-1)/3^{m-1}.$ From (\ref{eq5}),
$$(-3^{m-1}+1)c+2d(3^{m-1}-1)=(3^m+1)(3^{m-1}-1).$$ Dividing both sides
by $3^{m-1}-1,$ we get $c-2d=-3^m-1=\xi  3^m-1.$ In this case,
$2d=3^m+1+c.$ Since $c\geq 0,$ $2A=2+2c+3^m+1+c\geq 3^m+3>3^m+1.$

$(3)$ If $\xi=+,$ $r=t,$ then by Lemma \ref{ma17}, $r=-3^{m-1},l=(3^m-1)(3^{m-1}+1),$ and
$l/{k}=2(3^{m-1}+1)/3^{m-1}.$ From (\ref{eq5}),  $$(-3^{m-1}+1)c+2d(3^{m-1}+1)=(3^m-1)(3^{m-1}+1)$$
or $$2c=(3^{m-1}+1)(3^{m}-1+c-2d)=(\xi  3^{m-1}+1)(\xi 3^{m}-1+c-2d),$$ and
$$2d(3^{m-1}+1)=(3^{m-1}+1)(3^m-1)+(3^{m-1}-1)c.$$ As $c\geq 0,$ $$2d(3^{m-1}+1)\geq
(3^{m-1}+1)(3^m-1),$$ or $2d\geq 3^m-1.$ Now $2A=2+2c+2d\geq 2+2c+3^m-1\geq 3^m+1.$

$(4)$ If $\xi=-,$  $r=s,$ then by Lemma \ref{ma17},
$r=3^{m-1},l=(3^m+1)(3^{m-1}-1),$ and
${l}/{k}=2(3^{m-1}-1)/3^{m-1}.$  From (\ref{eq5}),

$$2c=(3^{m-1}-1)(3^{m}+1-c+2d) =(\xi  3^{m-1}+1)(\xi 3^{m}-1+c-2d).$$ Therefore
$$(3^{m-1}+1)c=(3^{m-1}-1)(3^m+1)+2d(3^{m-1}-1).$$ As $d\geq 0,$ we have $c\geq
(3^{m-1}-1)(3^m+1)/(3^{m-1}+1).$ Thus $$A\geq
1+(3^{m-1}-1)(3^m+1)/(3^{m-1}+1)=(3^{2m-1}-3^{m-1})/(3^{m-1}+1).$$ Since
$$\frac{3^{2m-1}-3^{m-1}}{3^{m-1}+1}-\frac{3^m+1}{2}=\frac{(3^m-1+3^{m+1}(3^{m-2}-1)}{2(3^{m-1}+1)}\geq
0,$$ we obtain $A\geq (3^m+1)/2>3^{m-1}.$  \end{proof}
\subsection{Permutation characters of maximal subgroups in $\mathcal{C}$}
 By \cite[Proposition $2.6.2$]{kl}, $A=I=S\times\la -1\ra$ and $S=\Omega\la r_\square r_\boxtimes\ra.$
 Let $M\in \mathcal{C}(G)$ be a maximal subgroup of $G.$  Let $M_\Xi\in \C(\Xi)$
 be such that $M=G\cap M_\Xi.$  Then $M_I=M_{\Delta}=M_\Xi,$ and  $M_\Omega\leq M\leq M_I.$
\subsubsection{The reducible subgroups ${\mathcal C}_1$}
The reducible subgroups in $\mathcal{C}_1(\Xi)$ are all the groups
$M_\Xi$ of the forms $N_\Xi(W),$ where $dimW=\al,$  $1\leq \al\leq
2m+1$ and $W$ is either non-degenerate or totally singular. The
corresponding subgroups are of type $O_\al^{\ep_1}(3)\perp
O_{2m+1-\al}^{\ep_2}(3)$ or $P_\al,$  where
$\ep_1=sgn(W),\ep_2=sgn(W^\perp).$
The subgroup $M_\Xi$ is maximal in $\Xi$ except when $M_\Xi$ is of
type $O_2^{+}(3)\perp O_{2m-1}(3),$ as it is contained in the
subgroups of type $O_1(3)\perp O_{2m}^{\pm}(3).$
\begin{prop}\label{odd1}
Assume $M$ is of type $O_\al^{\ep_1}(3)\perp O_{2m+1-\al}^{\ep_2}(3).$ There is an $M$-orbit
on $\m_\xi(V)$ such that equation (\ref{eq1}) does not hold unless $M$ is of type
$O_1(3)\perp O_{2m}^{\ep_2}(3).$ In this case $M$ is in Table \ref{in2}.
\end{prop}
\begin{proof} As $M$ is of type $O_\al^{\ep_1}(3)\perp O_{2m+1-\al}^{\ep_2}(3),$ there exists a
non-degenerate subspace $W\leq V$ of dimension $\al$ such that
$M=N_G(W).$ Put $W_1=W,W_2=W^\perp, \ep_i=sgn(W_i).$ Write
$X_i=X(W_i),$ where $X$ ranges over the symbols $\Omega,S$ and $I.$
By \cite[Lemma $4.1.1$]{kl}, we have $M_I=I_1\times I_2$ and
$\Omega_1\times \Omega_2\leq M_\Omega.$ Without loss of generality,
we can assume dim$W_1=2b+1,$ where $0\leq  b<m.$ If $b=0,$ then the
proposition holds. Assume $b\geq 1.$ Then $W_1$ contains
non-singular points of both types. Let $\la x_i\ra, i=1,2$ be
non-singular points in $W_1$ of different types. By \cite[Lemma
$2.10.5$]{kl}, $\la x_i\ra \Omega_1=\la x_i\ra I_1, i=1,2.$
Moreover, as $I_2$ centralizes $x_i, i=1,2,$ we have $\la x_i\ra
M_\Omega=\la x_i\ra M_I, $ so that $\la x_i\ra M_\Omega=\la x_i\ra
M.$  Thus it is sufficient to compute parameters $c,d$ for subgroup
$M_\Omega$ in $L.$ Let $x\in \{x_1,x_2\},$  $\eta=\rho_{W_1}(x)$ and
$\xi=\rho_V(x).$ Since $\la x\ra M_\Omega=\m_{\eta}(W_1),$ it
follows that  $c=l(W_1)$ and $d=k(W_1),$ the parameters for
$\m_{\eta}(W_1).$ By Lemma $\ref{ma17},$
$d=\frac{1}{2}3^{b-1}(3^b-\eta  1),~c=(3^b-\eta 1)(3^{b-1}+\eta 1),$
and so $c-2d=\eta  3^{b}-1.$ Suppose that equation (\ref{eq1}) holds
for some $r=s,t$ and any $M$-orbits on $\m_\xi(V).$ If (\ref{eq2})
holds then $\eta  3^{b}-1=\xi  3^m-1.$ This implies that $b=m,$ a
contradiction. Thus (\ref{eq3}) holds. Then $(3^{m-1}+\xi 1)(3^m-\xi
2+\xi  \eta  3^{b})=2(3^{b}-\eta  1)(3^{b-1}+\eta  1).$ Observe that
$m-1\geq b.$ Assume first that $m-1=b$ and $\xi=-\eta.$ We have
$3^{m-1}+\xi  1=3^{b}-\eta  1\mbox{ and }3^m-\xi  2+\xi \eta
3^{b}=3^{b+1}+\eta  2-3^b=2(3^b+\eta  1).$ Clearly $2(3^b+\eta
1)>2(3^{b-1}+\eta  1),$ and hence equation (\ref{eq3}) does not hold
in this case. Next assume that $m-1=b$ and $\xi=\eta.$ Then
$(3^{m-1}+\xi  1)(3^m-\xi  2+\xi \eta  3^{b})=(3^b+\xi  1)(
3^{b+1}-\xi  2+3^b)=2(3^b+\eta  1)(2\cdot 3^b-\eta 1)>2(3^{b-1}+\eta
1)(3^b-\eta  1),$ a contradiction. Finally assume that $m-1> b$ so
that $m-1\geq b+1.$ Then $3^{m-1}+\xi\geq 3^{b+1}+\xi 1>3^b-\eta  1$
and $3^m-\xi  2+\xi  \eta  3^{b}\geq 9\cdot 3^b-\xi 2-3^b=8\cdot
3^b-\xi 2>2(3^{b-1}+\eta  1).$ Multiplying these two inequalities
side by side will lead to a contradiction. Thus equation (\ref{eq3})
does not hold. Therefore $b=0$ and so $W$ is a point or a
hyperplane.
 \end{proof}
Let ${\f}_q$ be a finite field of size $q$ and let  $W$ be a totally
singular subspace of $V$ of dimension $\alpha.$ By Witt's Lemma, we
can assume that $W$ has a basis $\{e_1, \dots,e_\alpha\},$ where the
vectors $e_i$ are taken from a standard basis  of $V.$ Let
$X=\langle e_{\alpha+1},\dots, e_m,f_{\alpha+1},\dots,f_m,a\rangle$
and $Y=\langle f_1,\dots, f_\alpha\rangle.$ Then $V=(W\bigoplus
Y)\perp X,$ $W^\perp=W\perp X.$ Finally let
$U=C_{I(V)}(W,W^\perp/W,V/W^\perp),$ and $N=N_{I(V)}(W,Y,X).$ The
relation among these groups is given in the following lemmas.
\begin{lem}\label{ma19a} Let $(V,{\f }_q,Q)$ be a classical orthogonal geometry with $sgn(Q)=+.$ Assume that $dim(V)=2m$
 and  $\beta=\{e_1, \dots, e_m,f_1, \dots, f_m\}$
is a  standard basis of $V.$ Let $W_1=\langle e_1, \dots,
e_m\rangle, W_2=\langle f_1, \dots, f_m\rangle,$
and $T_0=N_{I(V)}(W_1, W_2),$ where $V=\langle \beta\rangle.$ Then\\
$(i)$ $T_0\cong GL_m(q)$ and $T_0$ acts naturally on $W_1.$\\
$(ii)$  As $T_0$-modules we have $W_2\cong W_1^{*}.$ \\
$(iii)$ $T_0\cap \Omega(V)=\{x\in T_0\;|\;\mbox{det}_{W_1}(x)\in
({\f}_q^*)^2\}.$
\end{lem}
\begin{proof} This is a special case of \cite[Lemma $4.1.9$]{kl}).
\end{proof}
The next lemma is also a special case of \cite[Lemma $4.1.12$]{kl}.
\begin{lem}\label{ma19} Let $W$ be a totally singular subspace of $V.$ Keeping the notation above, we have:\\
$(i)$ $M_I=U:N;$ \\
$(ii)$ $N=T_0\times I(X),$ where
$GL_\alpha(q)\simeq T_0\leq I(W\oplus Y),$ and $T_0$ acts naturally on $W;$ and as $T_0$-modules we have $Y\cong W;$\\
$(iii)$ $U$ is a $p$-group and $U\leq \Omega(V).$
\end{lem}
It follows from $(i)$ and $(iii)$ of Lemma \ref{ma19} that $M_\Omega=U(N\cap \Omega(V)).$

\begin{prop}\label{odd2}
Assume $M$ is of type $P_\al.$  Then $M$ has at most two orbits in $\m_\xi(V)$ so that $1_P^G\not\leq 1_M^G$
by Corollary \ref{smallorbits} and so $M$ is in Table \ref{in2}.
\end{prop}

\begin{proof} We will show that $M_\Omega$ has at most two orbits on $\m_\xi(V),$ and hence we deduce that
$M$ also has at most two orbits on $\m_\xi(V).$ Assume the notation
above and  denote by $f$  the associated bilinear form of $Q.$  From
definition $M$ stabilizes a totally singular subspace $W$ of
dimension $\al.$ Let $Y$ and $X$ be defined as above. By
\cite[Proposition $2.6.1$]{kl}, $X$ has a basis
$\beta_X=\{x_1,\dots,x_s\},$ with  $s=2m+1-2\alpha=2(m-\al)+1,$ such
that $[f_{\downarrow_X} ]_{\beta_X}=\lambda {\mathbf{I}}_s,$ where
$D(X)\equiv \lambda (\mbox{mod $(\f^*)^2$}).$ Let
$\beta=\{e_1,\dots,e_\alpha,x_1,\dots,x_s,f_1,\dots,f_\alpha\}.$
Let $x\in X$ be a non-singular point with  $\xi=\rho_X(x).$ As $sgn(W\oplus Y)=+,$ $\xi=\rho_V(x).$
If $\alpha<m,$ then $s=dim(X)=2(m-\al)+1\geq 3,$ so that $X$ contains both plus and minus points.
Otherwise, $X$ has no minus points.  As $\Omega(X)\leq M_\Omega,$ and $U\leq M_\Omega,$ we have
$x(\Omega(X)U)\subseteq xM_\Omega.$ For any $v\in x\Omega(X),~ w\in W,$ we will show that there
exists $u\in U$ such that $vu=v+w,$ which implies that $|\la x\ra\Omega(X)U|=|\m_{\xi}(X)||W|.$
Thus $|\la x\ra M_\Omega|\geq \frac{1}{2}|x\Omega(X)U|=|W|
|\m_{\xi}(X)|=3^\al\frac{1}{2}3^{m-\al}(3^{m-\al}+\xi 1))= \dfrac{1}{2}3^m(3^{m-\alpha}+\xi 1).$
Therefore
\begin{equation}\label{odde1} |\la x\ra M_\Omega|\geq
\frac{1}{2}3^m(3^{m-\alpha}+ \xi 1).
\end{equation}
Let $\hat{B}=[f]_\beta.$ Then $$\hat{B}=\begin{pmatrix}0& 0& I_\alpha\\
       0& \lambda I_s& 0\\
    I_\alpha &0 &0\end{pmatrix}.$$
We have $[v]_{\beta}=(0,a,0)$ and $[w]_{\beta}=(b,0,0),$ where $a,b$ are row vectors in $\f^s$ and $\f^{\alpha},$ respectively.
Since $v$ is non-singular, $v$ is non-zero. Choose $B\in M_{s\times \al}({\f})$ such that $aB=b.$ Let $C=-\lambda^{-1}B^t,
A+A^t=-\lambda^{-1}B^tB=-\lambda CC^t,$ and
$$u=\begin{pmatrix}I_\alpha & 0 & 0\\
B&I_s& 0\\A&C&I_\alpha\end{pmatrix}.$$
Then $u\hat{B}u^t=\hat{B},$
$u$ centralizes $W, W^\perp/W$ and $V/W^\perp.$ Thus $u\in U$  and
$ vu=v+w.$

Let $z=\eta  e_1+f_1\in V,$ where $\eta=Q(x).$ Then $Q(z)=Q(x)$  so
that $\la x\ra$ and $\la z\ra$ belong to the same $\Omega$-orbit of
non-singular points in $V.$ Let $T=\frac{1}{2}T_0.$ By Lemma
\ref{ma19a}, $SL_\alpha(3)\simeq T\leq N\cap \Omega\leq M_\Omega,$
and so $TU\leq M_\Omega,$ where $T_0, N$ are as in Lemma \ref{ma19}.
Thus $z(TU)\subseteq zM_\Omega.$ We find the stabilizer $(TU)_z$ in
$TU$ of vector $z.$ For any $g\in TU,$ there exist $h\in T, u\in U$
such that $g=hu.$ We have
$$[h]_{\beta}=\begin{pmatrix}D&0&0\\0&I_s&0\\0&0&D^*
\end{pmatrix}, \mbox{ and $[u]_{\beta}=\begin{pmatrix}I_\alpha & 0 & 0\\
B&I_s& 0\\A&C&I_\alpha\end{pmatrix},$}$$ where $D=(d_{ij})\in
SL_{\alpha}(3),$ $D^*$ its inverse transpose, and
$C=-\lambda^{-1}B^t, A+A^t=-\lambda^{-1}B^tB=-\lambda CC^t.$ Thus
$[g]_{\beta}=[hu]_{\beta}=[h]_{\beta}[u]_{\beta}.$ Suppose $g\in
(TU)_z.$ Write $E_1=(\underbrace{\eta   1,0,\dots,0}_{\alpha})
\mbox{ and }F_1=(\underbrace{1,0,\dots,0}_{\alpha}).$
Then $[z]_{\beta}=(E_1,0,F_1).$ As $zg=z,$ we have $$(E_1,0, F_1)\begin{pmatrix}D & 0 & 0\\
B&I_s& 0\\D^*A&D^*C&D^*\end{pmatrix}=(E_1,0,F_1),$$ or
  $(E_1D+F_1D^*A,F_1D^*C,F_1D^*)=(E_1,0,F_1),$ hence
$$\left\{\begin{array}{lll }
      F_1D^*&=&F_1    \\
      F_1D^*C&=&0\\
      E_1D+F_1D^*A&=&E_1
\end{array}\right.$$ Since $F_1D^*=F_1,$ $$D^*=\begin{pmatrix}1&0\\b&D_1^*\end{pmatrix},D=\begin{pmatrix}1&-b^tD_1\\0&D_1\end{pmatrix},
\mbox{ and }D^{-1}=\begin{pmatrix}1&b^t\\0&D_1^{-1}\end{pmatrix},$$
where $D_1\in SL_{\alpha-1}(3)$ and $b$ is a column vector of size
$\alpha-1.$ As $  F_1D^*C=0$ and $F_1D^*=F_1,$ we have $F_1C=0.$
Hence $$C=\begin{pmatrix}0&0\\c_0&C_1\end{pmatrix},
A=\begin{pmatrix}a_{11}&a\\a_0&A_1\end{pmatrix},$$ where $C_1\in
M_{\alpha-1,s-1}(3),$  $c_0$ is a column vector of size $\alpha-1,$
$A_1\in M_{\alpha-1}(3), a_{11}\in {\f}$ and $a,a_0$ are row, column
vectors, respectively, of size $\alpha-1.$ Now as $
E_1D+F_1D^*A=E_1D+F_1A=E_1,$ we have $$(\xi,
0)\begin{pmatrix}1&-b^tD_1\\0&D_1\end{pmatrix}+(1,0)\begin{pmatrix}a_{11}&a\\a_0&A_1\end{pmatrix}=
(\eta   1,0).$$ It follows that $(\eta  1+a_{11},a-\eta  b^t
D_1)=(\eta   1, 0).$ Therefore $a_{11}=0$ and $a=\eta  b^tD_1.$
Finally as $A+A^t=-\lambda CC^t,$ we have
$$\begin{pmatrix}0&a_0^t+\eta  b^t D_1\\a_0+\eta  D_1^t
b&A_1+A_1^t\end{pmatrix}=-\lambda^{-1}\begin{pmatrix}0&0\\0&c_0c_0^t+C_1
C_1^t\end{pmatrix},$$ hence $a_0=-\xi D_1^t b,$  $A_1+A_1^t=-\lambda
(c_0c_0^t+C_1 C_1^t)$ and $$A=\begin{pmatrix}0&\eta  b^t D_1\\-\eta
D_1^t b&A_1\end{pmatrix}.$$
In summary, for any $g\in (TU)_z,$ we have $$[g]_{\beta}=\begin{pmatrix}D & 0 & 0\\
B&I_s& 0\\D^*A&D^*C&D^*\end{pmatrix}=\begin{pmatrix}D&0&0\\0&I_s&0\\0&0&D^*
\end{pmatrix}\begin{pmatrix}I_\alpha & 0 & 0\\
B&I_s& 0\\A&C&I_\alpha\end{pmatrix},$$ where $$D=\begin{pmatrix}1&-b^tD_1\\0&D_1\end{pmatrix},
D^*=\begin{pmatrix}1&0\\b&D_1^*\end{pmatrix},C=\begin{pmatrix}0&0\\c_0&C_1\end{pmatrix},
A=\begin{pmatrix}0&\xi b^t D_1\\-\xi D_1^t b&A_1\end{pmatrix},$$ $B=-\lambda C^t\in M_{s,\alpha}(3),$
with $C_1\in M_{\alpha-1,s-1}(3),$  $A_1\in M_{\alpha-1}(3),b,c_0\in M_{\alpha-1,1}(3),$
$A_1+A_1^t=-\lambda^{-1}(c_0c_0^t+C_1 C_1^t),$ $A\in M_{\alpha}(3),D,D^*\in SL_{\alpha}(3),C\in M_{\alpha,s}(3).$

We see that the subgroup of $SL_{\alpha}(3)$ generated by all
matrices $D$ is isomorphic to $U_0:SL_{\alpha-1}(3),$ where $U_0$ is
an elementary abelian subgroup of order $3^{\alpha-1}.$ Given such
$D,$ there are $3^{\alpha-1+(\alpha-1)(s-1)}$ choices for $C$ and
$3^{\frac{1}{2}(\alpha-1)(\alpha-2)}$ choices for $A.$ Therefore
$$|(TU)_z|=3^{\alpha-1}  |SL_{\alpha-1}(3)|  3^{\alpha-1+(\alpha-1)(s-1)+\frac{1}{2}(\alpha-1)(\alpha-2)}.$$
Since $|U|=3^{\alpha  s+\frac{1}{2}\alpha(\alpha-1)},$ we have
$|z(TU)|=|TU:(TU)_z|=3^{s+\alpha-1}(3^{\alpha}-1).$ Thus
\begin{equation}\label{odde2}|\la z\ra M_\Omega|\geq \frac{1}{2}3^{s+\alpha-1}(3^{\alpha}-1).
\end{equation}

It follows from (\ref{odde1}) and (\ref{odde2}) that
$|\m_{\xi}(V)|\geq |\la x\ra M_\Omega|+|\la z\ra M_\Omega|\geq
\dfrac{1}{2}3^m(3^{m-\alpha}+\xi
1)+\frac{1}{2}3^{2m-\alpha}(3^{\alpha}-1)=\frac{1}{2}3^m(3^m+\xi
1)=|\m_{\xi}(V)|.$ Therefore $|\la x\ra M_\Omega|=|\la x\ra
\Omega(X)U|,$ $|\la z\ra M_\Omega|=|\la z\ra M_\Omega U|,$ so that
$\m_{\xi}(V)=\la x\ra M_\Omega\cup \la z\ra M_\Omega.$ Hence
$M_\Omega$ has at most two orbits on $\m_{\xi}(V).$  Clearly as
$M_\Omega\leq M,$ $\la x\ra M_\Omega\subseteq \la x\ra M$ and
similarly $\la z\ra M_\Omega\subseteq \la z\ra M.$ Since
$\m_{\xi}(V)=\la x\ra M_\Omega\cup \la z\ra M_\Omega$ and $\la x\ra
M\cup \la z\ra M\subseteq \m_{\xi}(V),$ it implies that
$\m_{\xi}(V)=\la x\ra M\cup \la z\ra M.$ Thus $M$ has at most two
orbits on $\m_{\xi}(V).$
  \end{proof}

\subsubsection{The imprimitive subgroups $\mathcal{C}_2$}
Let $V$ be a vector space over a finite field ${\f}_q$ with $n=dimV.$ A {\em subspace
decomposition} $ {\mathcal{D}}=\{V_1,\dots, V_b\}$ of $V$ is a set of subspaces $V_1,\dots, V_b$ of
$V$ with $b\geq 2$ such that $V=V_1\oplus V_2\oplus \dots \oplus V_b.$ Let $\mathfrak{G}$ be a
subgroup of $GL(V).$ The {\em stabilizer} in $\mathfrak{G}$ of $\mathcal{D}$ is the group
$N_\mathfrak{G}\{V_1,\dots, V_b\},$ which is the subgroup of $\mathfrak{G},$ permuting the spaces
$V_i$ amongst themselves and denoted by $\mathfrak{G}_{\mathcal{D} }$. The {\em centralizer in }
$\mathfrak{G}$ of ${\mathcal{D}},$ is the group
$\mathfrak{G}_{(\mathcal{D})}=N_\mathfrak{G}(V_1,\dots, V_b),$ which is a subgroup of
$\mathfrak{G}$ fixing each $V_i.$ We also define
$\mathfrak{G}^{\mathcal{D}}=\mathfrak{G}_{\mathcal{D}}/ \mathfrak{G}_{(\mathcal{D})}.$ If the
spaces $V_i$ in the subspace decomposition $\mathcal{D}$ all have the same dimension $\alpha,$ then
$\mathcal{D}$ is called an {\em $\alpha$-decomposition.}  If the $V_i's$ are non-degenerate and
pairwise orthogonal, then $\mathcal{D}$ is said to be {\em non-degenerate.}   For any vector $v\in
V,$ $v$ can be written in the form $v=v_1+v_2+\dots+v_b,$ where $v_i\in V_i.$ We define the
$\mathcal{D}-length$ of $v$ to be the number of non-zero vectors $v_i$ appearing in $v,$ and denote
by $\mathcal{D}^k_b,$ the number of all points of $\mathcal{D}-$length $k, 1\leq k\leq b.$ The
members of $\mathcal{C}_2(\Xi)$ are the stabilizers in $\Xi$ of $\alpha$-decomposition
$\mathcal{D}$ of $V$ such that $\mathcal{D}$ is non-degenerate and if $\alpha=1$ then $q=p,$ a
prime.
\begin{lem}\label{ma20}
Let $M_\Omega$ be the stabilizer in $\Omega(V)$ of a non-degenerate
$\alpha$-decomposition $\mathcal{D}$ and $n=b\al.$ Then

$(i)$ If $\alpha>1,$  then $M_\Omega\cong
\Omega(V)_{(\mathcal{D})}S_b.$

 $(ii)$ If $\alpha=1,$ and $q\equiv
\pm 3\;\mbox{$(${\em mod }$8$$)$},$ then $M_\Omega\cong
\Omega(V)_{(\mathcal{D})}A_n.$
\end{lem}

\begin{proof}

 This is  Propositions $4.2.14$ and  $4.2.15$ in \cite{kl}.
\end{proof}
\begin{lem}\label{ma20b} Assume $\mathcal{D}$ is a $1$-decomposition of $V$ with $q$ odd.
For   $1\leq k\leq n,$ $|\mathcal{D}^k_n|=(q-1)^{k-1}\binom{n}{k}.$
\end{lem}

\begin{proof} Without loss of generality, we can assume $V$ has an orthonormal basis $\beta=\{v_1, \dots, v_n\}.$
If $v\in V$  has $\mathcal{D}$-length $k$ then $v$ is a
linear combination of a set of $k$ basic vectors taken from the basis $\beta,$
 with coefficients in ${\f}^*_q.$ Clearly, there are $\binom{n}{k}$
 choices for  $k$-sets, and for each
$k$-set, there are $(q-1)^{k-1}$ points of length $k.$ The result follows.
 \end{proof}
\begin{lem}\label{1decom} Assume $G$ is nearly simple primitive rank $3$ of type $\Omega_n(3)$
and $M$ is of type $O_1(3)\wr S_n$ with $n=2m+1.$ Let $\beta=\{x_1,\dots, x_n\}$ be an orthogonal basis for $V.$\\
$(1)$ If $z=x_1$ then $d_1=d_{z}=n-1$ and $c_1=c_{z}=0;$\\
$(2)$ If $z=x_1+x_2$ then $d_2=d_{z}=n^2-5n+7$ and $c_2=c_z=4n-8.$
\end{lem}
\begin{proof}
By multiplying a suitable non-zero constant to the quadratic form
$Q,$ we can assume that $V$ has an orthonormal basis
$\beta=\{x_1,\dots, x_n\}.$ Setting $r_i=r_{x_i},$ the reflection
along vector $x_i.$ Then we have $I(V)_{(\mathcal{D})}=\la r_i|1\leq
i\leq n\ra \cong 2^n ,\mbox{ and } \Omega(V)_{(\mathcal{D})}=\la
r_ir_j|1\leq i, j\leq n\ra\cong 2^{n-1}.$ For $1\leq i\not=j\leq n,$
we see that $r_{x_i-x_j}$ permutes $\{x_i, x_j\},$ and fixes $x_t$
for any $t\not\in \{i,j\}$ and so $r_{x_i-x_j}$ acts as a
transposition $(i,j).$ Thus if we denote by $J,$ the group generated
by all reflections $r_{x_i-x_j},$ where $1\leq i\not=j\leq n,$ then
$J\cong S_n$ and hence $J_1,$ the subgroup of $J$ generated by
$r_{x_i-x_j}r_{x_r-x_s},$ with $i\not=j, r\not=s$ is isomorphic to
$A_n.$ By Lemma \ref{ma20}$(ii)$,
$M_\Omega=\Omega(V)_{(\mathcal{D})}J_1.$ For any $1\leq i\not=j\leq
n,$ as  $(x_i,x_j)=0,$ $r_j$ fixes $x_i.$ Hence
$\Omega(V)_{(\mathcal{D})}$ leaves invariant the point $\langle
x_1\rangle,$ and $\langle x_1+x_2\rangle
\Omega(V)_{(\mathcal{D})}=\{\langle x_1+x_2\rangle,\langle
x_1-x_2\rangle \},$ as $(x_1+x_2)r_2r_3=x_1-x_2.$ By
\cite[$(4.2.17)$]{kl}, we have $M_I=M_\Omega\la r_3,r_{x_3-x_4}\ra.$
Thus $\la x\ra M_\Omega=\la x\ra M_I$ for any $x\in
\{x_1,x_1+x_2\},$ so that it suffices to compute the parameters for
$M_\Omega$ in $L.$ Since $n\geq 5, A_n$ acts transitively on the set
$\{1, 2, \dots, n\}.$ Thus $\langle x_1\rangle M_\Omega=\{\langle
x_1\rangle, \dots, \langle x_n\rangle\}.$ Hence $1+c_1+d_1=n.$
Moreover as $(x_i,x_1)=0$ for any $i>1,$ we have
$d_1=|x_1^{\perp}\cap \la x_1\ra M_\Omega|=|\{\langle x_2\rangle,
\dots, \langle x_n\rangle\}|=n-1,$ and so $c_1=0.$ Similarly, as
$A_n$ acts doubly transitively on $\{1, \dots, n\},$ we have
$\langle x_1+x_2\rangle M_\Omega=\langle
x_1+x_2\rangle\Omega(V)_{(\mathcal{D})}J_1=\{\langle
x_1+x_2\rangle,\langle x_1-x_2\rangle \}J_1.$ Thus by Lemma
\ref{ma20b} $1+c_2+d_2=|\langle x_1+x_2\rangle
M_\Omega|=\mathcal{D}^2_n=n(n-1).$ For any $\langle v\rangle\in
\langle x_1+x_2\rangle^\perp\cap \langle x_1+x_2\rangle M_\Omega,$
$\langle v\rangle=\langle x_i\pm x_j\rangle$ for some $i\not=j\in
\{1, \dots, n\}$ and $(v,x_1+x_2)=0.$ Clearly $\langle v\rangle$ is
generated by $x_1-x_2$ or $x_i\pm x_j$ for some $i<j\in\{3, \dots,
n\}.$ By Lemma \ref{ma20b} again, $d_2=| \langle
x_1+x_2\rangle^\perp\cap \langle x_1+x_2\rangle
M_\Omega|=\mathcal{D}^2_{n-2}+1=n^2-5n+7,$ and so $c_2=4n-8.$
 \end{proof}

\begin{prop}\label{odd3}
Assume $M$ is of type $O_1(3)\wr S_n,$ with $n=2m+1.$ There is an $M$-orbit on $\m_\xi(V)$ such
that equation (\ref{eq1}) does not hold unless  $(n,\xi,r)=(5,+,t),(7,+,t)$ or $(5,-,s),$ in which
cases $M$ has two orbits on $\m_\xi(V)$ so that $1_P^G\not\leq 1_M^G$ by Corollary
\ref{smallorbits}, and hence $M$ is in Table \ref{in2}.
\end{prop}

\begin{proof} Retain the notation in the previous lemma.
By Propositions $2.5.10,$ and $2.5.13$ in \cite{kl},
$sgn(x_1^\perp)=(-)^m$ and $sgn(x_1+x_2)^\perp=(-)^{m+1},$ as the
discriminant of the corresponding subspaces is square or non-square
respectively. When $m$ is even  $x_1$ is a plus vector and $x_1+x_2$
is a minus vector and vice versa when $m$ is odd. Let $x\in
\{x_1,x_1+x_2\}.$ We consider the following cases:

$(i)$ $\langle x\rangle$ is a plus point. If $m$ is even then we
choose $x= x_1.$ By Lemma \ref{1decom} $d=d_1=n-1,c=c_1=0.$ Then
$k=-rd,$ so $r$ must be $t$ and so $3^{m}-1=4m.$ This equation holds
only when $m=2$ and hence $n=5.$ If $m$ is odd then choose
$x=x_1+x_2,$ and hence by Lemma \ref{1decom} again
$d=d_2=(n-2)(n-3)+1,c=c_2=4n-8.$ Then $c-2d=14n-2n^2-22.$ As $n\geq
5,$ $2n^2+22>14n$ so that $c-2d<0.$ Therefore, equation (\ref{eq2})
cannot hold. If equation (\ref{eq3}) holds then
$(3^{m-1}+1)(3^m-1+c-2d)=2c.$ It follows that
$(3^{m-1}+1)(3^m-1+14n-2n^2-22)=8(2m-1).$ If $m\geq 5,$ then
$3^{m-1}+1>8(2m-1),$ hence this equation cannot hold. For $2\leq
m\leq 4,$ the equation occurs only when $m=3.$ Thus equation
(\ref{eq1}) holds only when $r=t$ and $m=3$ or $n=7.$

$(ii)$ $\langle x\rangle$ is a minus point. If $m$ is odd then
$x=x_1$ and $d=n-1,c=0.$ Then $k=-rd.$ It follows that $r=t,$ and
hence $3^m+1=4m.$ Since $m\geq 2, 3^m+1>4m,$ so that  this equation
cannot hold. If $m$ is even, then $x=x_1+x_2$ and $d=(n-2)(n-3)+1,
c=4n-8.$ We have $2d-c=2n^2+22-14n.$ If equation (\ref{eq3}) holds
then $(3^{m-1}-1)(3^m+1+2d-c)=2c,$ hence
$(3^{m-1}-1)(3^m+1+2n^2+22-14n)=8(2m-1).$ Since $2n^2+22-14n>0,$ for
any $m\geq 3,$
$(3^{m-1}-1)(3^m+1+2n^2+22-14n)>(3^{m-1}-1)(3^m+1)>8(2m-1);$ when
$m=2,$ $(3^{m-1}-1)(3^m+1+2n^2+22-14n)=24=8(2m-1).$ If equation
(\ref{eq2}) holds then $2d-c=2n^2+22-14n=3^m+1.$ We can check that
$3^m+1> 2n^2+22-14n$ for any $m\geq 2.$ Thus equation (\ref{eq1})
holds only when $m=2$ or $n=5$ and $r=s.$

To finish the proof, we need to verify that when these cases happen then equation (\ref{eq1}) also
holds for all $\langle x\rangle \in \m_{\xi}(V).$ In view of Corollary \ref{smallorbits}, we will
show that there are only two orbits of non-singular points of specified types. Firstly, suppose
that $n=5.$ Then $m=2,$ and $|\m_{\xi}(V)|=\frac{1}{2}3^m(3^m+\xi  1).$ In this case, $\langle
x_1\rangle M_\Omega$ and $\la x_1+x_2+x_3+x_4\ra M_\Omega$ are two orbits of plus points with orbit
sizes $5$ and $2^3\binom{5}{4}=40,$ respectively. As $|\m_{+}(V)|=\frac{1}{2}3^2(3^2+1)=45=5+40,$
there are only two orbits of plus points. Similarly, $\la x_1+x_2\ra M_\Omega$ and $\la
x_1+x_2+x_3+x_4+x_5\ra M_\Omega$ are two orbits of minus points with orbit sizes $2\binom{5}{2}=20$
and $2^4\binom{5}{5}=16,$ respectively. Since $|\m_{-}(V)|=\frac{1}{2}3^2(3^2-1)=36=20+16,$ there
are exactly two orbits of minus points. Finally, suppose that $n=7,$ and $\xi=+.$ Then $m=3$ and
$|\m_{+}(V)|=\frac{1}{2}3^3(3^3+1)=378.$  In this case, $\la x_1+x_2\ra M_\Omega$ and $\la
x_1+x_2+x_3+x_4+x_5\ra M_\Omega$ are two orbits of plus points with orbit sizes $2\binom{7}{2}=42$
and $2^4\binom{7}{5}=336.$ Since $336+42=378,$ there are only two orbits of plus points. This
completes the proof.
 \end{proof}
We next consider the case when $\alpha>1.$ Since dim$V$ is odd, it
follows that $\alpha$ and $b$ are both odd. Write $\alpha=2a+1,b=2b_1+1.$

\begin{prop}\label{odd4}
Assume $M$ is of type $O_\al(3)\wr S_b,$ with  $\alpha>1$ odd. There
is an $M$-orbit on $\m_\xi(V)$ such that equation (\ref{eq1}) does
not hold  so that $M$ is not in Tables \ref{in2}-\ref{in4}.
\end{prop}

\begin{proof}
We can assume that $V$ has an orthonormal basis which is the union
of orthonormal bases of all  $V_i.$ Let $N=\Omega_1\times
\Omega_2\times \dots\times \Omega_b\leq M_\Omega.$  By  \cite[Lemma
$4.2.8$]{kl}, we have $M_I=I_1\wr S_b.$ Thus $\Omega_1\wr S_b\leq
M_\Omega\leq M_I.$ Since  $\alpha>1$ is odd,  $\alpha\geq 3$ and so
$V_1$ contains both plus and minus points. Let $x_\xi\in V_1$ be a
non-singular vector of type $\xi\in \{\pm\}.$ Clearly $\la x_\xi\ra
\Omega_1\wr S_b=\la x_\xi\ra I_1\wr S_b,$ we conclude that $\la
x_\xi\ra M_\Omega=\la x_\xi\ra M_I=\la x_\xi\ra NS_b.$ Thus we only
need to compute the parameter for $M_\Omega$ in $L.$ Since $S_b\leq
M_\Omega$ permutes the $V_i's,$ and $\Omega_i$ centralizes $V_1,$
for all $i>1,$ $\la x_\xi\ra M_\Omega=\la x_\xi\ra N S_b=\la
x_\xi\ra \Omega_1 S_b=\m_{\xi}(V_1)S_b=\cup_{i=1}^b \m_{\xi}(V_i).$
Thus $A=|\la x_\xi\ra M_\Omega|=\frac{b}{2}  3^a(3^a+\xi  1)$ by
Lemma \ref{ma17}(i). Hence $A\leq \frac{1}{2}b \cdot 3^a(3^a+1).$ In
view of inequality (\ref{eq4}), it suffices to show that
$\frac{1}{2}3^m\geq \frac{1}{2}b \cdot 3^a(3^a+1).$ Since
$m=ba+b_1,$ this inequality is equivalent to $3^{ba+b_1}\geq b
3^a(3^a+1).$ As $b\geq 3$ and $a\geq 1,$ $3^{ba}\geq 3^{3a}=3^a\cdot
3^{2a}\geq 3\cdot 3^{2a}>3^{2a}+3^a.$ If we can prove that
$3^{b_1}=3^{\frac{b-1}{2}}\geq b,$ then clearly $3^{ba}\cdot
3^{b_1}\geq b  3^{a}(3^a+1),$ and we are done. To show that
$3^{b_1}\geq b,$ we will argue by induction on $b\geq 3.$ When $b=3$
then $3^{\frac{b-1}{2}}=3\geq 3=b.$ Suppose that
$3^{\frac{b-1}{2}}\geq b.$ Then
$3^{\frac{(b+1)-1}{2}}=3^{\frac{b-1}{2}}  \sqrt{3}\geq b \sqrt{3},$
by induction assumption. We have $3b^2=b^2+2b^2\geq b^2+6b >
b^2+2b+1=(b+1)^2,$ as $b\geq 3.$ Thus $b  \sqrt{3}\geq b+1.$ Hence
$3^{\frac{b+1-1}{2}}\geq b+1.$ The result follows.
\end{proof}

\subsubsection{The field extension subgroups $\mathcal{C}_3$}
Let ${\f}_\sharp$ be a field extension of ${\f}={\f}_3$ of degree
$\al,$ where $\al$ is a prime divisor of $n=dimV.$ Then $V$ acquires
the structure of an ${\f}_\sharp$-vector space in a natural way.
Write $V_\sharp$ for $V$ regarded as a vector space over
${\f}_\sharp.$ Denote by $T$ the trace map from ${\f}_\sharp$ to
${\f}.$ If $Q_\sharp$ is a quadratic form on $(V_\sharp,\f_\sharp)$
then $Q=TQ_\sharp$ is a quadratic form on $(V,\f).$ Write
$f_{\sharp}$ for the associated bilinear form of $Q_{\sharp}.$
Denote by $N=N_{{\f}_{\sharp}/{\f}}$ the norm map of ${\f}_{\sharp}$
over ${\f}.$ Let $\mu,\nu$ be the generators for $\f_\sharp^*$ and
$Gal(\f_\sharp/\f),$ respectively. Also the trace map from
${\f}_{\sharp}$ to ${\f}$ defines a non-degenerate bilinear form on
${\f}_{\sharp}.$ Let $Q_T: \f_\sharp \rightarrow \f$ be a map
defined by $Q_T(x)=-T(x^2)$ for $x\in \f_\sharp.$ Then $Q_T$ is a
quadratic form on $\f_\sharp$ and $f_T(x,y)=T(xy)$ is the bilinear
form associated to $Q_T.$ Then $({\f}_\sharp, Q_T,\f)$ is an
orthogonal geometry.

\begin{lem}\label{fieldext}  Let $\beta_T=\{\zeta_1,\zeta_2, \dots ,\zeta_\al\}$ be an ${\f}$-basis of
${\f}_\sharp={\f}_{p^\al},$ where $p=3,$ and $v_{\sharp}\in
V_{\sharp}$ be such that
$f_{\sharp}(v_{\sharp},v_{\sharp})=\lambda\in {\f}_{\sharp}^*.$ Then

$(i)$  $D({\f}_{\sharp})\equiv det(f_{\beta_T})\equiv (-1)^{\al-1}
(\mbox{mod } ({\f}^*)^2);$

$(ii)$ span$_{{\f}_{\sharp}}(v_{\sharp})$ is a non-degenerate
$\al$-subspace in $V$ with discriminant
$D({\f}_{\sharp})N(\lambda).$
\end{lem}

\begin{proof} From definition, we have $$f_{\beta_T}=\begin{pmatrix}
     T(\zeta_1^2) & T(\zeta_1 \zeta_2)& \dots & T(\zeta_1 \zeta_\al)    \\
     T(\zeta_2\zeta_1)& T(\zeta_2^2)& \dots & T(\zeta_2 \zeta_\al)\\
     \vdots & \vdots&\ddots&\vdots\\
     T(\zeta_\al\zeta_1)& T(\zeta_\al\zeta_2)& \dots & T(\zeta_\al^2)\end{pmatrix}.$$
Let $$X=\begin{pmatrix}
     \zeta_1 & \zeta_2 &  \dots & \zeta_\al  \\
     \zeta_1^p & \zeta_2^p &  \dots & \zeta_\al^p\\
     \vdots & \vdots&\ddots&\vdots\\
     \zeta_1^{p^{\al-1}}& \zeta_2^{p^{\al-1}}& \dots & \zeta_\al^{p^{\al-1}}\end{pmatrix}$$
     and $E=diag(\lambda,\lambda^p, \dots,\lambda^{p^{\al-1}}).$
As $T(a)=\sum_{i=0}^{\al-1} a^{p^i}$ for any $a\in {\f}_{\sharp},$
$X^tX= f_{\beta_T},$ hence $det(f_{\beta_T})=det(X^tX)=det(X)^2.$
Since $det(f_{\beta_T})\in {\f}^*$ and $detX\in {\f}_{\sharp},$ if
$\al$ is odd then clearly $detX\in {\f}^*,$ as ${\f}_{\sharp}$ does
not have any subfields of degree $2$ over ${\f}.$ Thus $(detX)^2\in
({\f}^*)^2=\{1\},$ so $ det(f_{\beta_T})=(detX)^2=1.$ Now suppose
that $\al=2.$ Let $\zeta$ be a root of $x^2-x-1$ in
$\overline{{\f}},$ and let $\tau=\zeta+1.$ Then $\tau^2=-1,
{\f}_{\sharp}={\f}(\tau)$ and $T(\tau)=0.$ Choose
$\beta_T=\{1,\tau\}.$ Then $$f_{\beta_T}=\begin{pmatrix}
      T(1)&T(\tau)    \\
      T(\tau)&T(\tau^2)
\end{pmatrix}=\begin{pmatrix}
     2 & 0   \\
      0&-2
\end{pmatrix}.$$

Hence $ det(f_{\beta_T})=-1.$ This proves $(i).$ Let
$\beta=\{\zeta_1 v_{\sharp},\zeta_2 v_{\sharp}, \dots ,\zeta_\al
v_{\sharp}\}$ and $W=\la v_{\sharp}\ra_{{\f}_{\sharp}}.$ As
$(\zeta_i v_{\sharp},\zeta_j v_{\sharp})=T(f_{\sharp}(\zeta_i
v_{\sharp},\zeta_j v_{\sharp}))=
T(\zeta_i\zeta_jf_{\sharp}(v_{\sharp},v_{\sharp}))=T(\lambda\zeta_i\zeta_j),$
we have $$f_{\beta}=\begin{pmatrix} T(\lambda\zeta_1^2) &
T(\lambda\zeta_1 \zeta_2)& \dots & T(\lambda\zeta_1 \zeta_\al)    \\
     T(\lambda\zeta_2\zeta_1)& T(\lambda\zeta_2^2)& \dots & T(\lambda\zeta_2 \zeta_\al)\\
     \vdots & \vdots&\ddots&\vdots\\
     T(\lambda\zeta_\al\zeta_1)& T(\lambda\zeta_\al\zeta_2)& \dots & T(\lambda\zeta_\al^2)    \end{pmatrix}.$$ Obviously $X^tEX=f_{\beta}.$
     Therefore $det(f_{\beta})=det(X)^2N(\lambda)=D({\f}_{\sharp})N(\lambda). $
 \end{proof}
\begin{prop}\label{odd5}
Assume $M$ is of type $O_{\frac{n}{\al}}(3^\al)$ with $n=2m+1$ and
$\al\mid n.$  There is an $M$-orbit on $\m_\xi(V)$ such that
equation (\ref{eq1}) does not hold unless $\al=3.$ In this case $M$
has $3$ orbits on $\m_\xi(V)$ and equation (\ref{eq1}) holds for all
$M$-orbits on $\m_\xi(V)$ with $r=s,$ and hence $M$ is in Table
\ref{in2}.
\end{prop}

\begin{proof} Let $q=3^\al,$ and $\mu,\nu$ be the generators for
$\f_\sharp^*$ and $Gal(\f_\sharp/\f),$ respectively.   As $n$ is
odd,  it follows that $\al$ is also odd. Write $\al=2\al_1+1$ and
$\frac{n}{\al}=2b+1.$ Then $m=b\al+\al_1,$ where $n=2m+1.$
Multiplying by a suitable constant to the quadratic form $Q_\sharp,$
we can assume that $D(Q_\sharp)=\square.$ By \cite[Proposition
$2.6.1$]{kl}, there exists a basis $\beta_\sharp=\{w_1,w_2, \dots
,w_{2b+1}\}$ of $(V_\sharp,Q_\sharp)$ such that $f_{\beta_\sharp}=
I_{2b+1}.$  Define
$\phi_\sharp=\phi_{Q_\sharp,\beta_\sharp}=\phi_{\beta_\sharp}(\nu).$
Then $o(\phi_\sharp)=\al.$ We will show that $\phi_\sharp\in
\Omega.$  Let $\beta_n=\{\zeta,\zeta^3,\dots,\zeta^{3^{\al-1}}\}$ be
a normal basis of ${\f}_\sharp$ over ${\f},$ and
$\beta_i=\beta_n\otimes  w_i,$ Since  $(\zeta^{3^{j}}
w_i)\phi_\sharp=\zeta^{3^{j+1}} w_i,$ we obtain
$$(\phi_\sharp)_{\beta_i}= \begin{pmatrix}
                           0 & 1 & 0& \dots& 0 \\
                          0 & 0 & 1& \dots& 0\\
                           \vdots&\vdots&\vdots&\ddots &\vdots\\
                           0 &0 &0 & \dots &1\\
                           1 &0 &0 & \dots &0
                         \end{pmatrix}.$$
As $det(\phi_\sharp)_{\beta_i}=(-1)^{\al-1}=1,$
$det(\phi_\sharp)=1.$ So $\phi_\sharp\in S.$  As $[S:\Omega]=2$ and
$o(\phi_\sharp)=\al$ is odd, $\phi_\sharp\in \Omega.$ Let
$I=I(V,{\f},Q),$ $I_\sharp=I(V_\sharp,{\f}_\sharp,Q_\sharp),$ and
$\Omega_\sharp=\Omega(V_\sharp,{\f}_\sharp,Q_\sharp).$ Then by
\cite[$(4.3.11)$]{kl}, $M_I=I_\sharp\la \phi_\sharp\ra\cong
I_\sharp\Z_\al.$ Since $\Omega_\sharp$ is perfect,
$\Omega_\sharp\leq L,$ hence $\Omega_\sharp\leq M_\Omega
N_L({\f}_\sharp)\leq M_I=I_\sharp\la \phi_\sharp\ra.$ By
\cite[Proposition $4.3.17$]{kl}, $[M_\Omega:\Omega_\sharp]=\al.$ As
$\phi_\sharp\in \Omega\cap I_\sharp\la \phi_\sharp\ra,$
$\phi_\sharp\in M_\Omega,$ and hence $M_\Omega=\Omega_\sharp\la
\phi_\sharp\ra.$ It follows from \cite[ Lemma $2.10.5$]{kl} that
$\la z\ra \Omega_\sharp=\la z\ra I_\sharp $ for any non-singular
point $\la z\ra$ in $(V_\sharp,{\f}_\sharp,Q_\sharp),$ so that $\la
z\ra M_\Omega=\la z\ra M_I,$ hence we only need to compute
parameters for $M_\Omega$ in $L.$ We first claim the following:

$(1)$ $Q_\sharp(w\phi_\sharp)=Q_\sharp(w)^{\nu}$ for any $w\in V_\sharp.$

$(2)$ $\la z\ra M_\Omega=\{\la w\ra \in V_\sharp\;|\; Q_\sharp(w)\in
\{\gam,\gam^{\nu},\dots,\gam^{\nu^{\al-1}}\}\},$ where $z\in
V_\sharp$ with $\gam=Q_\sharp(z).$

For $(1),$  assume that $w=\sum_{i=1}^{2b+1}\lambda_i w_i\in V_\sharp.$ Then
$w\phi_\sharp=\sum_{i=1}^{2b+1}(\lambda_i w_i)\phi_\sharp=\sum_{i=1}^{2b+1}\lambda_i^{\nu} w_i,$
hence $Q_\sharp(w\phi_\sharp)=(\sum_{i=1}^{2b+1}-\lambda_i^{2})^{\nu}=Q_\sharp(w)^{\nu}.$
For  $(2),$ from $(1)$ we have $Q_\sharp(z\phi_\sharp)=\gam^\nu,$  hence
$\{ Q_\sharp(z\phi_\sharp^j)\}_{j=1}^\al=\{ \gam,\gam^\nu, \dots ,\gam^{\nu^{\al-1}}\}.$
Thus $\la z\ra M_\Omega=\{\la w\ra \in V_\sharp\;|\; Q_\sharp(w)\in \{\gamma,\gamma^{\nu},\dots,\gamma^{\nu^{\al-1}}\}\}.$
For  a non-zero vector $w\in V_\sharp,$ consider $span_{{\f}_\sharp}(w)$ as an $\al$-subspace in $V.$

 $(a)$ {\bf Case $\al > 3.$}
 Let $z\in \{w_1,w_1+w_2\}.$ Then $Q_\sharp(z)=\mp 1$ and $Q(z)=TQ_\sharp(z)=\mp \al\not=0$
so  $z$ is non-singular in $V.$ Also as $Q_\sharp(z)=\mp 1$ is fixed under $\nu,$ by $(2)$ we have
$\la z\ra M_\Omega=\la z\ra \Omega_\sharp, $ and by  Lemma \ref{pre10}, $|\la z\ra M_\Omega|=|\la
z\ra \Omega_\sharp|=\frac{1}{2}(q^{2b}+\ep q^b),$ with $\ep=sgn(z^{\perp}_{V_\sharp}).$   We have
$\la z\ra M_\Omega\cap z^{\perp}=\la z\ra \Omega_\sharp\cap z^\perp=\{v\in V_\sharp\;|\;
Q_\sharp(v)=Q_\sharp(z), Tf_\sharp(v,z)=0\}.$ For $w\in \la z\ra M_\Omega\cap z^{\perp}_V,$ write
$w=\varphi f_\sharp(z,z)^{-1}z+w_0,$ where $w_0\in z^{\perp}_{V_\sharp},$ and $T(\varphi)=0.$ Then
$f_\sharp (w,z)=\varphi$ and $Q_\sharp(w_0)=Q_\sharp(z)^{-1}(Q_\sharp(z)^2-\varphi^2).$ As $T(\pm
Q_\sharp(z))=T(\pm 1)\not=0,$ $Q_\sharp(w_0)\not=0$ for any $\varphi\in \f_\sharp$ with
$T(\varphi)=0.$ When $\varphi\in \mbox{ker}T$ is fixed, as
$dim_{{\f}_\sharp}(z^{\perp}_{V_\sharp})=2b$ and $sgn(z^\perp_{V_\sharp})=\ep,$ by Lemma
\ref{pre10}, there are $q^{2b-1}-\ep  q^{b-1}$ vectors $w_0$ with
$Q(w_0)=Q_\sharp(z)^{-1}(1-\varphi^2)\not=0.$   Also $dim_{\f}(\mbox{Ker}T)=\al-1,$ we conclude
that $d_z=|\la z\ra M_\Omega\cap z^{\perp}_V|=\frac{1}{2}3^{\al-1}(q^{2b-1}-\ep
q^{b-1})=\frac{1}{6}(q^{2b}-\ep  q^b).$ Thus $c_z=\frac{1}{3}(q^{2b}+\ep  2q^b)-1=(\ep
3^{b\al}-1)(\ep 3^{b\al-1}+1),\mbox{ and  }c_z-2d_z=\ep  q^b-1=\ep  3^{b\al}-1.$ Assume (\ref{eq2})
holds. Then $c_z-2d_z=\ep  3^{b\al}-1=\xi 3^{m-1}-1,$ where $\xi=sgn(z^{\perp}_{V}).$ The latter
equation yields $m=b\al+1.$ Recall that $m=b\al+\al_1,$ hence $\al_1=\frac{1}{2}(\al-1)=1,$ this
forces $\al=3,$ a contradiction. Suppose (\ref{eq3}) holds. Then $(3^{m-1}+\xi   1 )(3^m-\xi  1
+\xi (c-2d))=2c,$ hence $(3^{m-1}+\xi 1)(3^m-\ep \xi 3^{b\al}-\xi 2 )=2(3^{b\al}-\ep
1)(3^{b\al-1}+\ep  1 ).$ We will show that $3^{m-1}+\xi
>2(3^{b\al}-\ep  1 )$ and $3^m-\ep \xi  3^{b\al}-\xi  2>3^{b\al-1}+\ep   1,$
so that after multiplying these two inequalities side by side, we
get a contradiction. For the first inequality, we have $3^{m-1}+\xi
\geq 3^{b\al+\al_1-1}-1\geq 3^{\al_1-1}3^{b\al}-1\geq 3\cdot
3^{b\al}-1,$ as $\al_1\geq 2.$ Now $2(3^{b\al}-\ep  1)\leq
2(3^{b\al}+1).$ It suffices to show that $3\cdot
3^{b\al}-1>2(3^{b\al}+1).$ This inequality is equivalent to
$3^{b\al}>3.$ This is true because $b\al
>1.$ For the second inequality, as $3^{b\al}-3>0,$ we have
$3^m-\ep \xi  3^{b\al}-\xi  2 \geq
3^{b\al+\al_1}-3^{b\al}-2=(3^{\al_1}-1)3^{b\al}-2>2\cdot
3^{b\al}-2=(3^{b\al}+1)+(3^{b\al}-3)\geq 3^{b\al-1}+\ep  1.$

$(b)$ {\bf Case $\al=3.$} Let $\omega$ be a root of  $x^3-x+1$ in
$\overline{\f}.$ Then $\la \omega\ra={\f}_\sharp^*,$ and  $KerT$ has
a basis $\{1,\omega\}$ with $T(\omega^2)=-1.$ We have $m=3b+1,
q=3^3.$ Let $x_1=\omega w_1,x_2=\omega^2 w_1,x_3=\omega^4
(w_1+w_2)\mbox{ and }y_1=\omega
(w_1+w_2),y_2=\omega^2(w_1+w_2),y_3=\omega^4 w_1.$ For $i=1, \dots
,3,$ we have  $Q_\sharp(x_i)\not=0, Q_\sharp(y_i)\not=0,$ and
$Q(x_i)=1, Q(y_i)=-1$ and so $x_i, y_i$  are non-singular  in both
$V_\sharp$ and $V.$ Also all $x_i's$ ($y_i's$)  belong to different
$\Omega_\sharp$-orbits but they are in the same $\Omega$-orbits. For
each $i=1,2,$  we have ${x_i}^\perp_{V_\sharp}=\la w_2,
\dots,w_{2b+1}\ra\mbox{ and }{y_i}^\perp_{V_\sharp}=\la w_1-w_2,w_3,
\dots,w_{2b+1}\ra, $ so that $D({x_i}^\perp_{V_\sharp})=\square,
D({y_i}^\perp_{V_\sharp})=\boxtimes,$ and hence by \cite[Proposition
$2.5.10, 2.5.13$]{kl}, $sgn({x_i}^{\perp}_{V_\sharp})=(-)^b$ and
$sgn({y_i}^{\perp}_{V_\sharp})=(-)^{b-1},$ where
$dim{x_i}^\perp_{V_\sharp}=dim{y_i}^\perp_{V_\sharp}=2b.$ For $i=3,$
as computation above, we have
$sgn({x_3}^{\perp}_{V_\sharp})=(-)^{b-1}$ and
$sgn({y_3}^{\perp}_{V_\sharp})=(-)^{b}.$ We now determine the type
of $x_i$ and $y_i$ in $(V,{\f},Q).$  Let
$U=span_{{\f}_\sharp}(w_3)\perp \dots\perp
span_{{\f}_\sharp}(w_{2b+1})\leq V$ be an ${\f}$-subspace. We have
${x_1}^\perp_{V}=\la w_1,\omega^2 w_1\ra\perp
span_{{\f}_\sharp}(w_2)\perp U,$  ${x_2}^\perp_{V}=\la \omega
w_1,(\omega^2-\omega) w_1\ra\perp span_{{\f}_\sharp}(w_2)\perp U,$
${x_3}^\perp_{V}=\la (\omega+1) (w_1+w_2),(\omega^2-1)
(w_1+w_2)\ra\perp span_{{\f}_\sharp}(w_1-w_2)\perp U,$ and similarly
${y_1}^\perp_{V}=\la (w_1+w_2),\omega^2 (w_1+w_2)\ra\perp
span_{{\f}_\sharp}(w_1-w_2)\perp U,$  ${y_2}^\perp_{V}=\la \omega
(w_1+w_2),(\omega^2-\omega) (w_1+w_2)\ra\perp
span_{{\f}_\sharp}(w_1-w_2)\perp  U,$  ${y_3}^\perp_{V}=\la
(\omega+1) w_1,(\omega^2-1) w_1\ra\perp span_{{\f}_\sharp}(w_2)\perp
U.$ By Lemma \ref{fieldext},
$D(span_{{\f}_\sharp}(w_i))=N(f_\sharp(w_i,w_i))=N(1)=1=\square$ and
$D(span_{{\f}_\sharp}(w_1-w_2))=N(f_\sharp(w_1-w_2,w_1-w_2))=N(-1)=-1=\boxtimes.$
Thus $D({x_i}^{\perp}_V)=\boxtimes$ and $D({y_i}^{\perp}_V)=\square$
for all $i=1, \dots, 3,$ and so  by \cite[Proposition $2.5.11$]{kl},
$sgn({x_i}^{\perp}_{V})=(-)^{m-1}=(-)^{b}\mbox{ and
}sgn({y_i}^{\perp}_{V})=(-)^{m}=(-)^{3b+1}=(-)^{b-1}.$ Let $z\in
\{x_i, y_i\}$ and $\gam=Q_\sharp(z).$  We have $\la z\ra
M_\Omega=\{\la w\ra\in V_\sharp\;|\; Q_\sharp(w) \in
\{\gam,\gam^3,\gam^9\}\}=\bigcup_{j=1}^3 \la
z\phi_\sharp^j\ra\Omega_\sharp .$ Observe that in
$(V_\sharp,{\f}_\sharp,Q_\sharp)$ all vectors $z\phi_\sharp^j,j=1,
\dots ,3$ have the same type, say $\ep=sgn(z^{\perp}_{V_\sharp}).$
It follows from Lemma \ref{pre10} that
$1+c_z+d_z=3\frac{1}{2}(q^{2b}+\ep q^{b})=\frac{1}{2}(3^{6b+1}+\ep
3^{3b+1}).$ For any $w\in \la z\ra M_\Omega\cap z^{\perp}_V,$
$Q_\sharp(w) \in \{\gamma,\gamma^3,\gamma^9\}$ and $T(\varphi)=0,$
with $\varphi=f_\sharp(w,z).$ Write $w=\varphi
f_\sharp(z,z)^{-1}z+w_0,$ where $w_0\in  z^{\perp}_{V_\sharp}.$ We
have $f_\sharp(w,z)=\varphi$ and $Q_\sharp(w_0)=\gamma^{-1}(\gamma
Q_\sharp(w)-\varphi^2).$

Assume that $i=1,2.$ Let $z\in \{x_i,y_i\}$ and $\gamma=Q_\sharp(z).$ Then
$sgn({z}^{\perp}_{V})=sgn({z}^{\perp}_{V_\sharp})=\ep.$ We will show that $Q_\sharp(w_0)\not=0$ for
any $\varphi\in \mbox{ker}T.$ By way of contradiction, suppose that $Q_\sharp(w_0)=0.$ Then
$\varphi^2\in \{\gam^2,\gam^4,\gam^{10}\},$ hence $\varphi\in \{\pm\gam, \pm\gam^2,\pm\gam^5\}$ or
$\varphi\in \{\pm\omega^2,\pm\omega^4,\pm\omega^8,\pm\omega^{20},\pm\omega^{10}\}.$ As the trace
map is non-zero on these values, we get a contradiction. Thus $Q_\sharp(w_0)\not=0.$ By Lemma
\ref{pre10}, $d_z=3\frac{1}{2}(q^{2b-1}-\ep q^{b-1})\cdot 3^{\al-1}=\frac{1}{2}(q^{2b}-\ep  q^b)$
as $|\mbox{ker}T|=3^{\al-1}=3^2, \mbox{dim}_{{\f}_\sharp}(z^{\perp}_{V_\sharp})=2b,$ and
$\ep=sgn(z^{\perp}_{V_\sharp}).$ Then $c_z=q^{2b}+\ep  2q^b-1.$ Hence $c_z-2d_z=\ep  3q^{b}-1=\ep
3^{m}-1.$ Therefore equation (\ref{eq2}) holds.

Assume $z\in \{x_3,y_3\}.$ Then
$sgn({z}^{\perp}_{V})=-sgn({z}^{\perp}_{V_\sharp})=-\ep,$ $\gam=\pm
\omega^8$ and $Q_\sharp(w_0)=\gam^{-1}(\gam Q_\sharp(w)-\varphi^2).$
For any $w\in  \la z\ra M_\Omega,$ $Q_\sharp(w)\in
\{\gam,\gam^3,\gam^9\}.$ If $Q_\sharp(w)=\gam$ then
$Q_\sharp(w_0)\not=0$ as $T(\varphi)=T(\pm \gam)\not=0.$ If
$Q_\sharp(w)=\gam^3$ then $Q_\sharp(w_0)=0$ if and only if
$\varphi\in\{\pm \gam^2\},$ and similarly, if  $Q_\sharp(w)=\gam^9$
then $Q_\sharp(w_0)=0$ if and only if $\varphi\in \{\pm \gam^5\}.$
By Lemma \ref{pre10}, we have $2d_z=3^2\cdot (q^{2b-1}-\ep
q^{b-1})+2[2(q^{2b-1}+\ep  (q^b-q^{b-1}))+7(q^{2b-1}-\ep q^{b-1})].$
After simplifying, we get $2d_z=3^{6b}+\ep  3^{3b+1},$ and hence
$c_z=3^{6b}-1,$ so that $ c_z-2d_z=-\ep  3^{3b+1}-1=-\ep 3^m-1.$
Thus equation (\ref{eq2}) holds.

Let $\xi=(-)^b=sgn({x_i}^\perp_V)$ and $\eta=sgn({y_i}^\perp_V),i=1,
\dots,3.$ As $\sum_{i=1}^3|\la x_i\ra
M_\Omega|=\frac{1}{2}(3^{6b+1}+\xi
3^{3b+1})+\frac{1}{2}(3^{6b+1}+\xi
3^{3b+1})+\frac{1}{2}(3^{6b+1}-\xi  3^{3b+1})$
$=\frac{1}{2}3^{3b+1}(3^{3b+1}+\xi 1)=\frac{1}{2}3^m(3^m+\xi
1)=|\m_\xi(V)|,$ $M_\Omega$ has exactly three orbits on $\m_\xi(V).$
Thus equation (\ref{eq2}) holds for all points in $\m_\xi (V).$
Similarly $M_\Omega$ has three orbits on $\m_\eta(V),$ and so
equation (\ref{eq2}) holds for all points in $\m_\eta (V).$
\end{proof}

\subsubsection{The tensor product subgroups $\mathcal{C}_4$}
Let $V_i$ be vector spaces over ${\f}$ of dimension $n_i,$ $i=1, \dots, t.$ Let $V=V_1\otimes
\cdots\otimes V_t.$ For $g_i\in GL(V_i),$ the element $g_1\otimes  \dots\otimes g_t\in GL(V)$ acts
on $V$ as follows: $(v_1\otimes \dots\otimes v_t)(g_1\otimes  \dots\otimes g_t)=v_1g_1\otimes
 \dots\otimes v_tg_t (v_i\in V_i),$ and extends linearly. Now suppose
that $\mathbf{f}_i$ is a non-degenerate symmetric bilinear form on
$V_i,$ so that $(V_i,{\f},\mathbf{f}_i)$ is an orthogonal geometry.
We next define the bilinear form $\mathbf{f}=\mathbf{f}_1\otimes
\dots\otimes \mathbf{f}_t$ on $V_1\otimes \dots\otimes V_t$ by
$\mathbf{f}(v_1\otimes \dots\otimes v_t,w_1\otimes \dots\otimes
w_t)=\prod_{i=1}^t\mathbf{f}_i(v_i,w_i)$ and extend linearly. We
write $(V,\mathbf{f})=(V_1\otimes \dots\otimes
V_t,\mathbf{f}_1\otimes \dots\otimes \mathbf{f}_t)$ for such a
structure and call a {\em tensor decomposition} and denote by
$\mathcal{D}.$ The members of $\mathcal{C}_4(\Xi)$ is the stabilizer
of a tensor decomposition $\mathcal{D}$ such that

$(a)$ $(V,\mathbf{f})=(V_1\otimes V_2,\mathbf{f}_1\otimes\mathbf{f}_2),$

$(b)$ $(V_1,\mathbf{f}_1)$ is not similar to $(V_2,\mathbf{f}_2),$

$(c)$ $\mathbf{f}_i$ are symmetric, $(n_1,\ep_1)\not=(n_2,\ep_2),$
where $n_i=dimV_i\geq 3,$ $\ep_i=sgnV_i,$ and $dimV_1<dimV_2.$
\begin{prop}\label{odd6}
Assume $M$ is of type $O_{n_1}(3)\otimes O_{n_2}(3),$ with $n_1<n_2.$
There is an $M$-orbit on $\m_\xi(V)$ such that equation (\ref{eq1})
does not hold  so that $M$ is not in Tables \ref{in2}-\ref{in4}.
\end{prop}

\begin{proof}   Let $v=v_1\otimes v_2\in V,$ where $v_i\in V_i$ are
non-singular vectors. Then $\la v\ra$ is a non-singular point. By \cite[$(4.4.14)$]{kl},
$M_I=I_1\otimes I_2.$ We have $\Omega_1\times \Omega_2\unlhd M_\Omega\leq I_1\otimes I_2=M_I.$ As
all $\Omega_i$ act transitively on $\m_{\rho_{V_i}(v_i)}(V_i),$ it follows that $\la v\ra M_I=\la
v\ra M_\Omega=\la v_1\otimes v_2\ra(\Omega_1\times \Omega_2)=\la v_1\ra\Omega_1\otimes \la
v_2\ra\Omega_2.$ Therefore $\la v\ra M_\Omega=\la v\ra M$ since $M_\Omega\leq M\leq M_I.$ Write
$n_1=2a+1, n_2=2b+1.$ Then $|\la v\ra M_\Omega|\leq |\la v_1\ra \Omega_1||\la v_2\ra
\Omega_2|\leq\frac{1}{2}3^a(3^a+1)3^b(3^b+1)<3^{2(a+b)}$ as $3^a+3^b+1<3^{a+b}$ and $b>a\geq 1.$
Since $dimV=2m+1=(2a+1)(2b+1),$ $m=2ab+a+b.$ Then $m-1-2(a+b)=a(b-2)+b(a-1)+a-1\geq 0$ so that
$3^{m-1}\geq3^{2(a+b)}>|\la v\ra M_\Omega|.$ This violates (\ref{eq4}) so that (\ref{eq1}) cannot
hold.
 \end{proof}
\subsubsection{The tensor product subgroups $\mathcal{C}_7$}
Let $V_1$ be an $\al$-dimensional vector space over ${\f},$ and
assume that $\mathbf{f}_1$ is a non-degenerate symmetric bilinear
form on $V_1.$ For $i=1, \dots, b,$ let $(V_i,\mathbf{f}_i)$ be a
classical geometry which is similar to $(V_1,\mathbf{f}_1).$ For
each $i,$ denote by $\eta_i$ the similarity from
$(V_1,\mathbf{f}_1)$ to $(V_i,\mathbf{f}_i)$ satisfying
$\mathbf{f}_i(v\eta_i,w\eta_i)=\lam_i\mathbf{f}_1(v,w)$ for all
$v,w\in V_1,$ where $\lam_i\in {\f}^*$ is independent of $v$ and
$w.$ Thus we obtain a tensor decomposition $\mathcal{D}$ given by
$(V,\kappa)=(V_1,\mathbf{f}_1)\otimes \dots\otimes
(V_b,\mathbf{f}_b),$ where $V=V_1\otimes  \dots\otimes V_b$ and
$\kappa=\mathbf{f}_1\otimes \dots\otimes \mathbf{f}_b.$ Let
$X_i=X(V_i,\mathbf{f}_i)$ for $X\in \{\Omega,S,I,\Lambda,\Xi,A\}.$
Define $\Xi_{\mathcal{D}}=\Xi_{(\mathcal{D})}S_b.$ The members of
$\mathcal{C}_7(\Xi)$ are the groups $\Xi_\mathcal{D}$ with $b\geq 2$
described as above.

\begin{prop}\label{odd7}
Assume $M$ is of type $O_\al(3)\wr S_b$ with $\al\geq 5.$  There is an $M$-orbit
on $\m_\xi(V)$ such that (\ref{eq1}) does not hold  so that $M$ is not in Tables \ref{in2}-\ref{in4}.
\end{prop}

\begin{proof}  We have $\Omega_\al(3)\wr S_b\unlhd M_\Omega\leq
O_\al(3)\wr S_b.$ Let $v_1=v\otimes v\otimes \dots \otimes v$ and
$v_2=v_1+w\otimes w\otimes \dots \otimes w\in V,$ where $v\not=w$
belong to some orthogonal basis of $V_1.$ As $S_b$ fixes $v_i,$ $\la
v_i\ra M_\Omega=\la v_i\ra(\Pi_{i=1}^b \Omega_\al(3)).$ Also $\la
v_i\ra M_\Omega=\la v_i\ra M_I=\la v_i\ra M.$ For $i=1,2,$ the
stabilizers of $\la v_i\ra $ in $\Pi_{i=1}^b\Omega_\al(3)$ contain a
subgroup which is isomorphic to $\Pi_{i=1}^b\Omega_{\al-2}(3).$ Thus
$|\la v_i\ra M_\Omega|=|\la v_i\ra(\Pi_{i=1}^b \Omega_\al(3))|\leq
[\Omega_\al(3):\Omega_{\al-2}(3)]^b <3^{(4a-3)b},$ where
$a=\frac{\al-1}{2}\geq 2.$ Now $m-1=\frac{1}{2}((2a+1)^b-3).$
Consider the following function in variable $x\in [2,+\infty),$
where $b\geq 2,$ $f(x):=\frac{1}{2}((2x+1)^b-3)-(4x-3)b.$ We have
$f'(x)=b(2x+1)^{b-1}-4b\geq b(2x+1)-4b\geq b>0.$ Hence $f(x)\geq
f(2)=\frac{1}{2}g(b),$ where $g(b)=5^b-10b-3$ and $b\geq 2.$ By
induction on $b\geq 2,$ $g(b)>0.$ Thus $3^{m-1}>3^{b(4a-3)}>|\la
v_i\ra M_\Omega|.$ This contradicts (\ref{eq4}). Thus (\ref{eq1})
cannot hold.
 \end{proof}
These are all the maximal subgroups in $\mathcal{C}(G)$ of $G.$ We
next consider maximal subgroups in $\mathcal{S}(G).$
\subsection{Permutation characters of maximal subgroups in $\mathcal{S}$}
In this section, we consider the maximal subgroup
$M\in\mathcal{S}({G}).$ By definition of $\mathcal{S},$ $M$ is an
almost simple group and the socle $S$ of $M$ is a non-abelian simple
group.  Then the full covering group $\hat{S}$ of  $S$ acts
absolutely irreducible on $V,$ the natural module for $G$ and
preserves a non-degenerate quadratic form on $V,$ that is,
$ind(V)=+.$

We note that some of the small cases in this section will be handled
by using GAP \cite{GAP}. We will describe how to do this at the end
of the paper.


\subsubsection{Embedding of alternating and symmetric groups}
Recall the construction of the fully deleted permutation module for alternating groups in \cite[p.
185]{kl}. Let $\{\ep_1,\dots, \ep_n\}$ be a standard basis for ${\f}_p^n,$ and let $w_0=\ep_1+
\dots+\ep_n\in {\f}_p^n.$ Put $U=w_0^\perp, W={\f}_p w_0,$ and $V=U/(U\cap W).$ Define
$e_i=\ep_i-\ep_{i+1}, i=1,\dots, n-1.$ Then $\{e_i\}_{i=1}^{n-1}$ is a basis for $V$ if $p$ does
not divide $n,$ and $\{e_i+U\cap W\}_{i=1}^{n-2}$ is a basis for $V$ if $p|n.$ Define $\ep_p(n)$ to
be $1$ if $n$ is divisible by $p,$ otherwise, $\ep_p(n)=0.$ Then $dimV=n-1-\ep_p(n).$ Let $Q$ be
the quadratic form on $V$ induced from the quadratic form associated to the natural bilinear form
on ${\f}_p^n.$ Then $(V,{\f}_p,Q)$ is a classical orthogonal geometry and $A_n\leq \Omega(V).$ To
simplify the notation, we always write $e_i$ instead of  $e_i+U\cap W.$
\begin{lem}\label{minmod}

Assume $M$ is almost simple of type $A_n$  with $n\geq 10$ and $V$
is the fully deleted permutation module for $A_n$ in characteristic
$3.$ Let $v=\ep_1-\ep_2, w=\ep_1+\ep_2-\ep_3-\ep_4\in V.$  Then

$(1)$ $|\la v\ra
M|=\frac{1}{2}n(n-1),d_v=\frac{1}{2}(n-2)(n-3),\mbox{ and }
c_v=2n-4;$

$(2)$ $|\la w\ra M|=\frac{1}{8}n(n-1)(n-2)(n-3),
c_w=2n^3-25n^2+111n-172$ and

 $d_w=2+4(n-4)^2+\frac{1}{8}(n-4)(n-5)(n-6)(n-7).$
\end{lem}
\begin{proof}  Observe that $\la v\ra S_n=\{\la \ep_i-\ep_j\ra\:|\;i\neq j\in \{1, \dots,n\}\}$ and
$\la w\ra S_n$ is the set of all points of the form $\la
\ep_i+\ep_j-\ep_r-\ep_s\ra,$ where $i,j,r,s\in \{1, \dots,n\}$ and
pair-wise distinct. As $n\geq 10$ and $A_n$ is $(n-2)$-transitive on
the index set $\{1, \dots, n\},$ we have $\la x\ra A_n=\la x\ra S_n$
for any $x\in \{v,w\}.$ Hence $\la x\ra M=\la x\ra S_n$ for $x\in
\{v,w\}.$ Since $v=\ep_1-\ep_2,$ it is clear that if $g\in S_n$ and
$(\ep_1-\ep_2)g=\ep_1-\ep_2,$ then $g$ must fix indices $1$ and $2.$
Thus $(S_n)_v\simeq S_{n-2}.$ Similarly,
$(\ep_1+\ep_2-\ep_3-\ep_4)g=\ep_1+\ep_2-\ep_3-\ep_4$ implies that
$g$ must fix the partitions $\{1, 2\}, \{3,4\}.$ Thus $g\in
S_2\times S_2\times S_{n-4}.$ Therefore $|M:M_{\la
v\ra}|=\frac{1}{2}[S_n:S_{n-2}]=\frac{1}{2}n(n-1),$ and $|M:M_{\la
w\ra }|=\frac{1}{2}[S_n: S_2\times S_2\times
S_{n-4}]=\frac{1}{8}n(n-1)(n-2)(n-3).$

$(i)$ Parameters for $v.$ We have  $\la u\ra \in \la v\ra M\cap
v^\perp$ if and only if $u=\ep_i-\ep_j\not\in \mbox{${\f}_3$}v,
i\not=j $ and $(\ep_i-\ep_j,\ep_1-\ep_2)=0.$ This happens only if
$\{i,j\}\cap \{1,2\}=\emptyset,$ or $i, j\in \{3,4,\dots, n\}.$
There are $\binom{n-2}{2}$ such points $\la u\ra .$ Thus
 $d_v=|\la v\ra M\cap v^\perp|=\frac{1}{2}(n-2)(n-3),\mbox{ and }
c_v=2n-4.$

$(ii)$ Parameters for $w.$\\
We will show that $c=2n^3-25n^2+111n-172$ and
 $d=|\la w\ra M\cap
w^\perp|=2+4(n-4)^2+\frac{1}{8}(n-4)(n-5)(n-6)(n-7).$
 For any $\la u\ra \in \la w\ra M\cap w^\perp,$ $u=\ep_i+\ep_j-\ep_r-\ep_s,$ where
 $|\{i,j,r,s\}|=4,$ and $(u,w)=0.$ Denote by $supp(u)$ the
 set of non-zero indices of $\ep_i$ appearing in $u.$ We consider the
 cases:

 $(1)$ $supp(u)\cap supp(w)=\emptyset.$ Then $supp(u)\in
 \{5,6,\dots,n\}.$ Hence there are
 $(n-4)(n-5)(n-6)(n-7)/8\mbox{ points.}$

 $(2)$ $|supp(u)\cap supp(w)|=1.$ There are no such $u,$ since $(u,w)\not=0.$

 $(3)$ $|supp(u)\cap supp(w)|=2.$ Suppose that $i,j\in
 \{1,2,3,4\}.$ Then either $u=\ep_i-\ep_j+\ep_r-\ep_s,$ where $\ep_i-\ep_j\in\{\ep_1-\ep_2,\ep_3-\ep_4\},$
or $u=\ep_i+\ep_j-\ep_r-\ep_s,$ where
$\ep_i+\ep_j\in\{\ep_1+\ep_3,\ep_1+\ep_4,\ep_2+\ep_3,\ep_2+\ep_4\};$
and $r, s\in\{5,\dots, n\}.$ There are $2(n-4)(n-5)$ and
$4\binom{n-4}{2}$ points respectively. Thus there are
$4(n-4)(n-5)$ points in this case.

 $(4)$ $|supp(u)\cap supp(w)|=3.$ Suppose that $i,j,r\in
 \{1,2,3,4\}.$ Then $u=\pm \ep_i\pm\ep_j\pm\ep_r\pm\ep_s,$ where $s\in\{5,\dots,
 n\},$  $\ep_i,\ep_j,\ep_r$ with their signs appearing exactly
 as in $w,$ and sign of $\ep_s$ is chosen so that there are $2$
 minuses and $2$ pluses. There are $\binom{4}{3}(n-4)=4(n-4)$ such points.

  $(5)$ $|supp(u)\cap supp(w)|=4.$ There are just $2$ points in
 this case: $\{\ep_1-\ep_2+\ep_3-\ep_4,\ep_1-\ep_2-\ep_3+\ep_4\}.$

Therefore $d_w=\frac{1}{8}(n-4)(n-5)(n-6)(n-7)+2(n-4)(n-5)+4(n-4)(n-5)+2$
$=\frac{1}{8}(n-4)(n-5)(n-6)(n-7)+2(n-4)^2+2,$ and $c_w=\frac{1}{8}n(n-1)(n-2)(n-3)-d-1.$
 \end{proof}

\begin{prop}\label{odd8}
Assume $M$ is almost simple of type $A_n,$ with $n\geq 10,$ and $V$
is the fully deleted  permutation module for $A_n$ in characteristic
$p=3.$ Further assume that $n-1-\ep_3(n)=2m+1.$ There is an
$M$-orbit on $\m_\xi(V)$ such that equation (\ref{eq1}) does not
hold so that $M$ is not in Tables \ref{in2}-\ref{in4}.
\end{prop}
\begin{proof}
Let $v=\ep_1-\ep_2, w=\ep_1+\ep_2-\ep_3-\ep_4\in V.$ Then $Q(v)=1,$
and $Q(w)=-1.$ Hence $v, w$ are non-singular vectors in $V.$ We see
that  $n-1-\ep_3(n)$ is odd if and only if $n=6k+2, n=6k+3$ or
$n=6k+4.$ By Lemma \ref{minmod}$(1),$ $|\la v\ra
M|=\frac{1}{2}n(n-1).$

Assume that $n\geq 13.$ Then $m-1=\frac{1}{2}(n-2-\ep_3(n))-1\geq \frac{1}{2}(n-5), $ as
$\ep_3(n)\leq 1,$ and so $3^{m-1}\geq 3^{(n-5)/2}>n(n-1)/2=1+c+d,$ violating (\ref{eq4}) and so
equation (\ref{eq1}) cannot hold.


Hence we can assume that $10\leq n\leq 12.$ Since $n=6k+2,6k+3$ or $6k+4,$ it follows that $n=10.$
Then $n-1-\ep_3(n)=9, d_v=28,c_v=16, m=4.$ If equation (\ref{eq1}) holds, then either $c_v-2d_v=\xi
3^4-1=\xi  81-1=-40,$ or $(\xi 27+1)(\xi 81-41)=32.$ These equations clearly cannot hold with
$\xi=\pm .$ By Lemma \ref{minmod}$(2)$, $|\la w\ra M|=n(n-1)(n-2)(n-3)/8.$ If $n\geq 23,$ then
$(3^{m}+1)/2>n(n-1)(n-2)(n-3)/8$ so that equation (\ref{eq1}) cannot hold by (\ref{eq4}). Thus we
can assume that $10\leq n\leq 22.$ Then
$n\in\{10,14,15,16,20,21,22\}.$\\
$(a)$ $n=10.$ Then $n-\ep_3(n)=9, m=4,$ $d=191, c=438$ and $c-2d=56.$\\
$(b)$ $n=14.$ Then $n-\ep_3(n)=13, m=6,$ $d=1032, c=1970$ and $c-2d=-94.$\\
$(c)$ $n=15.$ Then $n-\ep_3(n)=13, m=6,$ $d=1476, c=2618$ and $c-2d=-334.$\\
$(d)$ $n=16.$ Then $n-\ep_3(n)=15, m=7,$ $d=2063, c=3396$ and $c-2d=-730.$\\
$(e)$ $n=20.$ Then $n-\ep_3(n)=19, m=9,$ $d=6486, c=8048$ and $c-2d=-4924.$\\
$(f)$ $n=21.$ Then $n-\ep_3(n)=19, m=9,$ $d=8298, c=9656$ and $c-2d=-6940.$\\
$(g)$ $n=22.$ Then $n-\ep_3(n)=21, m=10,$ $d=10478, c=11466$ and $c-2d=-9490.$\\
We can check that equation (\ref{eq1}) cannot hold in any of these cases.
 \end{proof}
\begin{prop}\label{odd9} %
Assume $M$ is almost simple of type $ A_n$ with $n\geq 12,$ and $V$
is not the fully deleted  permutation module for $A_n$ in
characteristic $p=3.$ There is an $M$-orbit on $\m_\xi(V)$ such that
equation (\ref{eq1}) does not hold so that $M$ is not in Tables
\ref{in2}-\ref{in4}.
\end{prop}
\begin{proof}
As $n\geq 12,$ by Lemma \ref{lowalt}, we have  $dim(V)=2m+1\geq
(n^2-5n+2)/2$ so that $m\geq (n^2-5n)/2.$ However when  $n\geq 12,$
$3^{m-1}\geq 3^{(n^2-5n-4)/{4}}>n!=|Aut(A_n)|.$  Thus equation
(\ref{eq1}) cannot hold in view of (\ref{eq4}).  \end{proof}

Let $\lambda=(\lambda_1^{a_1},\lambda_2^{a_2},
\dots,\lambda_h^{a_h})$ be a $p$-regular partition. Then $\lambda$
is called a $JS$- partition if
$\lambda_i-\lambda_{i+1}+a_i+a_{i+1}\equiv 0\mbox{ (mod $p$)}.$

\begin{prop}\label{odd10}
Assume $M$ is almost simple of type $S=A_n$ with $5\leq n\leq 11.$
There is an $M$-orbit on $\m_\xi(V)$ such that equation (\ref{eq1})
does not hold unless $n=9$ and  $V\cong D^{(8,1)},$ in which case
$M$ has at most $2$ orbits on $\m_\xi(V)$ so that $1_P^G\not\leq
1_M^G$ by Corollary \ref{smallorbits} and hence
$(L,S)=(\Omega_7(3),A_9)$ is in Table \ref{in3}.
\end{prop}

\begin{proof}
Using information on the $p$-modular representations of alternating groups and their covering
groups in \cite{abc}, we need to consider the cases given in Table \ref{ma30}.
\begin{table}
\centering \caption{Small degree representations of ssome
alternating groups.}\label{ma30}
\begin{tabular}{l|c|c|c|c|c|c|c|r}
\hline

$A_n$ &$A_6$&$A_7$ &$A_7$ &$A_8$&$A_8$ & $A_8$ &$A_9$&$A_9$\\\hline

dim$D^{\lambda}$ &9&13 &15&7&13 & 21 &21&7\\

$\lambda$&$(4,2)$&$(5,2)$&$(5,1^2)$&$(7,1)$&$(6,2)$&$(6,1^2)$&$(7,1^2)$&$(8,1)$\\

$m(\lambda)$&$(2^2,1^2)$&$(3,2,1^2)$&$(3,2^2)$&$(4,3,1)$&$(3^2,1^2)$&$(3^2,2)$&$(4,3,2)$&$(4^2,1)$\\\hline
\end{tabular}
\end{table}
$(i)$ Let $\lambda=(8,1).$ Then $m(\lambda)=(4^2,1)\neq \lambda.$ By
\cite[Theorem $2.1$]{ford}, $D^{\lambda}\downarrow_{A_9}$ is
irreducible. Thus $A_9\leq S_9\leq \Omega_7(3),$ and there are two
classes of $S_9$ in $\Omega_7(3).$ As $8-1+1+1=9\equiv 0 (\mbox{mod
$3$}),$ $\lambda$ is a $JS$-partition, and hence  by \cite[Theorem
$0.3$]{klesh1},
$D^{\lambda}\downarrow_{S_8}=D^{\lambda(1)}=D^{(7,1)}.$ Then since
$(7,1)\neq m(7,1)=(4,3,1),$ we have: $A_8\leq S_8\leq S_9.$ In this
case, $D^{\lambda}$ is the fully deleted permutation module for
$S_9$ over ${\f}_3.$ Then $n-1-\ep_3(n)=7,$ and $m=3.$ Let
$v=\ep_1-\ep_2, w=\ep_1+\ep_2-\ep_3-\ep_4\in V.$ There are only two
orbits of type $\rho_V(v),$ with representatives $v=\ep_1-\ep_2$ and
$\ep_1+\ep_2+\ep_3+\ep_4-\ep_5=\ep_1+\ep_2+\ep_3+\ep_4-\ep_5-\ep_6-\ep_7-\ep_8,$
and one orbit of type $\rho_V(w).$ Thus equation (\ref{eq1}) holds
for both types of  points.

$(ii)$ If $\lambda=(7,1)$ or $\lambda=(4,3,1),$ then
$A_8<S_8<S_9<\Omega_7(3),$ since
$D^{(8,1)}\downarrow_{S_8}=D^{(7,1)}$ and
$D^{(7,1)}\downarrow_{A_8}$ is irreducible.

$(iii)$ If $\lambda=(6,1^2)$ or $(3^2,2),$ then
$A_8<S_8<S_9<\Omega_{21}(3),$ since
$D^{(7,1^2)}\downarrow_{S_8}=D^{(6,1^2)}$ and
$D^{(6,1^2)}\downarrow_{A_8}$ is irreducible.

$(iv)$ If $\lambda=(5,2)$ or $(3,2,1^2),$ then
$A_7<S_7<S_8<\Omega_{13}(3).$

$(v)$ Now, if $\lambda=(n-2,1^2),$ where $n=7\mbox{  or } 9,$ then
$D^{\lambda}=\wedge^2(D^{(n-1,1)}).$ As $D^{(n-1,1)}$ is the fully deleted permutation module for
$S_n,$ we can apply the construction above for fully deleted module. Let $v=e_1\wedge
e_3=(\ep_1-\ep_2)\wedge(\ep_3-\ep_4),$ and $w=(e_1-e_2)\wedge (e_3-
e_4)=(\ep_1+\ep_2+\ep_3)\wedge(\ep_3+\ep_4+\ep_5).$ Then $v,w$ are non-singular points of different
types in  $D^{\lambda}.$ We have $|vS_n|=\frac{1}{2}\frac{n!}{2(n-4)!}\mbox{ and }
|wS_n|=\frac{1}{2}\frac{n!}{2(n-5)!}.$ We then get contradictions by using (\ref{eq4}).

$(vi)$  If $\lambda=(4,2)$ then $dimD^{\lambda}=9, m=4,$ and we have
an embedding $A_6\leq S_6\leq \Omega_9(3).$ As $A_6\cong
L_2(9)<A_{10}<\Omega_9(3),$  $M=N_{G}(A_6)$ is not maximal in  $G.$

$(vii)$ If $\lambda=(6,2),$ then $dimV=13, m=6.$ We have $A_8\leq
S_8\leq \Omega_{13}(3).$
Using \cite{GAP}, $S_8$ has two points $w_1, w_2$ of different types
with the same orbit sizes $315.$ We have $(c,d)=(230,84)$ or
$(c,d)=(212,102).$ We see that equations (\ref{eq2}) and (\ref{eq3})
cannot hold.
 \end{proof}

\subsubsection{Embedding of groups of Lie types in  cross-characteristic}

\begin{prop}\label{odd11}
Assume $M$ is almost simple of type $S,$ where $S$ is a finite
simple group of Lie type in cross-characteristic. There is an
$M$-orbit on $\m_\xi(V)$ such that equation (\ref{eq1}) does not
hold unless $(L,S)=(\Omega_7(3),PSp_6(2)),$ in which case $M$ has
only two orbits on $\m_\xi(V)$ so that $1_P^G\not\leq 1_M^G$ by
Corollary \ref{smallorbits}, and so $(L,S)$ is in Table \ref{in3}.
\end{prop}
\begin{proof}
Suppose equation (\ref{eq1}) holds for some $r\in \{s,t\}$ and some
$M$-orbit $\la x\ra M$ with $\la x\ra \in \m_\xi(V).$ Then $|\la
x\ra M|\geq \frac{1}{2}(3^m+1),$ by (\ref{eq4}). On the other hand,
using the lower bounds  for degrees of cross-characteristic
representations of finite simple groups of Lie type (see \cite{LaSe}
and \cite{SeZa}), we have $2m+1\geq e(S),$ so that $(3^{m}+1)/2\geq
(3^{(e(S)-1)/2})/2.$ Moreover $|\la x\ra M|\leq|M|\leq |Aut(S)|.$ It
follows that  $(3^{\frac{e(S)-1}{2}}+1)/2\leq |Aut(S)|.$ By this
condition, we get a finite list as in Table $\ref{Ta6}.$
\begin{table}
\centering
\caption{Small groups in cross-characteristic.}\label{Ta6}
\begin{tabular}{l|r}
\hline $S$ &$(3^{\frac{e(S)-1}{2}}+1)/2\leq |Aut(S)|$\\\hline
$L_2(q)$&$L_2(q), 3\leq q\leq 68$\\
$L_n(q), n\geq 3$&$L_3(2), L_4(2),L_5(2),L_3(4)$\\
$PSp_{2n}(q), n\geq 2$&$S_4(5),S_4(7),S_6(5)$\\
&$S_4(2),S_6(2),S_8(2),S_4(4)$\\
$U_n(q), n\geq 3$&$U_n(2),3\leq n\leq 7,U_3(4),U_3(5)$\\
$P\Omega_{2n}^+(q), n\geq
4$&$\Omega_8^+(2)$\\
$P\Omega_{2n}^-(q),n\geq 4$&$\Omega_8^-(2)$\\

$\Omega_{2n+1}(q),n\geq 3,q$ odd & \\

$E_6(q)$&\\
$E_7(q)$&\\
$E_8(q)$&\\
$F_4(q)$&$F_4(2)$\\
${}^2E_6(q)$&\\
$G_2(q)$&\\
$^3D_4(q)$&${}^3D_4(2)$\\
${}^2F_4(q)$&${}^2F_4(2)'$\\
$Sz(q)$&$Sz(8)$\\
${}^2G_2(q)$&\\\hline
\end{tabular}
\end{table}
Next we shorten the list by using information on cross-characteristic representations of small
groups together with \cite{abc} and \cite{GAP}. Further $V$ is an absolutely irreducible ${\f}_3\hat{S}$-module
which is self-dual with $ind(V)=+1,$ and $dimV$ is odd.\\
$(i)$ Case $S=L_2(q),2\leq q\leq 68.$ As $L_2(2), L_2(3)$ are not
simple, $L_2(4)\cong L_2(5)\cong A_5,$ and $q$ is a prime power, we
can assume that $7\leq q\leq 67.$

Case $q\equiv 1$ $(\mbox{mod $4$}).$ By \cite[Table $2$]{hm} and the
fact that $dimV$ is odd, we have $dimV\in \{(q+1)/2,q\}.$ Using
(\ref{eq4}) again, we only need to consider the following cases:

$(1)$ $q\in \{13,17,29,37,41,49\}$ and $dimV=(q+1)/2;$

$(2)$ $q=13,17$ and $dimV=q.$

If $dimV=(q+1)/2,$ then $q$ must be a square in ${\f}_3,$ so that
$q\equiv 1$ $(\mbox{mod $3$})$ and hence $q=13,37,49.$ If $dimV=q,$
then $3\nmid (q+1),$ so that $q=13.$

If $(S,dimV)=(L_2(13),7)$ then $L_2(13)<G_2(3)<\Omega_7(3)$ by
\cite{atlas} so $M$ is not maximal in $G.$

If $(S,dimV)=(L_2(13),13)$ then $L_2(13)<\Omega_{13}(3).$ Let $w$ be
a non-singular eigenvector of an element of order $13.$ Then $|\la
w\ra S|=14$ and hence $|\la w\ra M|\leq |Out(S)||\la w\ra S|=2\cdot
14<3^{m-1}=3^5$ so that equation (\ref{eq1}) cannot hold. For other
type of point, if $M=S$ then there exists  a point $\la u\ra$ with
$|\la u\ra S|=1092$ and $(c,d)=(734,357).$ If $M=S\cdot 2$ then
there exists  a point $\la v\ra$ with $|\la v\ra M|=2184$ and
$(c,d)=(1469,714).$ We check that equation (\ref{eq1}) cannot hold
in any of these cases.

If $(S,dimV)=(L_2(37),19)$ then $m=9$  and the eigenvector $w$ of an
element of order $37$ is non-singular and $|\la w\ra S|=38,$
$|Out(S)|=2$ and hence $|\la w\ra M|\leq 38\cdot 2<3^{m-1},$ so that
equation (\ref{eq1}) cannot hold for this point. For other type of
point, there exists  a point $\la u\ra$ with $|\la u\ra
S|=|S|=25308$ and $(c,d)=(16919,8388).$ We see that equation
(\ref{eq1}) cannot hold.

If $(S,dimV)=(L_2(49),25)$ then $m=12$ and there exist two non-singular vectors of different type
$u_+,u_-\in V$ which are eigenvectors of an element of order $7$ in $S$ such that $|\la u_\xi\ra
S|\leq 8400,\xi=\pm.$ As $|Out(S)|=4,$ the latter inequality yields  $|\la u_\xi \ra M|\leq
|Out(S)||\la u_\xi\ra S|\leq 4\cdot 8400<3^{m-1}.$ In view of (\ref{eq4}), equation (\ref{eq1})
cannot hold.


Case $q\equiv 3$ $(\mbox{mod $4$}).$ As in previous case, by
\cite[Table $2$]{hm}, we have $dimV=q$ and $3\nmid q+1.$ Using
(\ref{eq4}) again, we get $q\in \{7,19\}.$

If $(S,dimV)=(L_2(7),7)$ then $L_2(7)<\Omega_7(3)$ but $\Omega_7(3)$
has no maximal subgroup with socle $L_2(7)$ by \cite{atlas}.

If $(S,dimV)=(L_2(19),19)$ then $m=9$ and there exist two non-singular vectors of different type
$u_+,u_-\in V$ which are eigenvectors of an element of order $5$ in $S$ such that $|\la u_\xi\ra
S|\leq 342,\xi=\pm.$ As $|Out(S)|=2,$ this implies  $|\la u_\xi \ra M|\leq |Out(S)||\la u_\xi\ra
S|\leq 2\cdot 342<3^{m-1}.$ In view of (\ref{eq4}) equation (\ref{eq1}) cannot hold.


Case $q\equiv 0$ $(\mbox{mod $2$}).$ As in previous case, by
\cite[Table $2$]{hm}, we have $dimV\in \{q-1,q+1\}.$ Using
(\ref{eq4}) again, we have $q=8,16.$ By \cite{abc}, the only
possibility is $q=8$ and $dimV=7.$ However by \cite{atlas},
$L_2(8)\leq G_2(3)\leq \Omega_7(3).$


$(ii)$ Case $L_n(q), (n,q)\in \{(3,2),(4,2),(5,2), (3,4)\}.$ As
$L_3(2)\cong L_2(7), $ $L_4(2)\cong A_8,$ which have been done
above, we can exclude these groups.

If  $S=L_3(4),$ then $Out(S)\cong 2\times S_3\cong D_{12}.$ By
\cite{abc}, $dimV\in \{15, 19, 45,63\}.$ By using (\ref{eq4}), we
only need to consider the representations of degrees $15$ and $19.$


Assume first that $dimV=15.$ Then $m=7.$ If $M=L_3(4)$ then there
exists two non-singular vectors $u_i,i=1,2$ with $|\la u_i\ra
L_3(4)|=2016, $ $(c_1,d_1)=(1250,765),$ and $(c_2,d_2)=(1350,660).$
Similarly if $M=L_3(4)\cdot 2_1$ then $(c_1,d_1)=(2600,1431),$
$(c_2,d_2)=(1250,765);$ if $M=L_3(4)\cdot 2_2$ or $L_3(4).2_3$ then
$(c_1,d_1)=(2720,1311),$ $(c_2,d_2)=(2810,1221);$ if $M=L_3(4)\cdot
2^2$ then $(c_1,d_1)=(1250,765),$ $(c_2,d_2)=(5327,2736).$ We check
that equation (\ref{eq1}) cannot hold.

Assume that $dimV=19.$ Then $m=9.$ If $M=L_3(4)$   then there exist
two non-singular vectors of different type  $u_i,i=1,2$ which are
the eigenvectors of elements of order $7,5,$ respectively and
$(c_1,d_1)=(707,252),(c_2,d_2)=(1190,825).$ For the remaining
extensions of $L_3(4),$ there exist two non-singular vectors of
different type  $u_i,i=1,2,$ which are the eigenvectors of an
element of order $5$ such that the parameters $(c_i,d_i),i=1,2,$ are
as follows: If $M=L_3(4)\cdot 2_2$ or $M=L_3(4)\cdot 2_3$ then
$(c_1,d_1)=(1190,825),(c_2,d_2)=(1100,915).$ If $M=L_3(4)\cdot 2_3$
then $(c_1,d_1)=(2561,1470),(c_2,d_2)=(1100,915).$ If $M=
L_3(4)\cdot S_3$ or $L_3(4)\cdot D_{12}$ then
$(c_1,d_1)=(3842,2205),(c_2,d_2)=(4220,1827).$ If $M=L_3(4)\cdot 3$
then $(c_1,d_1)=(3932,2115),(c_2,d_2)=(4220,1827).$ We can check
that equation (\ref{eq1}) cannot hold in any of these cases.

Finally if $S=L_5(2),$ then $dimV\geq 155.$ But $(3^m+1)/2\geq
(3^{77}+1)/2>|Aut(L_5(2))|$ so that equation (\ref{eq1}) cannot
hold.

$(iii)$ Case $S\in \{S_4(5), S_4(7), S_6(5)\}.$ If $S\cong S_4(5)$
or $S_6(5)$ then  the smallest odd degree non-trivial  irreducible
representations of $S$ has degree $13,$ and $63,$ respectively.
However since the smallest field of definitions of  these
representations are  quadratic extensions of ${\f}_3,$ (cf.
\cite{hm}), $L$ cannot embed in $\Omega_{13}(3)$ and
$\Omega_{63}(3).$  By \cite[Theorem $2.1$]{gmst}, if $\Phi$ is a
representation of $S$ which is not the smallest representation, then
$dim\Phi\geq (q^n-1)(q^n-q)/(2(q+1)) ,$ which are $40$ and $1240,$
respectively. But then inequality (\ref{eq4}) cannot hold. If
$S=S_4(7),$ then the smallest non-trivial representation in
characteristic $3$ of $S$ is a Weil representation of degree $25.$
However, the Frobenius -Schur indicator is $\circ,$ (cf. \cite{hm}),
which means that $S_4(7)$ fixes no quadratic form. Thus $S_4(7)$
cannot embed in $\Omega_{25}(3).$ If $\Phi$ is a non-trivial
representation of $S_4(7),$ which is not the smallest representation
of $S,$ then $dim(\Phi)\geq 126,$ but this again violates
(\ref{eq4}). Thus equation (\ref{eq1}) cannot hold.


$(iv)$ Case $S\in \{S_4(2), S_6(2),S_8(2),S_4(4)\}.$ By the
isomorphism  $S_4(2)\cong S_6,$ it follows that $A_6\cong S_4(2)'.$
Thus we can exclude this case. If $S=S_4(4)$ then $dimV\geq 51$ and
$(3^m+1)/2\geq (3^{25}+1)/2>|Aut(S_4(4))|$ so equation (\ref{eq1})
cannot hold.  For $S_6(2),$ and $S_8(2),$ we need to consider the
following cases $(S_6(2),7), (S_6(2),21),$ $(S_6(2),27)$
$(S_8(2),35).$

If $S=S_6(2)$ and $dim(V)=7,$ then by \cite{atlas}, $S_6(2)$ is a
maximal subgroup of $\Omega_7(3)$ and it has only two orbits on
$\m(V)$ so that  equation (\ref{eq1}) holds for both types of
points.


If $(S,dimV)=(S_6(2),21)$ then $m=10,$  $Out(S)=1$ and $S_6(2)\leq
\Omega_{21}(3).$ There exist two non-singular vectors of different
type  $u_i,i=1,2$ which are the eigenvectors of elements of order
$5,12,$ respectively  and
$(c_1,d_1)=(212,165),(c_2,d_2)=(2132,1647).$


If $(S,dimV)=(S_6(2),27)$ then $m=13$ and $S_6(2)\leq
\Omega_{27}(3).$ There exist two non-singular vectors of different
type  $u_i,i=1,2$ which are the eigenvectors of an element of order
$5$   and $(c_1,d_1)=(96968,48183),(c_2,d_2)=(47912,24663).$


If $(S,dimV)=(S_8(2),35)$ then $m=17, Out(S_8(2))=1$ and $S_8(2)\leq
\Omega_{35}(3).$ There exist two non-singular vectors of different
type  $u_i,i=1,2$ which are the eigenvectors of an element of order
$5$   and $(c_1,d_1)=(119,0),(c_2,d_2)=(256094,129465)$ and so
equation (\ref{eq1}) cannot hold.


$(v)$ Case $S\in \{U_n(2), 3\leq n\leq 7, U_3(4), U_3(5)\}.$ As
$U_3(2)\cong 3^2.Q_8$ and $U_4(2)\cong S_4(3),$ we can rule out
these cases. If $S=U_5(2),$ then by \cite{abc}, the smallest odd
degree non-trivial $3$-modular representation of $S$ has degree
$55.$ Thus $(3^m+1)/2\geq (3^{27}+1)/2>|Aut(S)|.$ If $S=U_3(4),$
then  $dimV\in \{13,39,75\}.$  However if $dimV\geq 39$ then
(\ref{eq4}) cannot hold, and if  $dimV=13,$ then by \cite{hm},
$U_3(4)$ fixes no quadratic form. If $S=U_7(2),$ then by \cite{hm},
we have $dimV>250$ and so (\ref{eq4})  cannot hold.  For the
remaining cases, using \cite{abc}, we need to consider the following
cases: $(U_3(5), 21)$ and $(U_6(2),21).$


If $(M,dimV)=(U_3(5),21)$ then we have $m=10,$  $Out(S)=S_3$ and
$M\leq \Omega_{21}(3).$ For each extension $M$ of $S,$ there exist
two non-singular vectors of different type  with parameters
$(c_i,d_i),i=1,2,$ as follow: if $M=U_3(5)$ then
$(c_1,d_1)=(7033,3466),$ $(c_2,d_2)=(27145,14854);$ if
$M=U_3(5)\cdot 2$ then $(c_1,d_1)=(55437,28562),$
$(c_2,d_2)=(55371,28628).$

If $(S,dimV)=(U_6(2),21)$ then we have $m=10,$  $M\leq
\Omega_{21}(3)$ and $Out(S)=S_3.$ For each extension $M$ of $S,$
there exist two non-singular vectors of different type  with
parameters $(c_i,d_i),i=1,2,$ as follow: if $M=U_6(2)$ then
$(c_1,d_1)=(7033,3466),$ $(c_2,d_2)=(27145,14854);$ if
$M=U_6(2)\cdot 2$ then $(c_1,d_1)=(55437,28562),$
$(c_2,d_2)=(55371,28628).$ We can check that equation (\ref{eq1})
cannot hold in any of these cases.

$(vi)$ Case $S=O_8^+(2).$ By \cite{abc}, the smallest odd degree
non-trivial $3$-modular representation of $S$ has degree $35,$ and
the second smallest odd degree one has degree $147.$ Using
(\ref{eq4}) again, we only need to consider the $3$-modular
representation of $O_8^+(2)$ of degree $35.$  We have $m=17$ and
there are two non-singular points $v_i, i=1,2,$ of different type
with $|\la v_1\ra S|=|120$ and $|\la v_2\ra S|=90720.$ Then $|\la
v_i\ra M|\leq |Out(S)||\la v_i\ra S|<3^{m-1}$ and so equation
(\ref{eq1}) cannot hold.

$(vii)$ Case $S=O_8^-(2).$ By \cite{abc},  $dimV\geq 203.$ Then inequality (\ref{eq4}) cannot hold.

$(viii)$ Case $S=F_4(2).$ By \cite{hm},  $dimV>255.$ Then inequality (\ref{eq4}) cannot hold.

$(ix)$ Case $S=G_2(4).$ By \cite{abc}, $dimV\geq 649.$ Clearly, $3^{m-1}>4^{15}\geq |Aut(S)|.$

$(x)$ Case $S=Sz(8).$ By \cite{abc},  $dimV=35.$ Then inequality (\ref{eq4}) cannot hold.

$(xi)$ Case $S={}^3D_4(2).$ By \cite{abc}, either $dimV=25$ or
$dimV\geq 351.$ If $dimV\geq 351,$ then $m\geq 174$ and clearly
$3^{m-1}\geq 3^{174}>3\cdot 2^{29}\geq |Aut(S)|.$ When $dimV=25,$
the group $S$ is not maximal in $\Omega_{25}(3)$ as ${}^3D_4(2)\leq
F_4(3)\leq \Omega_{25}(3)$ (see \cite{finitesub}).

$(xii)$ Case $S={}^2F_4(2)'.$ By \cite{abc}, $2m+1\geq 77.$ Then $3^{m-1}\geq 3^{37}>|Aut(S)|.$
 \end{proof}

\subsubsection{Embedding of groups of Lie type in defining characteristic}

Let $k$ be an algebraically closed field of characteristic $p.$ Let
$\mathfrak{G}$ be a simply connected, simple algebraic group over
$k.$ Fix a maximal torus $T$ and a Borel subgroup $B$ containing
$T.$ Let $U$ be the unipotent radical of $B.$ Then $B=UT.$ Let
$\Phi$ be the root system of $\mathfrak{G},$ select a system of
positive roots $\Phi^+$ from $\Phi,$ with corresponding fundamental
roots $\Pi=\{\alpha_1,\dots,\alpha_{\ell}\}.$ Let
$\{\lambda_1,\dots,\lambda_{\ell}\}$ be the fundamental dominant
weights and $X^+$ be the set of dominant weights. If $q$ is a power
of $p,$ then we put $X_q=\{\sum_{i=1}^\ell c_i \lambda_i\;|\;c_i\in
\Z, 0\leq c_i\leq q-1\}.$ A {\em p-restricted dominant weight} is a
weight that lies in $X_p.$ Let $\mathfrak{L}$ be the simple Lie
algebra over $\Complex$ of the same type as $\mathfrak{G}.$ For each
dominant weight $\lambda\in X^+,$ there exists an irreducible
$\mathfrak{L}$-module $V(\lambda)$ of highest weight $\lambda,$ and
a maximal vector $v^+$ (unique up to scalar multiplication). Let
$\mathcal{U}=\mathcal{U}(\mathfrak{L})$ be the universal enveloping
algebra of $\mathfrak{L},$ and $\mathcal{U}_{\Z}$ be the Kostant
$\Z$-form of $\mathcal{U}$ (see \cite[$26.3$]{HumphLie}).  Now
$\mathcal{U}_{\Z}v^+$ is the minimal admissible lattice in
$V(\lambda),$ and $\mathcal{U}_{\Z}v^+\otimes_{\Z} k $ is a
$k\mathfrak{G}$-module of highest weight $\lambda,$ also denoted by
$V(\lambda),$ and called a {\em Weyl module} for $\mathfrak{G}$
(\cite[$27.3$]{HumphLie}). The Weyl module $V(\lambda)$ has a unique
maximal submodule $J(\lambda)$ and
$L(\lambda)=V(\lambda)/J(\lambda)$ is an irreducible
$k\mathfrak{G}$-module of highest weight $\lambda.$ We will label
the Dynkin diagram as in \cite{lub}.

Assume that $\hat{S}$ is simply connected of type $A_{\ell}$ or
${}^2A_\ell$ over ${\f}_{q},$ where $q=p^f,$ and $\mathfrak{G}$ be
the corresponding simply connected, simple algebraic group over $k,$
such that $\hat{S}=\mathfrak{G}_{\sigma}$ for some suitable
Frobenius map $\sigma.$ Let $N=N(\hat{S})$ be the natural module for
$\hat{S}.$ We collect here some information about $L(\lambda)$ for
some special dominant weights $\lambda.$

$(1)$ Let $0< c < p,$ and $\lambda=c\lambda_1$ or
$\lambda=c\lambda_{\ell}.$ Then $L(\lambda)$ has all weight spaces
of dimension $1.$ $L(\lambda)$ is isomorphic to the space of
homogeneous polynomials of degree $c,$ that is, $L(\lambda)\cong
S^c(N).$ In particular, dim$L(\lambda)=\frac{(\ell+c)!}{\ell!c!}$
(\cite[$1.14$]{Seitz87}).

$(2)$ If $\ell>1,$ then $L(\lambda_i)\cong \bigwedge^iN$ and
dim$L(\lambda_i)=\binom{\ell+1}{i}$ (see \cite{CPS}).

$(3)$ Let $\lambda=n_1\lambda_1+n_2\lambda_2+\dots+ n_{\ell}\lambda_{\ell}$ be a dominant weight.
Then $L(\lambda)$ preserves a non-degenerate bilinear form if and only if $n_1=n_{\ell},
n_2=n_{\ell-1},\dots.$ Thus if $\ell$ is even then $L(\lambda_i)$ leaves invariant no
non-degenerate bilinear form, and if $\ell$ is odd then $L(\lambda_i), i\neq \frac{\ell+1}{2}$ does
not preserve any such form. Let $\lambda=\lambda_{(\ell+1)/2}.$ If $\ell\equiv -1$ (mod $4$) then
$L(\lambda)$ fixes a symmetric bilinear form and it fixes an alternating bilinear form if
$\ell\equiv 1$ (mod $4$) (see   \cite{Bourbaki} Chapter VIII $\S 13$ Table $1,$ p. $217$).

The following constructions for adjoint modules of groups of type $A_\ell$
and ${}^2 A_\ell$ are taken from \cite{lie}, pp.$ 491-492.$

$(4)$ We  construct the irreducible module $L(\lam_1+\lam_\ell)$ as
follows: Let $V:=V_1/(V_1\cap V_2),$ where $V_1=\{A\in
M_{\ell+1}(q)\:|\: Tr(A)=0\},V_2=\{aI_{\ell+1}\:|\:a\in {\f}_q\}.$
Let $\hat{S}$ act on $V_1$ by conjugation. Then $V$ is an
irreducible $\hat{S}$-module of dimension
$\ell^2+2\ell-\ep_p(\ell+1).$ The bilinear form on $V_1$ is defined
as follows: for any $A,B\in V_1,$ $(A,B)=Tr(AB).$ We can check that
$\hat{S}$ preserves this bilinear form. Also, $V$ has a basis
consisting of $E_{i,j}, 1\leq i<j\leq \ell+1, E_{ii}-E_{i+1,i+1},
i=1, \dots, \ell-\ep_p(\ell+1).$

$(5)$ Let $\hat{S}=SU_n(q)$ and $\lam=\lam_1+\lam_\ell,$ where
$n=\ell+1.$  Let $V_2=\{aI_{n}\:|\:a\in {\f}_{q^2}\},$ $V_1=\{A\in
M_{n}(q^2)\:|\: Tr(A)=0, A=\overline{A}^t\},$ and set
$V:=V_1/(V_1\cap V_2),$ where the map $A\mapsto \overline{A}$ is the
map that raises each entry to its $q^{th}$-power. Let $\hat{S}$ act
on $V_1$ as in $(4)$. The bilinear form on $V_1$ is also defined as
in $(4).$ We can check that $\hat{S}$ preserves this bilinear form
and $L(\lam_1+\lam_\ell)\cong V.$ Moreover fix a generator $\mu$ of
${\f}_{q^2}^*,$ $V$ has a basis consisting of $E_{i,j}+E_{j,i}, \mu
E_{i,j}+\overline{\mu}E_{j,i} 1\leq i<j\leq \ell+1,
E_{ii}-E_{i+1,i+1}, i=1, \dots,\ell-\ep_p(\ell+1).$

\begin{prop}\label{odd12}
Assume $M$ is almost simple of type $S,$ where $\hat{S}$ is  simply
connected of type $A_{\ell}$ or $^2A_\ell$ over ${\f}_{3^f}.$ There
is an $M$-orbit on $\m_\xi(V)$ such that equation (\ref{eq1}) does
not hold  so that $M$ is not in Tables \ref{in2}-\ref{in4}.
\end{prop}

\begin{proof}
Assume  (\ref{eq1}) holds for some $r\in\{s, t\}$ and any $M$ orbits
in $\m_\xi(V).$ We consider the case $f=1$ and $f>1$ separately. We
use the notation $L^\ep_n(q)$ to denote $L_n(q)$ when the sign $\ep$
is $+$ and $U_n(q)$ when $\ep=-.$

{\bf Case $f=1.$} We can assume that $\ell\geq 2.$ By
\cite[$31.1$]{Humlinear} and \cite[$\S 13$]{St4}, there exists a
$3$-restricted dominant weight $\lambda\in X_3$ such that $V\cong
L(\lambda).$  As $|Aut(S)|=|Aut(L^\ep_{\ell+1}(3))|\leq
3^{(\ell+1)^2},$ where $\ep=\pm.$ It follows from (\ref{eq4}) that
$3^{m-1}<3^{(\ell+1)^2}.$ Hence $m<(\ell+1)^2+1$ and $dimV<
2(\ell+1)^2+3.$ We need to look for all dominant weights $\lambda\in
X_3$ such that $L(\lambda)$ is self-dual, has dimension less than
$2(\ell+1)^2+3$ and has odd degree. If $\ell\geq 18,$ then
${\ell^3}/{8}\geq 2(\ell+1)^2+3,$ and so by \cite[Theorem
$5.1$]{lub}, $\lambda$ is one of the following $3$-restricted
dominant weights $\{\lambda_1, \lambda_{\ell},
\lambda_2,\lambda_{\ell-1},2\lambda_1, 2\lambda_{\ell},
\lambda_1+\lambda_{\ell}\}.$ Since $L(\lambda)$ is self-dual, the
only possibility for $\lambda$ is $\lambda_1+\lambda_{\ell}.$ If
$\ell<18,$ then by Theorem $4.4,$ Appendix $A_6$ through $A_{21}$ in
\cite{lub}, either $\lambda=2\lambda_2$ when $\ell=3$ or
$\lambda=\lambda_1+\lambda_{\ell}$ for  $2\leq \ell\leq 17.$

Suppose that $\ell\geq 4.$ We have
$dimL(\lambda_1+\lambda_{\ell})=\ell^2+2\ell-\ep_p(\ell+1).$ As
$dimL(\lambda_1+\lambda_{\ell})$ is odd, it follows that
$\ell=6b_1+1, 6b_1+2 $ or $6b_1+3$ for $b_1\geq 1.$ Consequently
$\ell\geq 7.$ As constructed above,
$L(\lambda_1+\lambda_{\ell})\cong V:=V_1/(V_1\cap V_2).$ Let $U$ be
the subgroup of $\hat{S}$ consisting of all matrices of the form
$diag(I_2,A),$ where $A\in SL_{\ell-1}^\ep(3).$ Then $U\cong
SL^\ep_{\ell-1}(3).$ For $\xi=\pm , $ let $x_\xi=E_{1,2}+ \xi
E_{2,1}+V_1\cap V_2,$ when $\ep=+,$ and $x_+=E_{1,2}+
E_{2,1}+V_1\cap V_2, x_-=\mu E_{1,2}+ \overline{\mu} E_{2,1}+V_1\cap
V_2$ when $\ep=-.$ Then $x_\xi\in V$ and $Q(x_\xi)\not=0.$ It
follows that $\la x_\xi\ra$ is non-singular in $V,$ of plus or minus
type depending on $\xi$ and $\ell.$  As $V_1\cap V_2$ is fixed under
natural action of  $U,$ and clearly, $U$ centralizes $x_\xi,$ it
follows that $U\leq S_{\la x_\xi\ra},$ the stabilizer of $\la
x_\xi\ra$ in $S.$ We have $1+c+d\leq
|Aut(S):U|=[Aut(L^\ep_{\ell+1}(3)):SL^\ep_{\ell-1}(3)]\leq 2\cdot
3^{2\ell-1}(3^\ell+1)(3^{\ell+1}+1)<3^{4\ell+2}.$ As
$2m+1=\ell^2+2\ell-\ep_3(\ell+1)\geq \ell^2+2\ell-1,$ $m-1\geq
(\ell^2+2\ell-4)/2.$  We have
$(\ell^2+2\ell-4)/2-(4\ell+2)=(\ell(\ell-6)-8)/2.$ If $\ell\geq 8$
then $\ell(\ell-6)-8>0,$ hence $m-1>4\ell+2.$ If $\ell=7$ then
$2m+1=\ell^2+2\ell=63,$ or $m=31$ and $m-1=30=4\ell+2.$ Thus
$m-1\geq 4\ell+2$ for any $\ell\geq 4.$ Hence $3^{m-1}\geq
3^{4\ell+2}>1+c+d.$ This contradicts to inequality (\ref{eq4}).
Therefore equation (\ref{eq1}) cannot hold in this case.

We are left with the cases $\ell=2, \lambda=\lambda_1+\lambda_2,$
$\ell=3,\lambda=\lambda_1+\lambda_3,$ and $\ell=3,
\lambda=2\lambda_2.$ If the first case holds then
$\hat{S}=SL^\ep_3(3),$ and  $dimL(\lambda_1+\lambda_2)=7.$ However
by \cite{atlas}, $\Omega_7(3)$ has no maximal subgroup with socle
$L^\ep_3(3).$

Assume that  $(S,L)=(L_4(3),\Omega_{15}(3)).$ For each extension of
$S,$ using \cite{GAP}, we can find two non-singular points of
different type with parameters $(c,d)$ as follow:

if $M=L_4(3)$ then $(c,d)=(42524,20655),$ $ (1160,945);$

if $M=L_4(3)\cdot 2$ then $(c,d)=(311768,154791),$
$(505196,252963).$

Assume that $(S,L)=(U_4(3),\Omega_{15}(3)).$ As in case $L_4(3),$
for each extension $M$ of $U_4(3),$ we can find two non-singular
points of different type with the parameters $(c,d)$ as follow:

If $M=U_4(3)$ or $U_4(3)\cdot 2$ then $(c,d)=(435212,217971),$
$(2780,1755);$

if $M=U_4(3)\cdot 2^2$ then $(c,d)=(217970,108621),$
$(435212,217971);$

if $M=U_4(3)\cdot 4$ or $U_4(3)\cdot D_8$ then
$(c,d)=(217970,108621).$

If  $(S,L)=(L_4(3),\Omega_{19}(3))$ and $M=L_4(3),L_4(3)\cdot 2$ then there exist two non-singular
points of different types with $(c,d)=(2600,1611), (1070,1035).$

Assume $(S,L)=(U_4(3),\Omega_{19}(3)).$ For each extension $M$ of
$U_4(3),$ we can find two non-singular points of different type with
the parameters $(c,d)$ as follow:

if $M=U_4(3), U_4(3)\cdot 2$ then $(c,d)=(2690,1845),
(217700,108891);$

if $M=U_4(3)\cdot 4$ then  $(c,d)=(435752,217431),$
$(217700,108891);$

if $M=U_4(3)\cdot 2^2,U_4(3)\cdot D_8$ then
$(c,d)=(435752,217431),(2420,2115).$

We can check that equation (\ref{eq1}) cannot hold in any of these cases.

{\bf Case $f  \geq 2.$} First consider case $\ell=1.$ As
$SL_2(q)\cong SU_2(q),$ we can assume that $\ep=+.$ If $f=2,$ then
$S=SL_2(9).$ Then $\overline{S}=L_2(9)\cong A_6.$ Thus, we can
assume that $f\geq 3.$ If $\lambda$ is any $3$-restricted dominant
weight then $\lambda=c\lambda_1,$ where $0\leq c\leq 2,$
$dimL(c\lambda_1)=c+1$ and $L(c\lambda_1)$ is self-dual. By
\cite[Theorem $5.4.5$]{kl}, $dimV=(dim\Psi)^f,$ for some irreducible
$k\hat{S}$-module $\Psi.$ As $dimV=2m+1$ is odd, $dim\Psi$ is odd
and hence $dim\Psi\geq 3=dimL(2\lambda_1).$ It follows that
$2m+1\geq 3^f$ and hence $m-1\geq (3^f-3)/2.$ As $|Aut(L_2(3^f))|=f
3^f(3^{2f}-1)< 3^{4f},$ it follows from (\ref{eq4}) that
$(3^f-3)/2\leq 4f,$ with $f\geq 3.$ However by induction on $f\geq
3,$ this is not true. Thus equation (\ref{eq1}) cannot hold.

Consider case $\ell\geq 2.$ It is shown in case $f=1$ that if
$\lambda\in X_3$ such that $L(\lambda)$ is self-dual and has
smallest odd degree then $\lambda=\lambda_1+\lambda_{\ell}.$  By
\cite[Theorem $5.4.5$]{kl} and \cite[$1.11$]{Seitz87} again,
$2m+1=(dim\Psi)^f,$ for some self-dual  irreducible
$k\hat{S}$-module $\Psi$ of odd degree. Hence $dim\Psi\geq
dimL(\lambda_1+\lambda_{\ell}).$ It follows that $2m+1\geq
(\ell^2+2\ell-\ep_3(\ell+1))^f.$ We will show that $3^{m-1}>
|Aut(L^\ep_{\ell+1}(3^f))|.$ Then  (\ref{eq4}) cannot hold. As
$|Aut(L^\ep_{\ell+1}(3^f))|< 3^{f(\ell+1)^2}$ and $m-1\geq
((\ell^2+2\ell-\ep_3(\ell+1))^f-3)/2\geq ( (\ell^2+2\ell-1)^f-3)/2,$
it suffices to show that $ ( (\ell^2+2\ell-1)^f-3)/2>f(\ell+1)^2.$
This is true by induction. The proof is complete.
\end{proof}

Let $\hat{S}$ be a simply connected group of type $B_{\ell}$ over
${\f}_q,$ where $q=p^f,$ and $\mathfrak{G}$ be the corresponding
simply connected, simple algebraic group over $k,$ such that
$S=\mathfrak{G}_{\sigma}$ for some suitable Frobenius map $\sigma.$
Let $N$ be the natural module for $\hat{S}$ with the standard basis
$\beta=\{e_1, \dots, e_{\ell}, x, f_1,\dots, f_{\ell} \}.$
Multiplying some suitable constant to the symmetric bilinear form,
we can assume that  the representing matrix of the symmetric
bilinear form on $N$ has the form $$B=\begin{pmatrix}
    0  & 0& I_{\ell}   \\
      0& 1&0\\
      I_{\ell} & 0&0
\end{pmatrix}.$$

Let $T$ be the set of all matrices of the form $diag(d,1,d^{-1}),$
$d=diag(t_1,\dots,t_{\ell})$ $ \in GL_{\ell}(k).$ As $T\cong
(k^*)^{\ell},$ $T$ is a maximal torus of $\hat{S}.$ For $i=1,\dots,
\ell,$ define $\gamma_i:T\rightarrow k^*,$ by
$\gamma_i(d,1,d^{-1}))=t_i.$ Then $\{\gamma_i\}_{i=1}^{\ell}$ form
an orthonormal basis for $E.$ Also define
$\alpha_{\ell+1-i}=\gamma_{\ell+1-i}-\gamma_{\ell-i},$ for
$i=1,\dots, \ell-1,$ and $\alpha_1=\gamma_1.$ Then  $\{
\alpha_1,\dots, \alpha_{\ell}\}$ is a fundamental root system of
type $B_{\ell} ,$ and the corresponding $\Z$-basis of the
fundamental dominant weights is
$\{\lambda_1,\dots,\lambda_{\ell}\},$ defined as following, for
$i=1\dots, \ell-1$ $\lambda_1=\frac{1}{2}(\gamma_1+
s+\gamma_{\ell}),\mbox{ and
}\lambda_{\ell+1-i}=\gamma_{\ell}+\gamma_{\ell-1}+\dots+\gamma_{\ell+1-i},$

\begin{prop}\label{odd13}
Assume $M$ is almost simple of type $S,$ where $\hat{S}$ is  simply
connected of type $B_{\ell}$ over ${\f}_{3^f}.$ There is an
$M$-orbit on $\m_\xi(V)$ such that equation (\ref{eq1}) does not
hold so that $M$ is not in Tables \ref{in2}-\ref{in4}.
\end{prop}

\begin{proof} Assume  (\ref{eq1}) holds for some $r\in\{s, t\}$ and for any $M$ orbits in $\m_\xi(V).$

{\bf Case $f=1$.} First we claim that if  $\lambda$ is a
$3$-restricted dominant weight such that $dimL(\lambda)$ is odd and
greater than $dimN$ then $\lambda$ must be one of the following
weights:

$(i)$ $\lambda=\lambda_{\ell-1}, \ell\geq 3,\ell$  odd,  and
$dimL(\lambda)=2\ell^2+\ell;$

$(ii)$ $\lambda=2\lambda_{\ell}, \ell=6k+3, 6k+4$ or $6k+5,$ for some non-negative integer $k,$
 and $dimL(\lambda)=2\ell^2+3\ell, 2\ell^2+3\ell-1,2\ell^2+3\ell,$ respectively;

$(iii)$ $\ell=3, \lambda=2\lambda_1,$ and $dimL(\lambda)=35.$

From (\ref{eq4}), we have $3^{m-1}\leq
|Aut(S)|=|Aut(\Omega_{2\ell+1}(3))|=3^{\ell^2}\prod_{i=1}^{\ell}(3^{2i}-1)\leq
3^{2\ell^2+\ell}.$ Hence $dimL(\lambda)=2m+1\leq 4\ell^2+2\ell+3.$
Notice that if $\ell\geq 5$ then  $\ell^3-(4\ell^2+2\ell+3)\geq
5\ell^2-4\ell^2-2\ell-3=(\ell-1)^2-4>0,$ and so
$\ell^3>4\ell^2+2\ell+3.$  If $\ell>11,$ then $dimL(\lambda)\leq
4\ell^2+2\ell+3<\ell^3,$ and hence by \cite[Theorem $5.1$ ]{lub},
$\lambda$ is either $\lambda_{\ell-1}$ or $2\lambda_{\ell}. $ For
$2\leq \ell\leq 11,$ by \cite[Theorem $4.4$]{lub} and the upper
bound for dimension of $L(\lambda)$ above, again, $\lambda$ is one
of the weights above or $\ell=3,$ $\lambda=2\lambda_1,$ and
$dimL(2\lambda_1)=35.$ So case $(iii)$ holds. It remains to get the
restriction on $\ell$ in cases $(i)$ and $(ii).$ From the reference
above, we also have $dimL(\lambda_{\ell-1})=\ell(2\ell+1)$ and
$dimL(2\lambda_{\ell})=2\ell^2+3\ell-\ep_3(2\ell+1).$ Now case $(i)$
holds as $dimL(\lambda_{\ell-1})$ is odd if and only if $\ell$ is
odd. Suppose that $3\mid 2\ell+1.$ then as $2\ell+1$ is odd,
$2\ell+1=3(2t+1),$ hence $\ell=3t+1.$ Since
$dimL(2\lambda_{\ell})=2\ell^2+3\ell-1=\ell(2\ell+3)-1$ is odd, it
follows that $\ell=3t+1$ is even. Thus $t=2k+1$ and $\ell=6k+4.$
With the same argument, we can see that if $3\nmid 2\ell+1,$ then
$\ell=6k+3$ or $6k+4.$

Let $Q$ be the non-degenerate quadratic form associated with the non-degenerate symmetric bilinear
form on $N.$ Then for $v\in N, (v,v)=2Q(v).$ On the tensor product $N\otimes N,$ we can define a
non-degenerate symmetric bilinear form induced from the form on $N$ as follows: for $u_i,w_i\in N,
i=1,2,$ $(u_1\otimes u_2,w_1\otimes w_2)=(u_1,w_1)(u_2,w_2),$ and extend linearly on $N\otimes N.$
Recall  that if $u$ is a non-singular vector in $N,$ then the reflection $r_u:N\rightarrow N$ is
defined by $vr_u=v-\frac{(v,u)}{Q(u)}u,$ for any $v\in N.$  Let $\hat{S}$ act on $N\otimes N$ by
$(u\otimes v)g=(ug\otimes vg).$ Then $((u_1\otimes u_2)r_u,(w_1\otimes w_2)r_u)=((u_1\otimes
u_2,w_1\otimes w_2)$ for any non-singular vector $u.$ Thus, $\wedge^2N$ and $S^2(N)$ leave
invariant  symmetric bilinear forms induced from the one on $N\otimes N.$ We have
$L(\lambda_{\ell-1})\cong \wedge^2(N) \mbox{ and } L(2\lambda_{\ell})\cong w^{\perp}/(w^{\perp}\cap
\langle w\rangle),$ where $w=\sum_{i=1}^{\ell}(e_i\otimes f_i+f_i\otimes e_i)+x\otimes x\in
S^2(N).$

We now consider case $(i).$  As $L(\lambda_{\ell-1})\cong
\wedge^2N,$ $dimL(\lambda_{\ell-1})=\ell(2\ell+1)$ and
$L(\lambda_{\ell-1})$ has a basis consisting of $e_i\wedge e_j,
f_i\wedge f_j, 1\leq i<j\leq \ell,$ $e_i\wedge f_j, 1\leq i, j\leq
\ell$ and $e_i\wedge x, x\wedge f_i, 1\leq i\leq \ell.$ Also denote
by $Q$ the associated quadratic form on $N\otimes N.$ Then for
$\xi\in \{\pm 1\},$ let $v=e_1\wedge x+\xi x\wedge f_1=(e_1-\xi
f_{1})\wedge x.$ Since  $Q(e_1\wedge x)=0=Q(x\wedge f_1),$ we have
$Q(v)=(e_1\wedge x, x\wedge f_1)=\xi.$ Hence $v$ is a non-singular
point. Let $N_1$ be the subspace of $N$ generated by $\{e_1-\xi
f_{1}, x\}.$ As $N_1$ is non-degenerate, $N=N_1\perp N_1^{\perp}.$
Denote by $H$  the centralizer of $N_1$ in $\Omega(N)\cong
\Omega_{2\ell+1}(3).$ It follows that $H\cong \Omega_{2\ell-1}(3),$
and $H$ fixes $v.$  By (\ref{eq4})  we have $3^{m-1}\leq 1+c+d\leq
|Aut(\Omega_{2\ell+1}(3)):\Omega_{2\ell-1}(3)|=2\cdot
3^{2\ell-1}(3^{2\ell}-1)\leq 3^{{4\ell}}.$ Hence $m-1<4\ell,$ so
that $2\ell^2+\ell=2m+1<8\ell+3.$ As $\ell$ is odd and $\ell>1,$ the
above inequality holds only when $\ell=3.$ In this case, we have
$\Omega_7(3)\leq \Omega_{27}(3).$ Using \cite{GAP}, there are two
non-singular points of different type $\la x_i\ra,i=1,2,$ with
$(c_i,d_i)=(13040,9072),(26324,17901),$ we see that equation
(\ref{eq1}) cannot hold in this case.

In case $(ii),$ for $\xi\in\{\pm 1\},$ let $v=e_{1}\otimes e_{1}+\xi
f_1\otimes f_1+\la w\ra \cap w^\perp.$  Then $Q(v)=\xi,$ hence $v$
is  non-singular in $L(2\lambda_{\ell}).$ Let $N_{1}=\langle
e_{1},f_{1}\rangle.$ Then $N_{1}$ is a non-degenerate subspace of
$N.$ As in case $(i),$ let $H$ be the centralizer of $N_{1}$ in
$\Omega(N),$ as $H\cong \Omega_{2\ell-1}(3),$ we have $3^{m-1}\leq
1+c+d\leq
|Aut(\Omega_{2\ell+1}(3)):\Omega_{2\ell-1}(3)|<3^{{4\ell}},$ hence
$2m+1< 8\ell+3.$ As
$dimL(2\lambda_{\ell})=2\ell^{2}+3\ell-\ep_{3}(2\ell+1),$ it follows
that $2\ell^{2}+3\ell-1\leq
2\ell^{2}+3\ell-\ep_{3}(2\ell+1)<8\ell+3.$ If $\ell\geq 4$ then
$2\ell^{2}+3\ell-1\geq 2
4\ell+3\ell-1=8\ell+(3\ell-1)>8\ell+11>8\ell+3,$ and if $\ell=3,$
then $2\ell^{2}+3\ell-\ep_{3}(2\ell+1)=27=8\ell+3.$ Since $\ell\geq
3$ in this case, (\ref{eq4}) cannot hold.

Finally $\ell=3$ and $\lambda=2\lambda_{1}.$ In this case, we have
$\Omega_7(3)\leq \Omega_{27}(3).$ Using \cite{GAP}, there are two
non-singular points of different type  $\la x_i\ra,i=1,2,$  with
$(c_i,d_i)=(13850,8262),(26324,17901),$ we see that equation
(\ref{eq1}) cannot hold in this case.

{\bf Case $f\geq 2.$} By \cite[Theorem $5.4.5$]{kl} and
\cite[$1.11$]{Seitz87}, $dimV=(dim\Psi)^f,$ for some self-dual
irreducible $k\hat{S}$-module $\Psi.$ As $dimV$ is odd, so is
$dim\Psi.$ Firstly, suppose that $f\geq 3.$ Since $dim\Psi$ is at
least $2\ell+1,$ it follows that $2m+1\geq (2\ell+1)^f.$ By
(\ref{eq4}), we have $3^{m-1}\leq
|Aut(\Omega_{2\ell+1}(3^f))|<f\cdot 3^{f(2\ell^2+\ell)}\leq
3^{f(2\ell^2+\ell+1)}.$ As $2m+1\geq (2\ell+1)^f,$ we have
$\frac{1}{2}((2\ell+1)^f-3)<f (2\ell^2+\ell+1).$ Clearing fraction,
we get $(2\ell+1)^f-3-f (4\ell^2+2\ell+2)<0.$ By induction, this
inequality cannot happen.

Secondly, suppose $f=2$ and $dim\Psi>2\ell+1.$ It follows from case
$f=1$ that $dim\Psi\geq 2\ell^2+\ell.$ Arguing as above, we have
$(\ell(2\ell+1))^2-3<2 \cdot 2 (2\ell^2+\ell+1)=2(4\ell^2+2\ell+2).$
As $\ell\geq 2,$ $\ell^f (2\ell+1)^f-3- 2(4\ell^2+2\ell+2)\geq
2^2(2\ell+1)^2-2(2\ell+1)^2+4\ell-5\geq 2(2\ell+1)^2+3>0.$ Hence
$\ell^f (2\ell+1)^f-3>2(4\ell+2\ell+2),$ a contradiction.

Finally, suppose $f=2$ and $dim\Psi=2\ell+1.$ In this case, we can
assume that $\Psi\cong L(\lambda_{\ell})\cong N,$  and so
$V=N\otimes N^{(1)},$ where $N$ is the natural module for
$\Omega_{2\ell+1}(3^2),$ and $N^{(1)}$ denote the module received
from the twist action of $\hat{S}$ on $N.$ For any element $v\in N,$
denote by $v^{(1)}$ the corresponding element in $N^{(1)}.$ Notice
that if $p$ is odd then for any $a,b\in {\f}_p^f,$  $a+b=0$ or $1$
if and only if  $a^p+b^p=0$ or $1,$ correspondingly. This holds
because $a^p+b^p=(a+b)^p.$ Fix a standard basis $\beta=\{e_1, \dots,
e_{\ell},x,f_{\ell}, \dots, f_1\} $ of $S.$ Let $u=(\ep_1+\xi
f_1)\in N.$ Then for $g\in S,$ in the basis $\beta,$ we write
$g=(a_{i,j}).$ Assume that $ug=g.$ Then
$a_{11}+a_{2\ell+1,1}=1=a_{1,2\ell+1}+a_{2\ell+1,2\ell+1}$ and
$a_{1i}+a_{2\ell+1,i}=0, 1<i<2\ell+1.$ Hence, by the notice above,
we have
$a_{11}^p+a_{2\ell+1,1}^p=1=a_{1,2\ell+1}^p+a_{2\ell+1,2\ell+1}^p$
and $a_{1i}^p+a_{2\ell+1,i}^p=0, 1<i<2\ell+1.$ Therefore,
$u^{(1)}g=ug^{\nu}=u(a_{ij}^p)=u=u^{(1)}.$ This means that if $g\in
\hat{S}$ fixes $u$ then $g$ also fixes $u^{(1)}.$ Let $v=u\otimes
u^{(1)}\in N\otimes N^{(1)}.$ Let $H$ be the stabilizer of $u$ in
$\hat{S}.$ Then $H\cong \Omega^{\ep}_{2\ell}(3^2)$ with $\ep=\pm 1,$
and $H\leq \hat{S}_v.$ Hence by (\ref{eq4}), $3^{m-1}\leq
|Aut(\Omega_{2\ell+1}(3^2)):\Omega^{\ep}_{2\ell}(3^2)|\leq
|Aut(\Omega_{2\ell+1}(3^2)):\Omega^{+}_{2\ell}(3^2)|\leq
3^{4\ell+2}.$ Since $2m+1=(2\ell+1)^2,$ it follows that
$(2\ell+1)^2<8\ell+7.$ As $\ell\geq 2,$
$(2\ell+1)^2-8\ell-7=4\ell^2+4\ell+1-8\ell-7=4\ell(\ell-1)-6>0.$
This final contradiction finishes the proof.
 \end{proof}
Let $\hat{S}$ be a simply connected group of type $C_{\ell}$ over ${\f}_q,$ where $q=p^f,$ and $\mathfrak{G}$ be the corresponding simply connected, simple algebraic group over $k,$ such that $\hat{S}=\mathfrak{G}_{\sigma}$ for some suitable Frobenius map $\sigma.$ Let $N$ be the natural module for $\hat{S}$ with the standard basis $\beta=\{e_1, \dots, e_{\ell}, f_1,\dots, f_{\ell} \}.$ The representing matrix of the non-degenerate symplectic form on $N$ has the form $$B=\begin{pmatrix}
    0  &  I_{\ell}   \\
      -I_{\ell} &0
\end{pmatrix}.$$
From the isomorphisms $Sp_2(q)\cong SL_2(q),$ and $ Sp_4(q)\cong \Omega_5(q)$ for $q$ odd, we can assume that $\ell\geq 3.$

\begin{prop}\label{odd14}
Assume $M$ is almost simple of type $S,$ where $\hat{S}$ is  simply connected of type $C_{\ell}$ over ${\f}_{3^f}, $ with $\ell\geq 3.$
There is an $M$-orbit on $\m_\xi(V)$ such that equation (\ref{eq1}) does not hold unless $(\ell,\lam,dimV)=(3,\lam_2,13)$ or $(L,S,\lam)= (\Omega_{41}(3),PSp_8(3),\lam_1).$ If the first case holds then  $M$ has at most two orbits on $\m_\xi(V)$ so that $1_P^G\not\leq 1_M^G$ by Corollary \ref{smallorbits} and hence $(L,S)=(\Omega_{13}(3),PSp_6(3))$ is in Table \ref{in3}. The last case is in Table \ref{in4}.
\end{prop}

\begin{proof} Assume (\ref{eq1}) holds for some $r\in\{s, t\}$ and for any $M$ orbits in $\m_\xi(V).$

{\bf Case $f=1.$} Let $\lambda\in X_3$ be a $3$-restricted dominant
weight such that $L(\lambda)\cong V.$ We first get an upper bound
for $dimV.$ As $|Aut(PSp_{2\ell}(3))|=f
3^{\ell^2}\prod_{i=1}^{\ell}(3^{2i}-1)\leq 3^{2\ell^2+\ell},$ by
(\ref{eq4}), $3^{m-1}\leq 3^{2\ell^2+\ell},$ and hence $2m+1\leq
4\ell^2+2\ell+3.$ If $\ell\geq 5$ then $ 4\ell^2+2\ell+3<\ell^3, $
and so $dimV<\ell^3.$  If $\ell>11,$ then $dimL(\lambda)<\ell^3,$
and hence by \cite[Theorem $5.1$]{lub}, $\lambda$ is either
$\lambda_{\ell-1}$ or $2\lambda_{\ell}. $ For $2\leq \ell\leq 11,$
by \cite[Theorem $4.4$ ]{lub} and the upper bound for dimension of
$L(\lambda)$ above, again, $\lambda$ is one of the weights above or
$\ell=4,$ $\lambda=\lambda_1,$ and $dimL(\lambda_1)=41.$ We have
$dimL(2\lambda_{\ell})=\ell(2\ell+1), L(2\lambda_{\ell})\cong
S^2(N),$ and $dimL(\lambda_{\ell-1})=2\ell^2-\ell-1-\ep_p(\ell),
L(\lambda_{\ell-1})\cong w^{\perp}/(\langle w\rangle \cap
w^{\perp}),$ where $w=e_1\wedge f_1+ s+e_{\ell}\wedge f_{\ell}.$ In
these cases, $\hat{S}$ leaves invariant a quadratic form $Q$ induced
from the symplectic form on $N.$ In case $\lambda=2\lambda_{\ell},$
let $v=e_1\otimes e_1+\xi  f_1\otimes f_1\in L(2\lambda_{\ell}).$
Since $dimL(2\lambda_{\ell})=\ell(2\ell+1)$ is odd, $\ell$ must be
odd.

If $\ell=3$ then $dimL(2\lambda_{\ell})=21$ and $PSp_6(3)\leq \Omega_{21}(3).$ Using \cite{GAP}, there are two non-singular points of different type  $\la x_i\ra,i=1,2,$  with $(c_i,d_i)=(7075430,3538809),$ $(26324,17901),$ we see that equation (\ref{eq1}) cannot hold in this case.

Thus we assume that $\ell\geq 5.$ Let $H$ be the centralizer in
$\hat{S}$ of the subspace generated by $\{e_1, f_1\}.$ Then $H\cong
Sp_{2\ell-2}(3).$ By  (\ref{eq4}), we have $3^{m-1}\leq
|Aut(PSp_{2\ell}(3)):PSp_{2\ell-2}(3)|=2\cdot
3^{\ell^2}(3^{2\ell}-1)/3^{(\ell-1)^2}<3^{4\ell}.$ Hence
$2m+1<8\ell+3.$ As $2m+1=\ell(2\ell+1),$ it follows that
$\ell(2\ell+1)<8\ell+3.$ However as $\ell\geq 5,$ $\ell(2\ell+1)\geq
5(2\ell+1)=10\ell+5>8\ell+3,$ a contradiction. Next consider case
$\lambda=\lambda_{\ell-1}.$ As $dimL(\lambda_{\ell-1})=
2\ell^2-\ell-1-\ep_p(\ell)$ is odd, $\ell=6k+2, 6k+3$ or $6k+4.$

Let $v=e_1\wedge e_2+\xi  f_1\wedge f_2+\la w\ra\cap w^\perp,$ where
$\xi=\pm 1.$ Then $v$ is non-singular in $V.$ Let $N_1=\langle
e_1,e_2,f_1,f_2\rangle$ be a subspace of $N.$ Since $N_1$ is
non-degenerate, $N=N_1\perp N_1^{\perp}.$ Let $H, K$ be the
centralizers in $\hat{S}$ of $N_1,N_1^{\perp},$ respectively. Then
$H\cong Sp_{2(\ell-2)}(3),$ $K\cong Sp_4(3),$ and $H,K$ commute. In
the basis $\beta_1=\{e_1,e_2,f_1,f_2\},$ let $$g=\begin{pmatrix}
     1&0 & 1&0   \\
     0&1&0&0\\
     0&0&1&0\\
      0&-\xi&0&1
\end{pmatrix}, h=\begin{pmatrix}
      1&1&0&0\\
      0&1&0&0    \\
      0&0&1&0\\
      0&0&-1&1
\end{pmatrix}.$$ As $gBg^t=B, hBh^t=B,$ and $det(g)=1=det(h),$ $g,h\in Sp(V_1).$ Furthermore, $g,h$ are of order $3$ and $gh=hg,$ the subgroup generated by $g$ and $h$ are elementary abelian of order $9.$ Since $vg=v$ and $vh=h,$ it follows that $E=\langle g, h\rangle \leq K_v.$
Thus $E\times H\leq \hat{S}_v,$ hence $1+c+d\leq |Aut(S):(E\times H)|=|Aut(PSp_{2\ell}(3)):(E\times PSp_{2(\ell-2)}(3))|<3^{8\ell-7}.$ Hence $2m+1<16\ell-11.$ Since $2m+1=2\ell^2-\ell-1-\ep_p(\ell)\geq 2\ell^2-\ell-2,$ we have $ 2\ell^2-\ell-2<16\ell-11,$ or equivalent $2\ell^2-17\ell+9<0.$  As $\ell=6k+2, 6k+3,6k+4,$ if $\ell>4$ then $k\geq 1,$ and so $\ell\geq 8.$  Then  $2\ell^2-17\ell+9\geq 2\ell^2-16\ell-(\ell-8)+1=(2\ell-1)(\ell-8)+1>0,$ a contradiction. Thus $\ell\leq 4.$

If $\ell=4,$ then $dimL(\lambda_{\ell-1})=27,$ by using \cite{GAP}, equation (\ref{eq1}) cannot hold. If $\ell=3,$ then $dimL(\lambda_{\ell-1})=13.$ Using \cite{GAP}, there is only one orbit of minus points and two orbits of plus points. Hence equation (\ref{eq1}) holds for both types of points by Corollary \ref{smallorbits}. We are left with case $\ell=4,$  $\lambda=\lambda_1, dimL(\lambda_1)=41.$

{\bf Case $f\geq 2.$} It follows from case $f=1$ that if $\lambda$ is a $3$-restricted dominant weight such that $L(\lambda)$ is self-dual and is of odd degree then $dimL(\lambda)\geq 2\ell^2-\ell-1-\ep_p(\ell)\geq 2\ell^2-\ell-2.$ By \cite[Theorem $5.4.5$]{kl} and \cite[$1.11$]{Seitz87}, $2m+1=(dim\Psi)^f,$ for some self-dual  irreducible $k\hat{S}$-module $\Psi$ of odd degree. Thus $2m+1=(dim\Psi)^f\geq (2\ell^2-\ell-2)^f.$ By (\ref{eq4}), $3^{m-1}\leq |Aut(PSp_{2\ell}(3))|\leq 3^{f(2\ell^2+\ell+1)},$ so that $2m+1<2f(2\ell^2+\ell+1)+3.$ Combining these inequalities, we have $(2\ell^2-\ell-2)^f<2f(2\ell^2+\ell+1)+3,$ where $\ell\geq 3, f\geq 2.$
However by induction this inequality cannot happen. \end{proof}

Let $\hat{S}$ be a simply connected group of type $D_{\ell}$ or ${}^2 D_\ell$ over ${\f}_q,$ where $q=p^f,$ and $\mathfrak{G}$ be the corresponding simply connected, simple algebraic group over $k,$ such that $S=\mathfrak{G}_{\sigma}$ for some suitable Frobenius map $\sigma.$ Let $N$ be the natural module for $\hat{S}$ with the standard basis $\beta.$ The representing matrix of the non-degenerate symmetric bilinear form on $N$ has the form $$B=\begin{pmatrix}
    0  &  I_{\ell}   \\
      I_{\ell} &0
\end{pmatrix}.$$
From the isomorphisms $\Omega_4^+(q)\cong SL_2(q)\circ SL_2(q),  \Omega_4^-(q)\cong L_2(q)$ and $ P\Omega_6^{\pm}(q)\cong L_4^{\pm}(q),$ we can assume that $\ell\geq 4.$

\begin{prop}\label{odd15}
Assume $M$ is almost simple of type $S,$ where $\hat{S}$ is  simply
connected of type $D_{\ell}$ or $^2D_{\ell}$ over ${\f}_{3^f}.$
There is an $M$-orbit on $\m_\xi(V)$ such that equation (\ref{eq1})
does not hold so that $M$ is not in Tables \ref{in2}-\ref{in4}.
\end{prop}
\begin{proof} Assume  (\ref{eq1}) holds for some $r\in\{s, t\}$ and for any $M$ orbits in $\m_\xi(V).$

{\bf Case $f=1.$} Let $\lambda\in X_3$ be a $3$-restricted dominant
weight such that $L(\lambda)\cong V.$ By inequality (\ref{eq4}),
$3^{m-1}\leq |Aut(P\Omega_{2\ell}^{\ep}(3))|\leq
3^{2\ell^2-\ell+2}.$ Thus $2m+1\leq 4\ell^2-2\ell+7.$ By
\cite[Theorem $5.1$ and $4.4$]{lub},  either $\ell\geq 4$ and
$\lam=\lambda_{\ell-1},2\lambda_{\ell}$ or $\ell=4$ and
$\lam=2\lam_1,2\lam_2.$

If $\lambda=\lambda_{\ell-1},$ then $dimL(\lambda)=2\ell^2-\ell$ and
$L(\lambda)\cong \wedge^2N.$ Since $dimL(\lambda)$ is odd, $\ell$
must be odd. Thus $\ell\geq 5.$ Let $a=e_1-\xi  f_1, b=e_2+f_2\in N$
and $z=a\wedge b\in \wedge^2N,$  where $\xi=\pm 1.$ Then $\la z\ra$
is non-singular in $V=\wedge^2N.$ Also, let $N_1$ be a subspace of
$N$ generated by $a$ and $b.$  Clearly, $N_1$ is non-degenerate,
$dim(N_1)=2,$ and $sgn(N_1)=\xi.$ Since $N=N_1\perp N_1^{\perp},
dim(N_1^{\perp})=2\ell-2$  and $sgn(N)=sgn(N_1)\cdot
sgn(N_1^{\perp}),$ it follows from \cite[ Proposition $2.5.11$]{kl}
that $sgnN_1^{\perp}=\ep \xi.$  Since  the discriminant $D(Q)\equiv
detB=(-1)^{\ell}=-1 \;( \mbox{mod} ({\f}^*)^2),$ as $\ell$ is odd,
by \cite[Proposition $4.1.6$]{kl}, $H=\Omega(N_1^{\perp})\cong
\Omega_{2\ell-2}^{\ep \xi}(3)\leq \Omega(N)$ centralizes $N_1.$
Hence $H\leq M_{\la z\ra}.$ Therefore $1+c+d\leq
|Aut(P\Omega_{2\ell}^{\ep}(3)):\Omega_{2\ell-2}^{\ep \xi}(3)|=4\cdot
3^{2\ell-2}(3^\ell-\ep )(3^{\ell-1}+\ep \xi )\leq 3^{4\ell-1}.$
Since $\ell\geq 5,
3^{m-1}=3^{(2\ell^2-\ell-3)/2}>3^{4\ell-1}>1+c+d,$ a contradiction
to (\ref{eq4}).

If $\lambda=2\lambda_{\ell},$ then
$dimL(\lambda)=2\ell^2+\ell-1-\ep_3(\ell)$ and $V=L(\lambda)\cong
w^{\perp}/(w^{\perp}\cap \langle w\rangle),$ where
$w=\sum_{i=1}^{\ell}(e_i\otimes f_i+f_i\otimes e_i)\in S^2N$ if
$\hat{S}$ is of type $D_\ell$ and $w=\sum_{i=1}^{\ell-1}(e_i\otimes
f_i+f_i\otimes e_i)+x\otimes x+y\otimes y$ otherwise. Let
$z_\xi=e_1\otimes e_1+\xi  f_1\otimes f_1+\la w\ra \cap w^\perp,
\xi=\pm 1,$ and $N_1=\langle e_1, f_1\rangle\leq N.$ Then $z_\xi$ is
non-singular in $V.$ Let $N_1=\langle e_1, f_1\rangle\leq N$ and $H$
be the centralizer of $N_1$ in $\Omega(N)\cong
\Omega_{2\ell}^\ep(3).$ Since $sgn(N_1^\perp)=\ep,$  $H\cong
\Omega_{2\ell-2}^\ep(3)$ and $H$ fixes $\la z_\xi\ra.$ Thus
$1+c+d\leq |Aut(P\Omega_{2\ell}^\ep(3)):\Omega_{2\ell-2}^\ep(3)|\leq
12\cdot 3^{2\ell-2}(3^{\ell}-\ep )(3^{\ell-1}+\ep )\leq 3^{4\ell}.$
If $\ell\geq 5,$ then $3^{m-1}\geq
3^{(2\ell^2+\ell-5)/2}>3^{4\ell}>1+c+d,$ contradicts to (\ref{eq4}).
If $\ell=4,$ then $2m+1=35,$ hence $m-1=16.$ As
$3^{m-1}=3^{16}=3^{4\ell}>1+c+d,$ we also get a contradiction to
(\ref{eq4}).

{\bf Case $f\geq 2.$} It follows from case $f=1$ that if $\lambda$
is a $3$-restricted dominant weight such that $L(\lambda)$ is
self-dual and is of odd degree then $dimL(\lam)\geq 2\ell^2-\ell.$
As $V$ is an absolutely irreducible $k\hat{S}$-module, by
\cite[Theorem $5.4.5$]{kl} and \cite[$1.11$]{Seitz87},
$2m+1=(dim\Psi)^f,$ for some self-dual  irreducible
$k\hat{S}$-module $\Psi$ of odd degree. Thus $2m+1=(dim\Psi)^f\geq
(2\ell^2-\ell)^f.$ By  (\ref{eq4}), $3^{m-1}\leq
|Aut(P\Omega_{2\ell}^\ep(3))|\leq 3\cdot 3^{f(2\ell^2-\ell+1)},$ so
that $2m+1<2f(2\ell^2-\ell+1)+5.$ Combining these two inequalities,
we have $(2\ell^2-\ell)^f<2f(2\ell^2-\ell+1)+5,$ where $\ell\geq 4,
f\geq 2.$ By induction, this cannot happen.
 \end{proof}

Assume that $\hat{S}$ is simply connected of exceptional type defined over a field of characteristic $3.$ Then $\hat{S}$
is one of the following types: $G_2,F_4,E_6,E_7,$ $E_8,{}^2E_6,{}^3D_4,{}^2G_2.$

\begin{prop}\label{odd16}
Assume $M$ is almost simple of type $S,$ where   $\hat{S}$ is of
exceptional type above and is defined over ${\f}_{3^e}$ with $e\geq 1.$
There is an $M$-orbit on $\m_\xi(V)$ such that equation (\ref{eq1})
does not hold unless $(L,S,\lam)=(\Omega_7(3),G_2(3),\lam_i),$ $
i=1,2,$ or $(F_4(3),\Omega_{25}(3),\lam_4),$ and $M$  has one or at
most two orbits on $\m(V),$ respectively, so that $1_P^G\not\leq
1_M^G$ by Corollary \ref{smallorbits} and so  they are in Table
\ref{in3} or $(L,S,\lam)=(\Omega_{77}(3),E_6(3),\lam_2),$
$(\Omega_{133},E_7(3),\lam_1)$ and they are in Table \ref{in4}.
\end{prop}
\begin{proof}Assume  (\ref{eq1}) holds for some $r\in\{s, t\}$ and for any $M$ orbits in $\m_\xi(V).$

$(a)$ {\bf Case $G_2.$} By (\ref{eq4}) $3^{m-1}\leq
|Aut(G_2(3^e))|\leq 3^{15e}.$ Thus $2m+1\leq 30e+3.$ Assume first
that $e=1.$ Then $2m+1\leq 33.$ Let $\lam\in X_3$ be a
$3$-restricted dominant weight with $V\cong L(\lam).$ From Appendix
$A.49$ in \cite{lub}, $\lam$ is one of the following weights
$\lam_1,\lam_2,2\lam_1,2\lam_2.$ If $\lam$ is $\lam_i, i=1,2,$ then
$dimV=7.$ In these cases, we have $G_2(3)\leq \Omega_7(3).$ By
\cite{atlas}, these are maximal embeddings and $G_2(3)$ has only one
orbit in $\m_\xi(V).$ If $\lam=2\lam_1$ or $2\lam_2$ then $dimV=27.$
We have $G_2(3)\leq \Omega_{27}(3).$ But this is not a maximal
embedding as $G_2(3)\leq \Omega_7(3)\leq \Omega_{27}(3),$ where the
last embedding arises from the symmetric square of the natural
module for $\Omega_7(3)$ (see Proposition \ref{odd13}). Assume that
$e\geq 2.$ By \cite[Theorem $5.4.5$]{kl} and \cite[$1.11$]{Seitz87},
$2m+1=(dim\Psi)^e,$ for some self-dual irreducible $k\hat{S}$-module
$\Psi$ of odd degree. It follows that $2m+1\geq 7^e.$ Combining with
$2m+1\leq 30e+3,$ we have $e=2.$ Then $dimV=7^2=49.$ However
$G_2(9)$ is not maximal in $\Omega_{49}(3)$ since $G_2(9)\leq
\Omega_7(9)\leq \Omega_{7^2}(3)$ where the first embedding arises as
in previous case, while the second is the twisted tensor product
embedding.

$(b)$ {\bf Case $F_4.$} By  (\ref{eq4}), $3^{m-1}\leq |Aut(F_4(3^e))|\leq 3^{53e}.$ So $2m+1\leq 106e+3.$
  Assume that $e=1.$ From Appendix $A.50$ in \cite{lub}, $\lam=\lam_4$ and $dimL(\lam)=25.$ In this case, we have an embedding $F_4(3)\leq \Omega_{25}(3).$ By \cite{coco}, $F_4(3)$ has $5$ orbits of points in $V.$ But there are two orbits of singular points, and so there are at most two orbits for each types of non-singular points. Assume that $e\geq 2.$ we have $2m+1\geq 25^e,$ so that $25^e\leq 106e+3.$ But this cannot happen for any $e\geq 2.$

$(c)$ {\bf Case ${}^\ep E_6.$} By  (\ref{eq4}), $3^{m-1}\leq |Aut({}^\ep E_6(3^e))|\leq 3^{79e}.$ Thus $2m+1\leq 158e+3.$
  Assume that $e=1.$ From Appendix $A.51$ in \cite{lub}, $\lam=\lam_2$ and $dimL(\lam)=77.$ (Note that $dimL(\lam_1)=dimL(\lam_6)=27$ but these modules are not self-dual). In this case, we have  ${}^2 E_6(3)\leq E_6(3)\leq \Omega_{77}(3),$ and $V=L(E_6)/Z(L(E_6)),$ where $L(E_6)$ is the Lie algebra of $E_6$ over ${\f}_3.$ If $e\geq 2,$ then $77^e\leq 158e+3.$ However this cannot happen for any $e\geq 2.$.

$(d)$ {\bf Case $E_7.$} By (\ref{eq4}), $3^{m-1}\leq |Aut(E_7(3^e))|\leq 3^{134e}.$ Thus $2m+1\leq 268e+3.$
  Assume that $e=1.$ From Appendix $A.52$ in \cite{lub}, $\lam=\lam_1$ and $dimL(\lam)=133.$ We have  $E_7(3)\leq \Omega_{133}(3),$ and $V=L(E_7),$ the Lie algebra of $E_7$ over ${\f}_3.$
Assume that $e\geq 2.$ We have  $133^e\leq 268e+3.$ This cannot happen for $e\geq 2.$

$(e)$ {\bf Case $E_8.$} By  (\ref{eq4}), $3^{m-1}\leq |Aut(E_8(3^e))|\leq 3^{249e}.$ Thus $2m+1\leq 498e+3.$
Assume that $e=1.$ From Appendix $A.53$ in \cite{lub},  $dimL(\lam)\geq 3875>498.1+3=501.$
Assume that $e\geq 2.$ Clearly  $3875^e> 498e+3$ for any $e\geq 2.$

$(f)$ {\bf Case ${}^3D_4.$} By  (\ref{eq4}), $3^{m-1}\leq
|Aut({}^3D_4(3^e))|\leq 3^{30e}.$ Thus $2m+1\leq 60e+3.$ Assume that
$e=1.$ From Appendix $A.53$ in \cite{lub},
$\lam=2\lam_1,2\lam_2,2\lam_4$ and $dimL(\lam)=35.$ Since the
splitting field for  ${}^3D_4(3)$ is ${\f}_{3^3},$ $dimV=3\cdot
35=105>63.$ Assume that $e\geq 2.$ Clearly  $35^e > 60e+3$ for any
$e\geq 2.$

$(g)$ {\bf Case ${}^2G_2.$} By  (\ref{eq4}), $3^{m-1}\leq |Aut({}^2G_2(3^{2e+1}))|\leq 3^{8(2e+1)}.$ Thus $2m+1\leq 32e+19.$ By \cite[Theorem $5.4.5$]{kl}, $2m+1\geq 7^{2e+1},$ and so $7^{2e+1}\leq 32e+19.$ This cannot happen.
 \end{proof}

\subsubsection{Embedding of Sporadic groups}
Let ${S}$ be a sporadic simple group and let $\hat{S}$ be the
universal covering group of $S.$ Define
$g_3(S):=[2log_3(|Aut(S)|)+4]$ and $\Re_3(S)$ to be the minimal
degrees of irreducible faithful representations
 of $S$ and its covering groups over ${\f}_3.$

 \begin{lem}\label{spor1} Let $S$ be a simple sporadic group. Then $\Re_3(S)$ and $g_3(S)$ are given in Table \ref{spor}.
  \begin{table}
\centering

 \caption{Minimal degrees of representations for sporadic groups in characteristic $3$.}\label{spor}
 \begin{tabular}{l|l|l|l|l|r}\hline
$S$ &$\Re_3(S)$&$g_3(S)$&$S$ &$\Re_3(S)$&$g_3(S)$\\ \hline
$M_{11}$&$5$&$20$&$Suz$&$64$&$54$\\
$M_{12}$&$10$&$26$&$2.Suz$&$12$&$54$\\
$2.M_{12}$&$6$&$26$&$O'N$&$154$&$54$\\
$J_1$&$56$&$25$&$Co3$&$22$&$53$\\
$M_{22}$&$21$&$28$&$Co2$&$23$&$61$\\
$2.M_{22}$&$10$&$28$&$Fi_{22}$&$77$&$63$\\
$J_{2}$&$13$&$29$&$2.Fi_{22}$&$176$&$63$\\
$2.J_{2}$&$6$&$29$&$HN$&$133$&$65$\\
$M_{23}$&$22$&$33$&$Ly$&$651$&$74$\\
$HS$&$22$&$37$&$Th$&$248$&$75$\\
$2.HS$&$56$&$37$&$Fi_{23}$&$253$&$82$\\
$J_3$&$18$&$37$&$Co1$&$276$&$82$\\
$M_{24}$&$22$&$39$&$2.Co1$&$24$&$82$\\
$McL$&$21$&$42$&$J_4$&$1333$&$87$\\
$He$&$51$&$45$&$Fi_{24}'$&$781$&$106$\\
$Ru$&$378$&$50$&$B$&$4371$&$144$\\
$2.Ru$&$28$&$50$&$2.B$&$96256$&$144$\\
$M$&$196882$&$229$&&&\\
\hline
\end{tabular}
\end{table}
 \end{lem}

\begin{prop}\label{odd17}
Assume $M$ is almost simple of type $S,$ where $S$ is a simple sporadic group.  There is an $M$-orbit on $\m_\xi(V)$ such that equation (\ref{eq1})
does not hold so that $M$ is not in Tables \ref{in2}-\ref{in4}.
\end{prop}
\begin{proof} Assume  equation (\ref{eq1}) holds for some $r\in\{s, t\}$ and for any $M$ orbits in $\m_\xi(V).$ Then by (\ref{eq4}),
$|Aut(S)|\geq 3^{m-1}.$ It follows that $\Re_3(S)\leq 2m+1\leq
2log_3(|Aut(S)|)+3\leq g_3(S).$ Recall that $V$ is an absolutely
irreducible ${\f}_3\hat{S}$-module with $dimV=2m+1$ and the
Frobenius-Schur indicator of $V$ is $+.$ By \cite{hm} and Lemma
\ref{spor1}, we need to consider the following cases:
$(S,dimV)=(M_{22},21),$ $(McL,21),$ and $(Co_2,23).$
  If $(S,dimV)=(M_{22},21),$ then $M_{22}<A_{22}<\Omega_{21}(3),$ hence $N_{G}(S)$ is not maximal in $G.$
 If $(S,dimV)=(McL,21),$ then $S$ has two orbits with representatives $\la x\ra,\la y\ra$ and
 stabilizers $S_{\la x\ra}=L_3(4):2_2$ and $S_{\la y\ra}=M_{11}$ which are in different $G$-orbits.
 We also have $c_x=12194,d_x=10080$ and $c_y=72809,d_y=40590,$ and we can check that  (\ref{eq1}) cannot hold
 in any of these cases. For $S:2,$ we also get the same result. Finally if $(S,dimV)=(Co_2,23),$ then there exist
 two non-singular points in different $G$-orbits with representative $\la x\ra ,\la y\ra$ with $S_{\la x\ra}=2^{10}:M_{22}:2$
 and $S_{\la y\ra}=HS:2$ with sizes $|\la x\ra S|=46575<3^{m-1}=3^{10}, |\la y\ra S|=476928.$ The parameters for $\la y\ra S$
 are $c_y=296450,d_y=180477.$ We can check that (\ref{eq1}) cannot hold in this case.
 \end{proof}
\subsubsection{Computation with GAP} In this section we briefly
demonstrate how to compute the parameters $c$ and $d$  defined in Section \ref{higman3} by using
GAP \cite{GAP}. We still assume the hypotheses and notation from Section \ref{hypo}. Recall from
section \ref{charcon} that if $\la x\ra\in \m_{\xi}(V)$ then
$$\Delta(x)=\m_\xi(V)\cap x^\perp \mbox{ and
}\Gamma(x)=\m_\xi(V)\cap(V-x^\perp-\{\la x\ra\}).$$ Let $M$ be a
subgroup of $G,$  where $G$ is almost simple with socle
$L=\Omega_{2m+1}(3),m\geq 3.$ Note that $\m_{\xi}(V),\xi=\pm,$ is an
$L$-orbit of non-singular points of type $\xi.$ The parameters $c$
and $d$ are defined as follows: $$d=|\la x\ra M\cap \Delta(x)|\mbox{
and }c=|\la x\ra M\cap \Gamma(x)|=|\la x\ra M|-d-1.$$ As an example,
let $M\cong S_8,$ a symmetric group of degree $8$ and
$L=\Omega_{13}(3).$ Then $M$ embeds into $L$ via the Specht module
$V$ corresponding to the partition $(6,2)$ (see Proposition
\ref{odd10}). Let $U$ be the permutation module for $S_8$ in
characteristic $3.$ We know that $V$ is a composition factor of the
tensor product $U\otimes U.$ Using MeatAxe to decompose $U\otimes U$
to obtain $V.$ From this module, we can obtain the embedding of $M$
into $\Omega(V)$ and also the quadratic form of $V.$ Using this
information, it is very easy to compute $c$ and $d.$

{\em $>$U:=PermutationGModule(SymmetricGroup($8$),$GF(3)$);;

$\sharp$ $GF(3)$ is the field of size $3.$

$>$T:=TensorProductGModule(U,U);

$\sharp$ The tensor product  $U\otimes U.$

$>$decomp:=MTX.CompositionFactors(T);

$>$V:=decomp[Position(List(decomp,x$->$MTX.Dimension(x)),$13$)];

$\sharp$ find the module of dimension $13.$

$>$f:=MTX.InvariantBilinearForm(V);

$\sharp$ The bilinear form on $V.$

$>$gens:=MTX.Generators(V);

$\sharp$ The generators for $M$ when embedded in $\Omega(V).$

$>$M:=Group(gens);;

$>$ob:=OrbitsDomain(M,GF(3) $\hat{}$(MTX.Dimension(V)));;

$\sharp$ all orbits of $M$ on $V.$

$>$leng:=List(ob,x$->$Size(x));;

$\sharp$ The lengths of all orbits of $M$ on $V$;

$>$Positions(leng, $315$);

$\sharp$ find the positions of orbits of length $315.$

$\sharp$ for example, the positions are $[54,65,72,73].$ These
positions are not fixed.

$>$ $x:=ob[54][1];;$ $x*f*x;$

$\sharp$ $x$ is a representative for the orbit $ob[54].$

$\sharp$ $x*f*x;$ is the norm of $x.$

$\sharp$ The different norms mean the points are in different
$L$-orbits.

$>$ op:=Orbit(G,GF(3) $\hat{}$(MTX.Dimension(V)),x);;

$\sharp$ op is the $M$-orbit $\la x\ra M.$

$>$ d:=0;;for i in [1..Size(op)] do if $op[i]*f*x=0*Z(3)$ then
$d:=d+1;$ fi;od;

$\sharp$ $Z(3)$ is a generator for the field $GF(3).$

$>$ c:=Size(op)-1-d;}

For matrix groups, given the matrix generators `gens', which can be
obtained from the package `atlasrep' or from the Atlas website
\cite{Atlasweb}. The module is then constructed by using the GAP
command `V:=GModuleByMats(gens,GF(3))'.

{\em $>$ LoadPackage(`atlasrep');

$>$ DisplayAtlasInfo(`$G_2(4)$');

$>$ gens:=AtlasGenerators(`$G_2(4)$',9).generators;;

$\sharp$ This is the $64$-dimensional representation of $G_2(4)$ in
characteristic $3.$

$>$ V:=GModuleByMats(gens);;

$>$ M:=Group(gens);}

\subsection*{Acknowledgment} This work is part of my Ph.D thesis. I
would like to thank my thesis advisor Professor Kay Magaard for his
guidance and his help.

\bibliographystyle{amsplain}

\end{document}